\PassOptionsToPackage{dvipsnames}{xcolor}
\documentclass[11pt,reqno]{amsart}
\usepackage[letterpaper,hmargin=1.15in,vmargin=1.1in]{geometry}
\usepackage[T1]{fontenc}
\usepackage{lmodern,microtype}
\usepackage{amssymb,mathtools,mathrsfs}
\usepackage{needspace,etoolbox}
\usepackage{placeins}
\usepackage{enumitem}
\usepackage{tikz}
\usetikzlibrary{arrows.meta,positioning,calc}
\usepackage[dvipsnames]{xcolor}
\usepackage[colorlinks=true,linkcolor=MidnightBlue,citecolor=OliveGreen,urlcolor=MidnightBlue]{hyperref}
\hypersetup{pdftitle={A rank-zero 2-converse to the Gross-Zagier and Kolyvagin theorem for elliptic curves with abelian two-division field},pdfauthor={Ashay Burungale and Ye Tian}}

\setlist[enumerate]{label=\textup{(\arabic*)},leftmargin=*,itemsep=2pt,topsep=4pt}
\newcommand{\Q}{\mathbb Q}
\newcommand{\Z}{\mathbb Z}
\newcommand{\F}{\mathbb F}
\newcommand{\RR}{\mathbb R}
\newcommand{\C}{\mathbb C}
\DeclareFontEncoding{OT2}{}{}\DeclareFontSubstitution{OT2}{wncyr}{m}{n}
\DeclareFontFamily{OT2}{wncyr}{}\DeclareFontShape{OT2}{wncyr}{m}{n}{<->wncyr10}{}
\DeclareSymbolFont{cyrletters}{OT2}{wncyr}{m}{n}
\DeclareMathSymbol{\Sha}{\mathalpha}{cyrletters}{"58}
\newcommand{\Gal}{\operatorname{Gal}}
\newcommand{\length}{\operatorname{length}}
\newcommand{\Frac}{\operatorname{Frac}}
\newcommand{\Sel}{\operatorname{Sel}}
\newcommand{\Cl}{\operatorname{Cl}}
\newcommand{\Cone}{\operatorname{Cone}}
\newcommand{\Rgam}{R\Gamma}
\newcommand{\Cplus}{\mathcal C_S^+}
\newcommand{\CS}{\mathcal C_S}
\newcommand{\tensL}{\mathbin{\otimes}^{\mathbf L}}
\newcommand{\val}{v_2}
\newcommand{\ellp}{\ell_{\mathfrak p}}
\newcommand{\bsd}{\mathfrak d_2}
\newcommand{\Lalg}{L^{\mathrm{alg}}}
\newcommand{\LambdaQ}{\Lambda_{\Q_2}}
\newcommand{\kLambda}{k_\Lambda}
\newcommand{\Dat}{\mathcal D}
\newcommand{\Kt}{\widetilde K}
\newcommand{\FE}{\mathscr F_E}
\newcommand{\typ}{\operatorname{typ}}

\newcommand{\RF}{\hyperref[def:RF]{\textup{(RF)}}}
\newcommand{\CaseI}{\hyperref[case:I]{\textup{(I)}}}
\newcommand{\CaseII}{\hyperref[case:II]{\textup{(II)}}}
\newcommand{\CaseIII}{\hyperref[case:III]{\textup{(III)}}}
\theoremstyle{plain}
\newtheorem{maintheorem}{Theorem}

\newtheorem{theorem}{Theorem}[section]
\newtheorem{proposition}[theorem]{Proposition}
\newtheorem{lemma}[theorem]{Lemma}
\newtheorem{corollary}[theorem]{Corollary}
\theoremstyle{definition}
\newtheorem{definition}[theorem]{Definition}

\theoremstyle{remark}
\newtheorem{remark}[theorem]{Remark}
\numberwithin{equation}{section}
\BeforeBeginEnvironment{maintheorem}{\Needspace{4\baselineskip}}
\BeforeBeginEnvironment{theorem}{\Needspace{6\baselineskip}}
\BeforeBeginEnvironment{proposition}{\Needspace{6\baselineskip}}
\BeforeBeginEnvironment{lemma}{\Needspace{5\baselineskip}}
\BeforeBeginEnvironment{definition}{\Needspace{5\baselineskip}}

\begin{document}

\title[A rank-zero $2$-converse to Gross--Zagier and Kolyvagin]{A rank-zero $2$-converse to the Gross--Zagier and Kolyvagin theorem for elliptic curves with abelian two-division field}

\author{Ashay Burungale}
\address{Department of Mathematics, The University of Texas at Austin, 2515 Speedway, Stop C1200, Austin, TX 78712, USA}
\email{ashayk@utexas.edu}
\author{Ye Tian}
\address{Academy of Mathematics and Systems Science, Morningside Center of Mathematics, Chinese Academy of Sciences, Beijing 100190, China}
\email{ytian@math.ac.cn}

\date{November 7, 2026}
\subjclass[2020]{Primary 11G40; Secondary 11G05, 11R23, 11F37}

\begin{abstract}
Let \(E/\Q\) be an elliptic curve with abelian two-division field.
We prove that finiteness of \(\Sel_{2^\infty}(E/\Q)\) implies
\(L(E,1)\ne0\). The proof follows the quadratic twist deformation
approach of OpenAI: Beilinson--Kato elements are interpolated over the
cyclotomic \(\Z_2\)-extension and over families of quadratic twists,
and congruences in the resulting group ring transfer non-vanishing
from the twists to \(E\). The integral control of this interpolation
is governed by one cohomological \(\mu\)-invariant, which depends only
on \(E[2]\) and vanishes for abelian two-division fields. Together
with Kato's cyclotomic divisibility and the fixed lattice in his
construction, it bounds the denominator by a fixed power of \(2\),
independently of the twisting primes and the family dimension.
Congruences among coefficients of a half-integral-weight form, in the
style of Ono--Skinner, construct the required twist families.
\end{abstract}

\maketitle
{\hypersetup{linkcolor=black}\tableofcontents}

\section{Introduction}\label{sec:intro}

Let \(E/\Q\) be an elliptic curve of conductor \(N\), and put
\[
 s_2(E)=\operatorname{corank}_{\Z_2}\Sel_{2^\infty}(E/\Q).
\]
By the theorem of Gross--Zagier and Kolyvagin \cite{GZ,Kolyvagin},
\(L(E,1)\ne0\) implies \(s_2(E)=0\). For odd primes of good ordinary
reduction, the corresponding rank-zero \(p\)-converse
is known under residual and ramification hypotheses through Skinner--Urban
\cite[Theorem~2(b)]{SU}; for CM curves and \(p>3\), it follows from
Rubin \cite[\S11]{Rubin91}. We prove the \(2\)-converse when the
two-division field is abelian.

\begin{maintheorem}\label{thm:A}
Let \(E/\Q\) be an elliptic curve with \(\Q(E[2])/\Q\) abelian.
If \(s_2(E)=0\), then \(L(E,1)\ne0\).
\end{maintheorem}

An abelian subgroup of \(\operatorname{GL}_2(\F_2)\simeq S_3\)
has order one, two, or three, so the theorem covers curves with a
rational point of order two and curves with cyclic cubic two-division
field. Its CM case is due to the authors \cite{BT}.

We give a self-contained proof of Theorem~\ref{thm:A}; no result
of the OpenAI manuscripts is invoked. The proof follows the
quadratic twist deformation approach of \cite{OAPW,OAG}, retaining
the cyclotomic variable while varying the quadratic character. It
combines Kato's Euler system, which already contains the quadratic
deformation; Iwasawa theory over quadratic twist families; the
vanishing of one cohomological \(\mu\)-invariant, which depends only
on \(E[2]\) and holds for abelian two-division fields; and congruences
among Fourier coefficients of a half-integral-weight form, in the
manner of Ono--Skinner. The \(\mu\)-invariant organizes the integral
comparison. In the cyclic cubic case this gives a substantial
simplification: no Heegner points are needed. We hope that our
account makes the argument of \cite{OAPW,OAG} more transparent, and offers
some different viewpoints on it.\footnote{The OpenAI manuscripts
appear difficult to read; this has been our experience, and a number
of researchers in the field have communicated the same to us.}

In view of Smith's work on the Selmer counterpart of Goldfeld's
conjecture \cite{Smith17,Smith25}, Theorem~\ref{thm:A} and the theorem
of Gross--Zagier and Kolyvagin give the following consequence
(cf.~\cite[Theorem~1.2]{OAG}, \cite[Theorem~1.2]{BT}).\footnote{The
passage from integer parameters to signed square-free parameters
uses squareclass invariance and M\"obius inversion; the tail from
square divisors \(m^2\) with \(m>M\) is \(O(X/M)\).}

\begin{corollary}\label{cor:Goldfeld}
For an elliptic curve as in Theorem~\ref{thm:A}, the quadratic twists
of analytic rank zero have natural density \(1/2\) among signed
square-free twisting parameters, ordered by absolute value.
\end{corollary}

\subsection{Kato's main conjecture}\label{subsec:progression}

Let
\[
 \Q_\infty=\bigcup_{n\ge0}\Q(\zeta_{2^{n+2}})^+,\qquad
 \Gamma=\Gal(\Q_\infty/\Q)\simeq\Z_2,
\]
where the superscript \(+\) denotes the maximal real subfield.
Thus \(\Q_\infty/\Q\) is the real cyclotomic \(\Z_2\)-extension.
Choose a topological generator \(\gamma\), put \(t=\gamma-1\), and write
\[
 \Lambda=\Z_2[[\Gamma]]=\Z_2[[t]],\qquad
 \LambdaQ=\Lambda[1/2].
\]
The central specialization is \(t=0\).

Let \(\mathcal H_E^i\) denote primitive rational Iwasawa cohomology
of \(V_E=T_2(E)\otimes_{\Z_2}\Q_2\), and let
\(Z_E\subset\mathcal H_E^1\) be the module generated by the
Beilinson--Kato elements
(see Subsection~\ref{sec:kato} for definitions).\footnote{Primitive
cohomology is taken over \(\Z[1/2]\), as in \cite[\S12]{Kato}
and \cite[\S2.2]{BT}. Enlarging the support introduces the
corresponding localization factors.}
Kato proved that \(\mathcal H_E^1\) is free of rank one over
\(\LambdaQ\), while \(\mathcal H_E^2\) and
\(\mathcal H_E^1/Z_E\) are torsion
\cite[Theorems~12.4(1)--(2), 12.5(2)]{Kato}. The rational form of
his main conjecture is
\begin{equation}\label{eq:intro-main-conjecture}
 \operatorname{char}_{\LambdaQ}(\mathcal H_E^1/Z_E)
 =\operatorname{char}_{\LambdaQ}(\mathcal H_E^2).
\end{equation}
Selmer finiteness and this equality imply that the zeta element
does not vanish at the central specialization. Kato's reciprocity
law then gives \(L(E,1)\ne0\). This deduction from the main
conjecture is established in \cite{BT}.

Kato proved one divisibility. At every height-one prime
\(\mathfrak p\) of \(\LambdaQ\), that is, every height-one prime of
\(\Lambda\) other than \((2)\),
\begin{equation}\label{eq:intro-kato-divisibility}
 \ell_{\mathfrak p}(\mathcal H_E^2)
 \le \ell_{\mathfrak p}(\mathcal H_E^1/Z_E)
       +\ell_{\mathfrak p}(H^2_{\mathrm{Iw}}(\Q_2,V_E)),
\end{equation}
where \(\ell_{\mathfrak p}\) denotes length over
\((\LambdaQ)_{\mathfrak p}\) and
\(H^2_{\mathrm{Iw}}(\Q_2,V_E)\) is local Iwasawa cohomology.
The divisibility bounds the index of the zeta module from below.
It does not exclude an additional zero at \(t=0\), and it carries no
information at \((2)\): Kato's integral refinement
\cite[Theorem~12.5(4)]{Kato} assumes \(p\ne2\). The quadratic twist
deformation supplies the integrality and the congruences that exclude
this zero.

\subsection{Quadratic twist families and the uniform comparison}\label{subsec:statements}

Fix an embedding and normalize the valuation by
\[
 \overline\Q\hookrightarrow\overline\Q_2,\qquad
 \val(2)=1,\qquad \val(0)=+\infty.
\]
Fix a finite set \(\Sigma\) containing the primes dividing \(2N\).
For square-free \(h\), write \(\chi_h\) for the quadratic character
of \(\Q(\sqrt h)/\Q\), and \(E^h\) for the corresponding twist.

\begin{definition}\label{def:intro-family}
A parameter \(h\) is \emph{admissible} if it is a positive odd
fundamental discriminant, prime to \(\Sigma\), and a square in
\(\Q_\ell^\times\) for every \(\ell\in\Sigma\); we allow \(h=1\).
A quadratic twist family \(\mathcal F\) on \(X_b=\F_2^b\) consists
of a finite set of primes \(P\) outside \(\Sigma\) and linear forms
\(\alpha_q\in X_b^*=\operatorname{Hom}(X_b,\F_2)\), with
\[
 h_x=\prod_{q\in P}(q^*)^{\alpha_q(x)},\qquad
 q^*=(-1)^{(q-1)/2}q,\qquad E_x=E^{h_x}.
\]
Its \emph{dimension} is \(b=\dim_{\F_2}X_b\).
It is admissible if every \(h_x\) is admissible.
Put \(G_b=\operatorname{Hom}(X_b,\{\pm1\})\); evaluation at \(x\)
is denoted by a subscript \(x\).
\end{definition}

Integrality of an element of \(\Q_2[G_b]\) is a system of
congruences among its values. When \(b\) is large, these congruences
force one value to be non-zero once all the others are non-zero and
not too divisible by \(2\). This principle underlies the approach
(cf.~\cite[Lemma~11.4]{OAG}, \cite[Lemma~3.4]{OAPW}).

\begin{maintheorem}\label{thm:B}
Let \(C,H\in\Z\) with \(C+H\ge0\), and let \(v\in\Q_2[G_b]\) satisfy
\(2^Cv\in\Z_2[G_b]\). If \(v_x\ne0\) and \(\val(v_x)\le H\) for
every \(x\ne x_0\), then
\[
 b\ge2^{C+H}\quad\Longrightarrow\quad v_{x_0}\ne0,
 \qquad\val(v_{x_0})\le H.
\]
\end{maintheorem}

Write
\[
 L^{\mathrm{alg}}(E^h,1)=\frac{L(E^h,1)}{\Omega^+_{E^h}},\qquad
 c(h)=\prod_{q\mid h}c_q(E^h),\qquad c(1)=1,
\]
where \(\Omega^+\) is the positive N\'eron period and \(c_q\) is the
Tamagawa number at \(q\). The quadratic deformation of the
Beilinson--Kato elements is simple. Every quadratic character factors
through a cyclotomic field, so Kato's Euler system for the fixed
newform attached to \(E\) already contains the deformation: with
\(m=\prod_{q\in P}q\), the projection
\(\Gal(\Q(\mu_{m2^\infty})/\Q)\to\Gamma\times G_b\) induces an
integral map of group rings, which assembles Kato's elements with
coefficients in \(\Lambda[G_b]\) without dividing by \(\#G_b\)
(Proposition~\ref{prop:class}; cf.~\cite[\S3.3]{OAPW}).
The resulting interpolation element \(\mathscr L_{\mathcal F}\) lies
a priori in the total quotient ring of \(\Lambda[G_b]\)
(Definition~\ref{def:interpolation-element}). Write
\(\mathscr L_{\mathcal F,x}(t)\) for its value at \(x\).
The following theorem is the central result of the paper.

\begin{maintheorem}\label{thm:C}
Let \(E/\Q\) have abelian two-division field. There are integers
\(C_E\ge0\) and \(\beta_E\), independent of the twisting primes and
the family dimension,\footnote{The curve and its Beilinson--Kato
normalization are fixed (Subsection~\ref{sec:kato}).}
such that for every admissible quadratic twist family \(\mathcal F\)
the interpolation element lies in \(\LambdaQ[G_b]\) and satisfies
\begin{equation}\label{eq:intro-uniform}
 2^{C_E}\mathscr L_{\mathcal F}\in\Lambda[G_b].
\end{equation}
At a member \(x\) of Selmer corank zero,
\[
 \mathscr L_{\mathcal F,x}(0)=0\quad\Longleftrightarrow\quad L(E_x,1)=0.
\]
If also \(h_x\ne1\) and \(L(E_x,1)\ne0\), then
\begin{equation}\label{eq:introcomparison}
 \val(\mathscr L_{\mathcal F,x}(0))
 =\val\left(\frac{L^{\mathrm{alg}}(E_x,1)}
 {\#\Sha(E_x/\Q)[2^\infty]c(h_x)}\right)+\beta_E.
\end{equation}
\end{maintheorem}

\begin{remark}\label{rem:intro-RF}
The main text proves this theorem under the weaker condition \RF{}:
the degree-two cyclotomic cohomology with the real-place condition
has \(\mu\)-invariant zero
(Definition~\ref{def:RF}, Lemma~\ref{lem:RF-mu}).
This depends only on \(E[2]\) and holds for every abelian
two-division field (Proposition~\ref{prop:RF-cases}).
Condition \RF{} combines the vanishing of the \(\mu\)-invariant of
the fine Selmer group with a condition at the real places. For curves
with a rational point of order two, the former was obtained
by Lin--Yin \cite[Theorem~1.4]{LY26} (Remark~\ref{rem:fine-selmer}).
\end{remark}

\begin{remark}\label{rem:intro-BSD}
The Birch--Swinnerton-Dyer formula predicts that the valuation of
the ratio in \eqref{eq:introcomparison} is constant for non-trivial
admissible twists: admissibility fixes the old local factors, the
period comparison and the rational \(2\)-power torsion.
\end{remark}

Taking a one-dimensional family gives, for admissible \(h\ne1\)
with \mbox{\(L(E^h,1)\ne0\)},
\begin{equation}\label{eq:introlower}
 \val\left(\frac{L^{\mathrm{alg}}(E^h,1)}{c(h)}\right)
 \ge \val(\#\Sha(E^h/\Q)[2^\infty])-C_E-\beta_E.
\end{equation}

The Shimura--Shintani--Waldspurger correspondence
\cite{Shimura73,Shintani,Wald80} relates the weight-two newform
attached to \(E\) to cusp forms of weight \(3/2\).
Waldspurger's formula \cite{Wald} turns the central values
\(L(E^h,1)\) into squares of Fourier coefficients of such a form, and
\eqref{eq:introlower} into a lower bound for the coefficients
(cf.~\cite[Proposition~7.1, Lemma~7.5]{OAPW}, \cite[\S4.4]{OAG}).
Here the lower bound comes from Theorem~\ref{thm:C}.

\begin{maintheorem}\label{thm:D}
Suppose that \(E\) has root number \(+1\). There are an admissible
\(n_0\) with \mbox{\(L(E^{n_0},1)\ne0\)}, a character \(\psi\) of \(2\)-power
order, and a cusp form \(g=\sum a(n)q^n\) of weight \(3/2\), with
\(a(n_0)=1\), such that for admissible \(h\),\footnote{Characters are
extended multiplicatively to rational numbers prime to their
conductors: \(\psi(a/b)=\psi(a)\psi(b)^{-1}\).}
\[
 a(h)^2=\psi(h/n_0)\frac{L^{\mathrm{alg}}(E^h,1)}
 {L^{\mathrm{alg}}(E^{n_0},1)}.
\]
The numbers \(\val(a(h))-\tfrac12\val(c(h))\), for admissible \(h\)
with \(a(h)\ne0\), are integers. If \(E\) satisfies \RF{}, they
are bounded below and have a minimum
\[
 \delta=\min_{\substack{h\ \mathrm{admissible}\\a(h)\ne0}}
 \left\{\val(a(h))-\tfrac12\val(c(h))\right\}.
\]
\end{maintheorem}

For the family construction, take \(E\) with full rational two-torsion,
with one rational point of order two and no isogenous curve having
full rational two-torsion, or with cyclic cubic two-division field
(cases \CaseI{}--\CaseIII{} of Section~\ref{sec:symbol}). The
following family construction develops \cite[\S\S7--8]{OAPW} and
\cite[\S\S6, 8--9]{OAG}.

\begin{maintheorem}\label{thm:E}
For \(E\) in one of these three cases, with root number \(+1\), if
\(\val(a(1))>\delta\), then for every \(b\ge1\) there is an admissible
quadratic twist family on \(\F_2^b\) such that, for every \(x\ne0\),
\[
 a(h_x)\ne0,\qquad
 \val(a(h_x))-\tfrac12\val(c(h_x))=\delta.
\]
\end{maintheorem}

\subsection{About the proof}\label{subsec:strategy}

The essential step is to retain the cyclotomic variable while moving
in quadratic twist families. Iwasawa theory then controls one
integral determinant over the whole family, and coefficient
congruences choose its non-zero specializations. Uniformity in both
the twisting primes and the dimension permits the use of
Theorem~\ref{thm:B}.

\subsubsection{Cyclotomic cohomology and the integral comparison}
The cyclotomic variable absorbs the local terms that grow with the
family. At a twisting prime \(q\), the Euler factor is a polynomial
in the Frobenius of \(q\), which acts on \(\Lambda\) through a
non-zero cyclotomic exponent, since no odd prime splits completely in
\(\Q_\infty\). Its reduction is therefore non-zero in
\(\Lambda/2\Lambda\), so it is a unit in the discrete valuation ring
\(\mathcal O=\widehat{\Lambda_{(2)}}\), and the local cohomology at
\(q\) vanishes over \(\mathcal O\) (Lemma~\ref{lem:residual-rank}).
These terms return at \(t=0\): the Euler factors at inactive primes
cancel between the zeta element and global cohomology, and the active
primes contribute the Tamagawa numbers \(c(h_x)\) in
\eqref{eq:introcomparison}.

Away from \((2)\), Kato's divisibility, applied to each twist
separately, shows that \(\mathscr L_{\mathcal F}\) has no poles:
\(\mathscr L_{\mathcal F}\in\LambdaQ[G_b]\). For non-CM curves the
Euler-system theorem is applied to the rank-two lattice of each
twisted newform. The errors in normalizing its Beilinson--Kato
elements and comparing lattices are powers of \(2\) that need not be
bounded as the twist varies; they disappear after inverting \(2\).
For CM curves, \cite{BT} gives the same rational statement.

At \((2)\), the uniform estimate comes from the fixed newform. The
relevant vanishing is
\begin{equation}\label{eq:intro-mu}
 \RF{}\quad\Longleftrightarrow\quad
 \mu_\Lambda(\mathcal H^2_{E,+})=0,
\end{equation}
where \(\mathcal H^2_{E,+}\) is degree-two cyclotomic cohomology
with the real-place condition
(Definition~\ref{def:positive}, Lemma~\ref{lem:RF-mu}).
It removes degree-two cohomology and makes degree-one cohomology
free of rank two over \(\mathcal O[G_b]\). The Beilinson--Kato
elements of the fixed newform, pushed into the family, are integral
up to a fixed power of \(2\), which compares Kato's lattice for the
newform with \(T_2(E)\) and does not depend on the family.\footnote{Kato's
zeta morphism \cite[Theorem~12.5(1)]{Kato} is defined on the Betti
lattice of a single newform. The family elements arise from that of
the fixed newform through integral coefficient maps
(Proposition~\ref{prop:class}), which suggests the uniformity in
Theorem~\ref{thm:C}.} Combining the two comparisons gives
\[
 2^{C_E}\mathscr L_{\mathcal F}\in
 \mathcal O[G_b]\cap\LambdaQ[G_b]=\Lambda[G_b].
\]
Kato's reciprocity law and Poitou--Tate duality then compute the
central values, giving the rest of Theorem~\ref{thm:C}.

\subsubsection{Congruences and the \texorpdfstring{\(\mu\)}{mu}-invariant}
Greenberg--Vatsal show that, under their residual hypotheses and
analytic and algebraic \(\mu=0\), the cyclotomic main conjecture is
equivalent for congruent ordinary newforms. Congruences and Kato's
divisibility give this transfer \cite[Corollary~3.14]{GV}.
For Kato's main conjecture, see also \cite{KLP19,Nakamura}.
Here the \(\mu\)-invariant is that of cohomology, and it depends only
on \(E[2]\). Quadratic twists have the same \(E[2]\), so one
verification of \RF{} controls the entire deformation. The integral
element \(2^{C_E}\mathscr L_{\mathcal F}(0)\in\Z_2[G_b]\) encodes
congruences: its iterated differences of order \(k\) are divisible
by \(2^k\). Theorem~\ref{thm:B} turns them into non-vanishing at
the original curve, a horizontal form of the Greenberg--Vatsal
principle. The use of Kato's main conjecture for the converse goes
back to \cite{BT}.

The integral coefficient maps also suggest a zeta morphism for
quadratic deformations. Nakamura's zeta morphisms address universal
deformations under different residual hypotheses \cite{Nakamura}
(Remark~\ref{rem:kato-equivariant}).

\subsubsection{Families at the minimum}
Theorem~\ref{thm:B} requires, for every \(b\), a family on
\(\F_2^b\) whose non-zero members have non-zero central values of
bounded valuation. Waldspurger's formula and \eqref{eq:introlower}
bound the normalized coefficient valuations from below, hence give
the minimum \(\delta\). Reduction at the minimum attaches to each
admissible set \(U\) of primes a symbol \(A(U)\in\F_2\), equal to one
exactly when \(U\) attains \(\delta\). Following the
coefficient-congruence method of Ono--Skinner \cite{OS98a,OS98b},
multiplication by \(\theta\equiv1\pmod2\) relates this symbol to
integral-weight forms and their Galois representations.

For rational two-torsion, \(A(U)\) then depends only on quadratic
residue symbols among the primes of \(U\) and the fixed support, and
it obeys a contraction rule that merges two primes into one.
Theorem~\ref{thm:E} becomes a combinatorial problem on graphs of
residue symbols, to attain the minimum at every non-zero index of a
quadratic twist family; auxiliary primes and the removal of pairs of
inert primes solve it. In the cyclic cubic case, Theorem~\ref{thm:C}
also bounds the \(2\)-Selmer group of every twist at the minimum.
Up to a bounded correction, that group contains the kernel of a
matrix of local conditions. In a configuration attaining the minimum,
built from mutually uncoupled blocks of primes, only boundedly many
blocks therefore have singular matrix. A maximality argument over
disjoint blocks then gives simultaneous attainment. Coefficient
congruences thus construct the auxiliary family; the group-ring
congruences of Theorem~\ref{thm:B} transfer non-vanishing to its
missing member.

For these families, \eqref{eq:introcomparison} gives the second
bound needed in Theorem~\ref{thm:B}:
\begin{equation}\label{eq:mechanism}
 \begin{gathered}
 -C_E\le\val(\mathscr L_{\mathcal F,x}(0))
 =H-\val(\#\Sha(E_x/\Q)[2^\infty])\le H\qquad(x\ne0),\\
 H=2\delta+\val(L^{\mathrm{alg}}(E^{n_0},1))+\beta_E.
 \end{gathered}
\end{equation}
The same identity bounds the Shafarevich--Tate contribution by
\(C_E+H\), yielding the Selmer bound used in the cubic construction.
In particular, \(C_E+H\ge0\), independently of the dimension.

\begin{proof}[Proof of Theorem~\ref{thm:A}]
Assume \(s_2(E)=0\). Replace \(E\) by an isogenous curve with full
rational two-torsion if such a curve exists; this preserves
\(s_2(E)\) and \(L(E,s)\). Monsky's parity theorem \cite{Monsky}
gives root number \(+1\). If \(L(E,1)=0\), then \(a(1)=0\).
Choose the family of Theorem~\ref{thm:E} with \(b\ge2^{C_E+H}\).
For every \(x\ne0\), Waldspurger's formula gives \(L(E_x,1)\ne0\),
hence Selmer corank zero by Gross--Zagier and Kolyvagin.
Put \(v=\mathscr L_{\mathcal F}(0)\in\Q_2[G_b]\).
Theorem~\ref{thm:C} gives \(2^{C_E}v\in\Z_2[G_b]\), and
\eqref{eq:mechanism} gives \(v_x\ne0\) and \(\val(v_x)\le H\)
for every \(x\ne0\). Theorem~\ref{thm:B} gives \(v_0\ne0\).
Since \(s_2(E)=0\), Theorem~\ref{thm:C} implies \(L(E,1)\ne0\),
a contradiction.
\end{proof}

\begin{remark}[The trivial character]\label{rem:trivial-character}
Let \(E/\Q\) have abelian two-division field and root number \(+1\).
The proof also gives
\[
 \operatorname{char}_{\LambdaQ}(\mathcal H_E^1/Z_E)\,
 (\LambdaQ)_{(t)}
 =\operatorname{char}_{\LambdaQ}(\mathcal H_E^2)\,
 (\LambdaQ)_{(t)}.
\]
Indeed, if \(L(E,1)=0\), repeat the proof through the application
of Theorem~\ref{thm:B}. This produces a family with
\(\mathscr L_{\mathcal F,0}(0)\ne0\) without using \(s_2(E)=0\).
If \(L(E,1)\ne0\), Gross--Zagier--Kolyvagin and
Theorem~\ref{thm:C} give the same non-vanishing.
Lemma~\ref{lem:central-order} then identifies the two orders at
\(t=0\). In particular, \(\mathscr L_{\mathcal F,0}(0)\ne0\)
for every admissible family. The isogeny reduction in the proof
preserves the rational cohomology and zeta modules.
\end{remark}

\subsection{Relation to the OpenAI manuscripts}\label{subsec:comparison}

Our proof retains the quadratic twist deformation and coefficient
congruences of \cite{OAPW,OAG}. The cohomological \(\mu\)-invariant
gives a common explanation for the coefficient lower bound, the
Selmer bound at the minimum and integral transfer: all three follow
from one Beilinson--Kato comparison. For rational two-torsion, the
uniform integral bound is already proved in
\cite[(3.11)--(3.12), Lemmas~3.3 and~3.6]{OAPW}.

The main simplification concerns a cyclic cubic two-division field.
Ferrero--Washington verifies \RF{}, so the same comparison applies
at all new primes. No Heegner points or auxiliary twist of analytic
rank one enter our proof. By contrast, \cite{OAG} obtains the
coefficient lower bound from its genus Heegner-point bound and uses
a bounded clearing factor in the transfer, which requires a uniform
bound on \(v_2(q^2-1)\) at the new primes.\footnote{See
\cite[Theorem~3.6, Proposition~3.8(ii), \S13]{OAG} for the genus
bound and its use of the explicit Gross--Zagier formula, and
\cite[Theorem~3.1]{OAG} for the transfer.}

The cubic trace identities, least-weight transfer, Selmer matrices
and calibration adapt \cite[\S\S5.5--5.6, 8--9]{OAG}.
Bounded maximality on actual prime configurations replaces Ramsey
stabilization. We also shorten the rational two-torsion
constructions. The group-ring threshold \(b\ge2^{C+H}\) is that of
\cite[Lemma~11.4]{OAG}; Remark~\ref{rem:sharp} proves its sharpness.
The exact central formula of Theorem~\ref{thm:central} refines
\cite[Proposition~3.7]{OAPW} and \cite[Proposition~12.3]{OAG}.
The scope of \cite{OAG} is broader: it also treats \(S_3\)
two-division fields, rank one and Goldfeld's conjecture.
The other cases of the \(p\)-converse announced by OpenAI will be
discussed elsewhere.

\subsection{Leitfaden and antecedents}\label{subsec:organization}

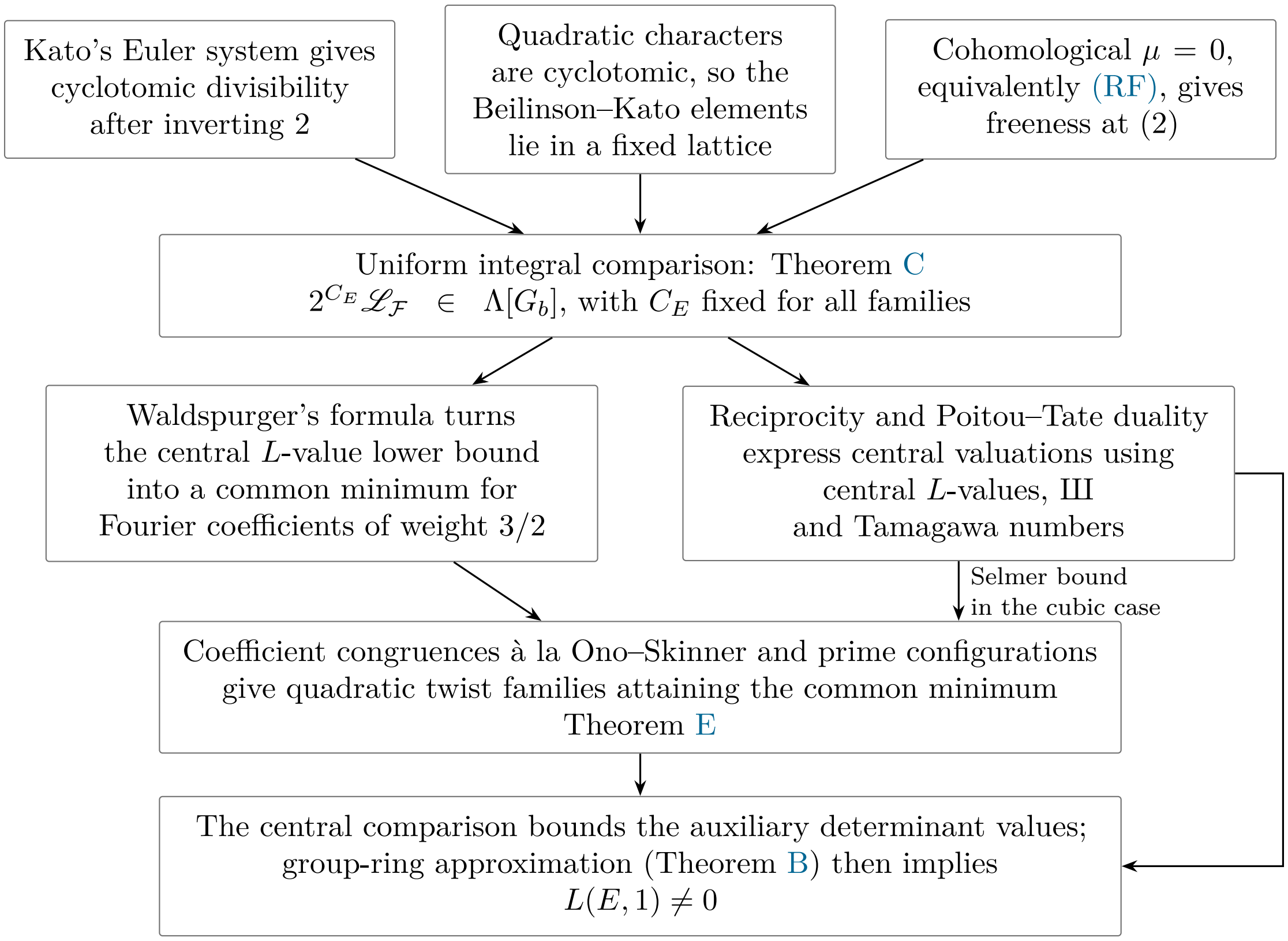
\begin{figure}[ht]
\centering
\begin{tikzpicture}[
 box/.style={draw=black!55,rounded corners=1pt,align=center,
             text width=4.05cm,inner sep=6pt,font=\small},
 arr/.style={-{Stealth[length=1.8mm]},line width=.6pt}]
\node[box] (a) at (-5.05,0) {Kato's Euler system gives\\cyclotomic divisibility\\after inverting \(2\)};
\node[box] (b) at (0,0) {Quadratic characters\\are cyclotomic, so the\\Beilinson--Kato elements\\lie in a fixed lattice};
\node[box] (c) at (5.05,0) {Cohomological \(\mu=0\),\\equivalently \RF{}, gives\\freeness at \((2)\)};
\node[box,text width=10.6cm] (d) at (0,-2.25) {Uniform integral comparison: Theorem~\ref{thm:C}\\
\(2^{C_E}\mathscr L_{\mathcal F}\in\Lambda[G_b]\), with \(C_E\) fixed for all families};
\node[box,text width=5.9cm] (e) at (-3.65,-4.4) {Waldspurger's formula turns\\the central \(L\)-value lower bound\\into a common minimum for\\Fourier coefficients of weight \(3/2\)};
\node[box,text width=5.9cm] (f) at (3.65,-4.4) {Reciprocity and Poitou--Tate duality\\express central valuations using\\central \(L\)-values, \(\Sha\)\\and Tamagawa numbers};
\node[box,text width=10.6cm] (g) at (0,-6.85) {Coefficient congruences \`a la Ono--Skinner and prime configurations\\give quadratic twist families attaining the common minimum\\Theorem~\ref{thm:E}};
\node[box,text width=10.6cm] (h) at (0,-8.9) {The central comparison bounds the auxiliary determinant values;\\group-ring approximation (Theorem~\ref{thm:B}) then implies\\\mbox{\(L(E,1)\ne0\)}};
\draw[arr] (a)--(d); \draw[arr] (b)--(d); \draw[arr] (c)--(d);
\draw[arr] (d)--(e); \draw[arr] (d)--(f);
\draw[arr] (e)--(g); \draw[arr] (f.south)--node[right,font=\scriptsize,align=left] {Selmer bound\\in the cubic case} (g.north -| f.south);
\draw[arr] (g)--(h);
\draw[arr] (f.east)--++(.55,0)|-(h.east);
\end{tikzpicture}
\caption{The proof of the converse. Theorem~\ref{thm:C} supplies
the lower bound defining the coefficient minimum and the Selmer
bound used in the cubic construction. On the resulting families,
the central formula gives the upper bound required by
Theorem~\ref{thm:B}.}\label{fig:flow}
\end{figure}

Section~\ref{sec:principle} proves Theorem~\ref{thm:B}.
Section~\ref{sec:families} constructs the Beilinson--Kato
interpolation and proves Theorem~\ref{thm:C}.
Section~\ref{sec:symbol} proves Theorem~\ref{thm:D}, and the next
three sections establish Theorem~\ref{thm:E}.
Appendix~\ref{sec:residual-appendix} gives further criteria for \RF{},
and Appendix~\ref{sec:congruent-example} illustrates the families and
dimension bound for congruent-number curves.
The flow of the proof is displayed in Figure~\ref{fig:flow}.

For an account of the congruent-number literature, see \cite{BT21}.
The valuation and residue-symbol methods originate in \cite{Zhao97,Zhao01,Feng,Tian14,TYZ}.
Related results on quadratic twists, central \(L\)-values and Iwasawa
theory may be found in \cite{SU,CGLS,CLTZ,CLZ,KL19,Zhai20,Zhai25}.

\subsection*{Acknowledgements}
The authors are grateful to Alex Smith for drawing their attention to
the uniform integrality of Beilinson--Kato elements in quadratic twist families.
The authors also thank David Ben-Zvi, Manjul Bhargava, Frank Calegari, Francesc Castella, Daniel Disegni,
Francesco Maggi, Marco Sangiovanni Vincentelli, Christopher Skinner,
Alex Smith and Wei Zhang for helpful discussions and comments.

\subsection*{Use of AI tools}
ChatGPT and Claude were used to assist in reading the OpenAI
manuscripts and in drafting and revising the exposition. The authors
are responsible for the mathematical arguments and the final text.

\FloatBarrier

\section{The approximation principle in group rings}\label{sec:principle}

Theorem~\ref{thm:B} follows the congruence method of
\cite[Lemma~11.4]{OAG} (cf.~\cite[Lemma~3.4]{OAPW}). We prove it
using the finite-difference criterion of
Lemma~\ref{lem:differences} and the congruence bound of
Lemma~\ref{lem:count}. Remark~\ref{rem:sharp} records the sharpness
of the threshold.

Throughout, \(X_b=\F_2^b\), \(G_b=\operatorname{Hom}(X_b,\{\pm1\})\), and \(v_x\) is the value at \(x\) of \(v\in\Q_2[G_b]\). Functions on \(X_b\) have unique multilinear polynomial representatives in \(x_1,\dots,x_b\), with \(x_i^2=x_i\). Put \(x_I=\prod_{i\in I}x_i\), let \(1_J\) be the indicator vector of \(J\), and let \(s_1,\dots,s_b\) be the dual basis of \(G_b\), so \(x(s_i)=(-1)^{x_i}=1-2x_i\).

\subsection{Integrality and finite differences}\label{subsec:integrality}

\begin{lemma}\label{lem:differences}
For \(y\in\Q_2[G_b]\) write \(y_x=\sum_{I}c_I(y)\,x_I\). Then \(y\in\Z_2[G_b]\) if and only if \(c_I(y)\in2^{|I|}\Z_2\) for every \(I\subseteq\{1,\dots,b\}\). Since \(c_I(y)=\sum_{J\subseteq I}(-1)^{|I\smallsetminus J|}y_{1_J}\), integrality of \(y\) means that the iterated differences of order \(k\) of its values are divisible by \(2^k\).
\end{lemma}

\begin{proof}
For \(I\subseteq\{1,\dots,b\}\) put \(f_I=\prod_{i\in I}(1-[s_i])\). Its value at \(x\) is \(\prod_{i\in I}(1-x(s_i))=2^{|I|}x_I\). As \(\Z_2[G_b]=\bigotimes_i\Z_2[s_i]/(s_i^2-1)\) and \(1,1-[s_i]\) is a basis of each factor, the \(f_I\) form a basis of \(\Z_2[G_b]\) over \(\Z_2\) and of \(\Q_2[G_b]\) over \(\Q_2\). Writing \(y=\sum_Id_If_I\), the values are \(y_x=\sum_I2^{|I|}d_Ix_I\), so \(c_I(y)=2^{|I|}d_I\), and \(y\) is integral exactly when every \(d_I\) is. For \(I=\varnothing\) the condition is \(y_0\in\Z_2\). The last formula is M\"obius inversion of \(y_{1_J}=\sum_{I\subseteq J}c_I(y)\).
\end{proof}

With \(\varepsilon_i=1-[s_i]\) one has
\[
 \Z_2[G_b]=\Z_2[\varepsilon_1,\dots,\varepsilon_b]/(\varepsilon_1^2-2\varepsilon_1,\dots,\varepsilon_b^2-2\varepsilon_b),\qquad f_I=\prod_{i\in I}\varepsilon_i,
\]
whose reduction is \(\F_2[\varepsilon_1,\dots,\varepsilon_b]/(\varepsilon_1^2,\dots,\varepsilon_b^2)\). All characters reduce to \(\varepsilon_i=0\), while the nilpotent variables retain the higher congruences \(c_I(y)\in2^{|I|}\Z_2\). Integral character values alone give only \(c_I(y)\in\Z_2\). For example, the top difference of \(1-2^{-b}\sum_g[g]\) is a unit, so this element is not integral.

\subsection{Congruences among the values}
The following degree bound is the essential congruence.

\begin{lemma}\label{lem:count}
Let \(y\in\Z_2[G_b]\), \(M\ge1\) and \(x_0\in X_b\). The indicator of
\[
 S=\{x\in X_b:y_x\equiv y_{x_0}\pmod{2^M}\}
\]
is a polynomial function on \(X_b=\F_2^b\) of degree at most \(2^{M-1}-1\). Consequently, if \(b\ge2^{M-1}\), then \(S\) has even cardinality and contains an element other than \(x_0\).
\end{lemma}

\begin{proof}
Replacing \(y=\sum_gy_g[g]\) by \(\sum_gy_gx_0(g)[g]\), whose value at \(x\) is \(y_{x+x_0}\), we may assume \(x_0=0\) (cf.~\cite[Lemma~11.4]{OAG}). By Lemma~\ref{lem:differences},
\[
 y_x-y_0=\sum_{\varnothing\ne I\subseteq\{1,\dots,b\}}2^{|I|}d_I\,x_I,\qquad d_I\in\Z_2 .
\]
Put \(Q(x)=(y_x-y_0)/2=\sum_I2^{|I|-1}d_Ix_I\) and \(m=2^{M-1}\). In \(\F_2[Z]/(Z^m)\) let \(U=1+Z\); since \(U^m=1+Z^m=1\), the power \(U^a\) is defined for \(a\in\Z_2\) through \(a\bmod m\). Let \(P(x)\) be the coefficient of \(Z^{m-1}\) in \(U^{Q(x)-1}\). If \(Q(x)\equiv0\pmod m\), this power is \((1+Z)^{m-1}\) and \(P(x)=1\); if \(Q(x)\equiv c\pmod m\) with \(1\le c<m\), it is \((1+Z)^{c-1}\), of degree less than \(m-1\), and \(P(x)=0\). Thus \(P\) is the indicator of \(S\). As a function of \(x\),
\[
 U^{Q(x)-1}=U^{-1}\prod_{I\ne\varnothing}\Bigl[1+\bigl(U^{2^{|I|-1}d_I}-1\bigr)x_I\Bigr],
\]
and \(U^{2^{|I|-1}}=1+Z^{2^{|I|-1}}\) in characteristic two, so the coefficient of \(x_I\) in its factor is divisible by \(Z^{2^{|I|-1}}\), while \(|I|\le2^{|I|-1}\); the factor \(U^{-1}\) does not involve \(x\). A product of such coefficients that contributes to \(Z^{m-1}\) therefore has degree at most \(m-1\), and \(\deg P\le m-1\). If \(b\ge m\), every monomial of \(P\) omits a variable, so its sum over \(\F_2^b\) vanishes in \(\F_2\); hence \(\#S\) is even, and as \(0\in S\) it contains another index.
\end{proof}

\begin{proof}[Proof of Theorem~\ref{thm:B}]
Put \(y=2^Cv\in\Z_2[G_b]\) and \(M=C+H+1\). For \(x\ne x_0\), \(y_x\ne0\) and \(\val(y_x)\le C+H<M\). Since \(b\ge2^{M-1}\), Lemma~\ref{lem:count} gives \(x\ne x_0\) with \(y_x\equiv y_{x_0}\pmod{2^M}\), so \(\val(y_{x_0})=\val(y_x)\le C+H\). Hence \(v_{x_0}\ne0\) and \(\val(v_{x_0})\le H\).
\end{proof}

\begin{remark}\label{rem:sharp}
The threshold is sharp. For \(C+H\ge1\) and
\(b=2^{C+H}-1\), the element
\(v=2^{-C}\sum_{i=1}^b(1-[s_i])\) satisfies
\(2^Cv\in\Z_2[G_b]\), \(v_0=0\), and
\(v_x=2^{1-C}|x|\ne0\) of valuation at most \(H\) elsewhere.
The bound on group-ring denominators cannot be dropped: the element
\(1-2^{-b}\sum_g[g]\) of Subsection~\ref{subsec:integrality} has
values zero and one.
\end{remark}

\section{Beilinson--Kato elements in quadratic twist families}\label{sec:families}
Theorem~\ref{thm:main} is the principal result of this section. We
construct the interpolated Beilinson--Kato element, prove its integral bound, and
compute the central specializations by reciprocity and global duality.
Proposition~\ref{prop:RF-cases} verifies the residual hypothesis
for every curve in Theorem~\ref{thm:A}.
All constants depend on \(E\) and fixed normalizations only.
Throughout, \(E\) has conductor \(N\), and \(E^h\) denotes its quadratic
twist by \(\chi_h\). We impose \RF{} only for the uniform integral bound.

\subsection{Quadratic twist families}

Let \(\Sigma_0\) be the set of primes dividing \(2N\), and let \(\Sigma\supseteq\Sigma_0\) be finite. Let \(X_b=\F_2^b\), \(G_b=\operatorname{Hom}(X_b,\{\pm1\})\), and write \(x(g)\) for their pairing.

\begin{definition}\label{def:admissible}
A \emph{quadratic twist family} \(\mathcal F\) on \(X_b\) is the data
\((P,(\alpha_q)_{q\in P})\) of Definition~\ref{def:intro-family}.
We call \(\alpha_q\) the \emph{ramification form} at \(q\); the
\emph{members} are \(E_x=E^{h_x}\), and the \emph{dimension} is
\(\dim_{\F_2}X_b=b\). It is \emph{admissible}, or
\emph{\(\Sigma\)-admissible} when the support is to be stressed, if
each \(h_x\) is admissible for \(\Sigma\) as in that definition.
\end{definition}

A prime \(q\) ramifies at \(x\), and is called \emph{active} at \(x\), exactly when \(\alpha_q(x)=1\). As the \(\alpha_q\) are linear, \(h_0=1\), so the original curve belongs to the family, and \(h_{x+y}=h_xh_y\) up to squares. The corresponding Galois representations satisfy\footnote{A quadratic
twisting isomorphism has Galois cocycle \(\chi_{h_x}\in\{\pm1\}\).
This gives the first isomorphism; \(-1\) acts trivially on
two-torsion, giving the second.}
\[
 T_2(E_x)\simeq T_2(E)\otimes\chi_{h_x},\qquad E_x[2]\simeq E[2].
\]
In an admissible family each \(E_x\) is locally isomorphic to \(E\) at \(\Sigma\) and at infinity. The quadratic deformation therefore preserves the local Tate modules at these places.

\begin{remark}\label{rem:admissibility}
Admissibility keeps the local representations at the primes dividing
\(2N\), the period comparison, and complex conjugation fixed. It
also places every twisting parameter in the same local squareclasses
for Waldspurger's formula. The condition therefore serves both the
cyclotomic comparison and the comparison of Fourier coefficients.
\end{remark}

\subsection{Central \texorpdfstring{\(L\)}{L}-values}

Let \(\omega_{E^h}\) be a minimal differential and \(\Omega^+_{E^h}>0\) its integral over a generator of the invariant part of \(H_1(E^h(\C),\Z)\) under complex conjugation; this is the period in the central \(L\)-value
\[
 \Lalg(E^h,1)=\frac{L(E^h,1)}{\Omega^+_{E^h}}\in\Q
\]
of Subsection~\ref{subsec:statements}. For admissible \(h\) one has \(\Omega^+_{E^h}=\Omega^+_E/\sqrt h\) (Lemma~\ref{lem:period}). A non-zero central \(L\)-value implies that \(\Sel_{2^\infty}(E^h)\) has corank zero, by the theorem of Gross--Zagier and Kolyvagin \cite{GZ,Kolyvagin}; see \cite[\S1.0.1]{BT}. For an odd prime \(q\nmid N\) let \(t_q=\dim_{\F_2}E[2]^{\mathrm{Fr}_q}\in\{0,1,2\}\), where \(\mathrm{Fr}_q\) is arithmetic Frobenius. At \(q\mid h\) the twist has reduction type \(I_0^*\) and Tamagawa number \(c_q(E^h)=\#E(\F_q)[2]=2^{t_q}\).

\begin{definition}\label{def:defect}
Let \(h\) be admissible, and put \(\tau(h)=\sum_{q\mid h}t_q\), so that \(\tau(h)=\val(c(h))\). If \(\Sel_{2^\infty}(E^h)\) has corank zero, so that \(E^h(\Q)\) and \(\Sha(E^h/\Q)[2^\infty]\) are finite, put \(s_{\Sha}(h)=\val(\#\Sha(E^h/\Q)[2^\infty])\). If moreover \(L(E^h,1)\ne0\), put
\begin{equation}\label{eq:defects}
 D(h)=\val\bigl(\Lalg(E^h,1)\bigr)-s_{\Sha}(h)-\tau(h).
\end{equation}
\end{definition}

\subsection{Integral interpolation}

Let \(\Q_\infty/\Q\) be the real cyclotomic \(\Z_2\)-extension and
\(\Gamma=\Gal(\Q_\infty/\Q)\). Choose a topological generator
\(\gamma\) of \(\Gamma\), and put \(t=\gamma-1\). Write
\begin{equation}\label{eq:rings}
 \begin{gathered}
 \Lambda=\Z_2[[\Gamma]]=\Z_2[[t]],\\
 R=\Lambda[G_b],\qquad \LambdaQ=\Lambda[\tfrac12].
 \end{gathered}
\end{equation}
For \(v=\sum_{g\in G_b}v_g(t)[g]\in\LambdaQ[G_b]\), put \(v_x(t)=\sum_gv_g(t)x(g)\). Its \emph{central specialization} at \(x\) is \(v_x(0)\).

Let \(f=\sum_{n\ge1}a_nq^n\) be the normalized newform attached
to \(E\). Fix its cohomological lattice and the comparison with
\(T_2(E)\). These determine the normalization of the
Beilinson--Kato elements used below.

\begin{theorem}\label{thm:main}
Assume \RF{}. Fix the normalization of the Beilinson--Kato element in Subsection~\ref{sec:kato}. There are integers \(C_E\ge0\) and \(\beta_E\), depending only on \(E\) and this normalization, with the following property. For every admissible family, the interpolation element \(\mathscr L_{\mathcal F}\) of Definition~\ref{def:interpolation-element} lies in \(\LambdaQ[G_b]\) and satisfies
\[
 2^{C_E}\mathscr L_{\mathcal F}\in\Lambda[G_b];
\]
at every index \(x\) for which \(\Sel_{2^\infty}(E^{h_x})\) has corank zero,
\[
 \mathscr L_{\mathcal F,x}(0)=0\iff L(E^{h_x},1)=0,
\]
and, if \(h_x\ne1\) and \(L(E^{h_x},1)\ne0\),
\begin{equation}\label{eq:exactcentral}
 \val\bigl(\mathscr L_{\mathcal F,x}(0)\bigr)=D(h_x)+\beta_E.
\end{equation}
\end{theorem}

\begin{corollary}\label{cor:sha}
Assume \RF{}. For every admissible \(h\ne1\) with \(L(E^h,1)\ne0\),
\[
 \val\bigl(\Lalg(E^h,1)\bigr)\ \ge\ s_{\Sha}(h)+\tau(h)-C_E-\beta_E .
\]
\end{corollary}
\begin{proof}
Take \(b=1\), let \(P\) be the set of primes dividing \(h\), and let every \(\alpha_q\) be the non-zero form. Then \(h_1=\prod_{q\mid h}q^*=h\), since both sides are congruent to \(1\) modulo \(4\) and have the same absolute value. The twist \(E^h\) has Selmer corank zero, and Theorem~\ref{thm:main} gives \(D(h)+\beta_E=\val(\mathscr L_{\mathcal F,1}(0))\ge-C_E\).
\end{proof}

\begin{remark}\label{rem:earlier-valuation-bounds}
For earlier bounds on central \(L\)-values of quadratic twists, see
\cite{Zhao97,Zhao01,Tian14,Zhai20,Zhai25}.
\end{remark}

\subsection{Cyclotomic freeness}\label{sec:positive}

Proposition~\ref{prop:free} proves freeness of rank two for
degree-one cohomology over \(\mathcal O[G_b]\).
Its residual hypothesis is equivalent to cohomological
\(\mu=0\) (Lemma~\ref{lem:RF-mu}) and is verified in
Proposition~\ref{prop:RF-cases}.

\subsubsection{The modified global complex}

Put \(T_E=T_2(E)\), \(V_E=T_E\otimes\Q_2\) and \(S=\{2\}\cup\{\ell:\ell\mid N\}\cup P\), and let \(G_{\Q,S}\) be the Galois group of the maximal extension of \(\Q\) unramified outside \(S\) and infinity. The quadratic characters define a homomorphism \(\chi_{\mathcal F}\colon G_\Q\to G_b\) by \(x\circ\chi_{\mathcal F}=\chi_{h_x}\). The cyclotomic and quadratic deformation of \(T_E\) is the free \(R=\Lambda[G_b]\)-module \(\mathbb T=T_E\otimes_{\Z_2}R\), equipped with the action
\[
 \sigma(v\otimes r)=\sigma v\otimes\chi_{\mathcal F}(\sigma)\,\bar\sigma^{-1}r,\qquad\bar\sigma\in\Gamma .
\]
The inverse cyclotomic factor is the Shapiro action. At a finite
layer, Galois acts on the group-ring coefficient by
\([\tau]\mapsto[\bar\sigma^{-1}\tau]\); passage to the inverse
limit gives \(\bar\sigma^{-1}\). Thus arithmetic Frobenius at an
odd prime \(q\in S\) acts on this factor by \(\gamma_q^{-1}\),
where \(\gamma_q\in\Gamma\) is its cyclotomic image. Specialization at \(x\) gives the cyclotomic deformation of \(T_E\chi_{h_x}\); further specialization at \(t=0\) gives \(T_2(E_x)\).

At \(2\), the global and real cohomology complexes can have
periodic tails. The following fibre removes their common tail and
places \(H^0(\RR,\mathbb T)\) in degree one.\footnote{The zero local
condition at the real place is imposed by a mapping fibre, as in
the determinant formalism of \cite{BF01} and
\cite[\S5.7]{Nekovar}.}

\begin{definition}\label{def:positive}
The \emph{modified global complex} of the family is
\[
 \Cplus=\Cone\bigl(\CS\longrightarrow\Rgam(\RR,\mathbb T)\bigr)[-1],\qquad \CS=\Rgam(G_{\Q,S},\mathbb T).
\]
The real complex is ordinary cohomology of \(\Gal(\C/\RR)\) in
non-negative degrees. For a number field \(K\) and coefficients \(W\),
write \(H^i_+(K,W)\) for the cohomology of the analogous cone, with
restriction to all real places and finite support induced by \(S\).
\end{definition}

Since every \(h_x\) is positive, complex conjugation \(\jmath\) acts on \(\mathbb T\) through \(T_E\). Hence \(H^0(\RR,\mathbb T)=T_E^{\jmath=1}\otimes R\) is free of rank one, and \(H^1(\RR,\mathbb T)\) is killed by \(2\). As \(H^0(G_{\Q,S},\mathbb T)=0\), the cone gives an exact sequence
\begin{equation}\label{eq:realseq}
 0\longrightarrow H^0(\RR,\mathbb T)\longrightarrow H^1(\Cplus)\longrightarrow H^1(\CS)\longrightarrow H^1(\RR,\mathbb T).
\end{equation}
The generator of \(H^0(\RR,\mathbb T)\) gives a class in
\(H^1(\Cplus)\). In Subsection~\ref{sec:kato} we take its exterior
product with a lift of twice the Beilinson--Kato element to define
the interpolation element (cf.~\cite[(3.1)]{OAPW},
\cite[(12.1)]{OAG}).

Put
\[
 \mathcal O=\widehat{\Lambda_{(2)}},\qquad
 \kLambda=\mathcal O/2\mathcal O
          =\operatorname{Frac}(\Lambda/2\Lambda).
\]
Thus \(\kLambda\) is the fraction field of the residual cyclotomic
Iwasawa algebra. Localizing at \((2)\) inverts \(t\);
\(\mathcal O\) is a discrete valuation ring with uniformizer \(2\). The local ring \(\mathcal O[G_b]\) is flat over \(R\) and has
residue field \(\kLambda\), with every \(g\in G_b\) mapping to one.

\begin{lemma}[{cf.~\cite[Lemma~3.5]{OAPW}}]\label{lem:perfect}
The complex \(\Cplus\) is a perfect complex of \(R\)-modules of Tor amplitude \([1,2]\). Its formation commutes with derived change of coefficients along the maps used below: evaluation at a character of \(G_b\), specialization at \(t=0\), reduction to the residue field, and the flat maps \(R\to\mathcal O[G_b]\) and \(R\to\LambdaQ[G_b]\).
\end{lemma}

\begin{proof}
For a finite \(2\)-primary \(G_{\Q,S}\)-module \(W\), restriction
to the real place is injective in degree zero and an isomorphism
in degrees at least three. In degree two, the terminal part of
Poitou--Tate duality is
\[
 H^2(G_{\Q,S},W)\longrightarrow
 \bigoplus_{v\in S}H^2(\Q_v,W)\oplus H^2(\RR,W)
 \longrightarrow H^0(G_{\Q,S},W^\vee(1))^\vee\longrightarrow0,
\]
with \({}^\vee\) Pontryagin duality. The map from
\(H^2(\Q_2,W)\) to the last group is surjective: its dual is the
injection of global into local invariants at \(2\).
Given \(c\in H^2(\RR,W)\), choose a dyadic class cancelling its
image in the last group. The resulting tuple lifts to global
cohomology. Thus restriction is surjective in degree two, and the
fibre has cohomology only in degrees one and two.

Apply the continuous-cochain finiteness and coefficient-change
construction to this fibre
\cite[\S4.2, 5.7.1.6--5.7.1.7]{Nekovar}.
The preceding bounds hold uniformly for every finite coefficient
module. They give finite Tor amplitude for the fibre by the
projection-formula argument of \cite[proof of 4.2.9]{Nekovar};
the finiteness criterion \cite[4.2.8]{Nekovar} then gives a perfect
complex. In particular,
\[
 \Cplus\tensL_R\F_2\simeq
 R\Gamma_+(\Q,E[2])\in D^{[1,2]}(\F_2).
\]
Here \(R\Gamma_+\) denotes the fibre of
Definition~\ref{def:positive}. A minimal finite free model over the
local ring \(R\) therefore has terms only in degrees one and two.
The same coefficient-change construction gives all the asserted
specializations and flat base changes
(cf.~\cite[Lemma~3.5]{OAPW}).\footnote{At \(p=2\), ordinary global
cohomology and the real complex need not be perfect separately; their
fibre \(\Cplus\) is.}
\end{proof}

\subsubsection{The residual rank calculation}
The following rank and specialization calculations are standard in
Iwasawa cohomology.\footnote{For their relation to weak Leopoldt and
the \(\mu\)-invariant, see \cite[\S1.3 and Appendix~B]{PR95}. The
proofs below retain the real-place terms needed at \(2\).}
\begin{definition}[The residual modified global complex]\label{def:residual-complex}
Put \(\Lambda_0=\F_2[[t]]\). For a finite
\(\F_2[G_{\Q,S}]\)-module \(W\), let \(\mathcal K_S(W)\) be the
modified global complex of Definition~\ref{def:positive} with coefficients
\(W\otimes_{\F_2}\Lambda_0\), where \(\sigma\) acts on the cyclotomic
factor by \(\bar\sigma^{-1}\). Its cohomology is concentrated in
degrees one and two and is finitely generated over \(\Lambda_0\).
\end{definition}

Write \(\Q_n=\Q(\zeta_{2^{n+2}})^+\),
\(d_n=[\Q_n:\Q]=2^n\), and \(\Gamma_n=\Gal(\Q_n/\Q)\).

\begin{lemma}\label{lem:residual-euler}
For every such \(W\),
\[
 \dim_{\F_2}H^1_+(\Q_n,W)-\dim_{\F_2}H^2_+(\Q_n,W)
 =d_n\dim_{\F_2}W.
\]
\end{lemma}
\begin{proof}
Put \(d=\dim W\), and write \(h^i\) for ordinary global cohomology
dimensions over \(\Q_n\). Tate's Euler characteristic formula
\cite[I, Theorem~5.1]{MilneADT} gives
\[
 h^0-h^1+h^2=-d_nd+
 \sum_{v\mid\infty}\dim H^0(\Q_{n,v},W).
\]
For a finite module acted on by complex conjugation, the two Tate
cohomology groups have the same order. Consequently
\(\dim H^1(\RR,W)=\dim H^2(\RR,W)\).
Restriction is injective in degree zero, surjective in degree two, and
an isomorphism in degrees at least three, as in Lemma~\ref{lem:perfect}.
The long exact sequence defining cohomology with the real-place condition therefore gives
\[
 \dim H^1_+-\dim H^2_+
 =-(h^0-h^1+h^2)+\sum_{v\mid\infty}\dim H^0(\Q_{n,v},W).
\]
Substitution proves the assertion.
\end{proof}

\begin{lemma}\label{lem:residual-rank}
Put \(r_i=\operatorname{rank}_{\Lambda_0}H^i(\mathcal K_S(W))\).
Then
\[
 r_1-r_2=\dim W,\qquad
 \dim H^2_+(\Q_n,W)=d_nr_2+O(1).
\]
In particular, \(H^2(\mathcal K_S(W))\) is torsion over \(\Lambda_0\)
if and only if \(\dim H^2_+(\Q_n,W)\) remains bounded. This
condition is unchanged when \(S\) is enlarged.
The finite-layer groups \(H^2_+(\Q_n,W)\) themselves can depend on \(S\).
\end{lemma}
\begin{proof}
Derived specialization modulo \(t^{d_n}\) computes the modified global complex
over \(\Q_n\). For \(H^i=H^i(\mathcal K_S(W))\), the
universal-coefficient sequence is
\[
0\longrightarrow H^i/t^{d_n}H^i\longrightarrow H^i_+(\Q_n,W)
\longrightarrow H^{i+1}[t^{d_n}]\longrightarrow0.
\]
A finitely generated \(\F_2[[t]]\)-module of rank \(r\) has quotient
modulo \(t^d\) of dimension \(rd+O(1)\), and its \(t^d\)-torsion has
bounded dimension. The formula for \(r_1-r_2\) follows from
Lemma~\ref{lem:residual-euler}. The torsion criterion follows as well.
At an added odd prime the residual localization complex becomes
acyclic over \(\kLambda\): its Frobenius differential has non-zero
determinant. Indeed, if its cyclotomic exponent is \(a\ne0\), then
\((1+t)^a\) is a nonconstant element of \(\F_2[[t]]\), hence
transcendental over \(\F_2\). For the finite-order Frobenius matrix
\(A\), the polynomial \(\det(1-AZ)\) is non-zero, so it cannot
vanish at \(Z=(1+t)^{-a}\). Thus enlarging
\(S\) does not change these generic cohomology groups.
\end{proof}

\begin{definition}\label{def:RF}
The curve \(E\) satisfies \textup{(RF)} if
\[
 H^2(\mathcal K_{\Sigma_0}(E[2]))\otimes_{\Lambda_0}\kLambda=0.
\]
By Lemma~\ref{lem:residual-rank}, this is equivalent to boundedness of
\(\dim H^2_+(\Q_n,E[2])\), and any fixed larger support may be used.\footnote{The
real-place condition matters here. When \(E(\RR)\) has two
components, its degree-two local Iwasawa cohomology is \(\Lambda/2\).
The unmodified global degree-two group surjects onto this module,
so its \(\mu\)-invariant is positive.}
\end{definition}

\begin{lemma}\label{lem:RF-mu}
Let \(\mathcal C_E^+\) be the complex of
Definition~\ref{def:positive} with coefficients
\(T_E\otimes_{\Z_2}\Lambda\), inverse cyclotomic action and support
\(\Sigma_0\). Put \(\mathcal H^2_{E,+}=H^2(\mathcal C_E^+)\).
This is a finitely generated torsion \(\Lambda\)-module, and
\[
 \RF{}\quad\Longleftrightarrow\quad
 \mu_\Lambda(\mathcal H^2_{E,+})=0,
 \qquad
 \mu_\Lambda(M)=\length_{\Lambda_{(2)}}M_{(2)}.
\]
\end{lemma}
\begin{proof}
Torsion follows from Kato's rational torsion theorem
\cite[Theorem~12.4(1)]{Kato} and localization; the cohomology at the real places in degrees one and two is
killed by \(2\). By Lemma~\ref{lem:perfect}, reduction
of \(\mathcal C_E^+\) modulo \(2\) is
\(\mathcal K_{\Sigma_0}(E[2])\). Its amplitude \([1,2]\) gives
\[
 H^2(\mathcal K_{\Sigma_0}(E[2]))
 \simeq\mathcal H^2_{E,+}/2\mathcal H^2_{E,+}.
\]
After tensoring with \(\kLambda\), the right side is the reduction
of \((\mathcal H^2_{E,+})_{(2)}\) modulo the maximal ideal of
\(\Lambda_{(2)}\). Nakayama's lemma identifies its vanishing with
\((\mathcal H^2_{E,+})_{(2)}=0\), equivalently \(\mu_\Lambda=0\).
\end{proof}

\subsubsection{Verification of the residual condition}
For constant coefficients the cohomology with the real-place condition of the layers is explicit. Let \(s_n\) be the number of finite places of \(\Q_n\)
above \(S\). These numbers are bounded for fixed \(S\): the dyadic prime
is totally ramified, and every odd prime has a non-zero cyclotomic
Frobenius and hence a decomposition group of finite index.

\begin{lemma}\label{lem:layers}
The class number of \(\Q_n\) is odd, and its units realize every real
signature. Its cohomology with the real-place condition with coefficients \(\F_2\) satisfies
\begin{equation}\label{eq:growth}
 \dim H^1_+(\Q_n,\F_2)=d_n-1+s_n,\qquad
 \dim H^2_+(\Q_n,\F_2)=s_n-1.
\end{equation}
\end{lemma}
\begin{proof}
In each quadratic step of the real cyclotomic tower only the dyadic
prime ramifies, and no archimedean place ramifies. The ambiguous class number formula, starting with
\(\Q\), makes the invariant class number odd. A non-zero finite
\(2\)-group with an involution has a non-zero fixed element, so the
whole class number is odd \cite{Iwasawa56}.

For \(n\ge1\), let \(\zeta\) have order \(2^{n+2}\). The unit
\(\varepsilon_n=1+\zeta+\zeta^{-1}\) has norm \(-1\): the minimal
polynomials of \(\zeta+\zeta^{-1}\) satisfy
\(f_1(Z)=Z^2-2\), \(f_{n+1}(Z)=f_n(Z^2-2)\), and \(f_n(-1)=-1\).
Its signature vector has odd augmentation in the local ring
\(\F_2[\Gamma_n]\), hence is a unit. Its conjugates span all real
signatures. For \(n=0\), use \(-1\)
(cf.~\cite[proof of Lemma~3.6]{OAPW}).

Kummer theory identifies global \(H^1\) with the \(S\)-units modulo
squares, of dimension \(d_n+s_n\). Restriction to real places is the
surjective signature map. The Brauer invariant sequence identifies the
kernel of restriction in degree two with a space of dimension
\(s_n-1\); restriction is surjective. In degree zero the cokernel has
dimension \(d_n-1\). The long exact sequence of the modified complex gives
\eqref{eq:growth}.
\end{proof}

Write \(d_2(A)=\dim_{\F_2}(A/2A)\) for a finite abelian group.
\begin{lemma}\label{lem:class-bound}
For a number field \(K\), let \(S_K\) contain the dyadic places
and let \(\Cl^+(K)\) be its narrow class group. Then
\begin{equation}\label{eq:positive-class-bound}
 \dim H^2_+(K,\F_2)\le2d_2\Cl^+(K)+\#S_K-1,
 \qquad d_2\Cl^+(K)\le d_2\Cl(K(i))+1.
\end{equation}
\end{lemma}
\begin{proof}
The long exact sequence of the modified complex bounds \(H^2_+\) by the
cokernel of signatures and the kernel of restriction in degree two.
The signature sequence bounds the former by \(d_2\Cl^+(K)\).
Kummer theory and the Brauer invariant sequence bound the latter
by \(d_2\Cl(K)+\#S_K-1\). This proves the first inequality.

For the second, let \(M/K\) be the elementary abelian extension
corresponding to \(\Cl^+(K)/2\). It is unramified at all finite
places. Its compositum with \(K(i)\) is everywhere unramified,
and the intersection \(M\cap K(i)\) removes at most one quadratic
direction. Ordinary class field theory gives the bound.
\end{proof}

\begin{proposition}\label{prop:RF-cases}
If \(E(\Q)[2]\ne0\), or if \(\Q(E[2])/\Q\) is cyclic cubic, then
\(E\) satisfies \RF{}.
\end{proposition}
\begin{proof}
In the first case, \(E[2]\) has two trivial successive quotients.
Lemmas~\ref{lem:layers} and~\ref{lem:residual-rank}, followed by the
long exact sequence, give the required vanishing
(cf.~\cite[Lemma~3.6]{OAPW}, \cite[Lemma~12.4]{OAG}).

In the cubic case put \(F=\Q(E[2])\). The permutation module on the
three non-zero two-torsion points is \(\F_2\oplus E[2]\). Shapiro's
lemma reduces the degree-two bound to constant coefficients over
\(F\Q_n\). Since \(F(i)/\Q\) is abelian, Ferrero--Washington gives
\(\mu=0\) for its cyclotomic tower \cite{FW}. Its inverse-limit
class group is therefore finitely generated over \(\Z_2\).
After finitely many layers a dyadic prime is totally ramified,
so the class-group norms are surjective and the finite-layer
\(2\)-ranks are bounded \cite[Chapter~13]{Washington}.
Lemma~\ref{lem:class-bound}, and boundedness of the number of places
above the fixed support, give the required degree-two cohomology bound.
\end{proof}

\begin{remark}\label{rem:fine-selmer}
Let \(\mathcal R_E\) be the Pontryagin dual of the fine Selmer group
of \(E[2^\infty]\) over \(\Q_\infty\), with zero local conditions at
the places above \(\Sigma_0\) and at the real places, as in
\cite{LY26}. By Poitou--Tate duality, \(\mathcal R_E\) is the kernel
of localization of global degree-two Iwasawa cohomology at these
places. As local degree-two Iwasawa cohomology at the finite places is
finitely generated over \(\Z_2\), the long exact sequence of the
modified complex gives
\[
 \mu_\Lambda(\mathcal H^2_{E,+})
 =\mu_\Lambda(\mathcal R_E)+\mu_\Lambda(\mathcal Q_E),
\]
where \(\mathcal Q_E\) is the cokernel of restriction from global
degree-one Iwasawa cohomology to \(H^1(\RR,T_E\otimes\Lambda)\), a
quotient of \(\Lambda/2\) that vanishes if \(E(\RR)\) is connected.
Thus \RF{} is the vanishing of the fine Selmer \(\mu\)-invariant
together with a condition at the real places. For \(E(\Q)[2]\ne0\),
the first vanishing is due to Lin--Yin \cite[Theorem~1.4]{LY26},
following \cite{CS05} and using Ferrero--Washington for
\(\Q(E[2])\), without hypothesis on the reduction at \(2\).
Proposition~\ref{prop:RF-cases} gives both vanishings for abelian
two-division fields, and so extends the first to cyclic cubic
two-division fields. For rational two-torsion its proof uses only the
odd class numbers and the units of the layers \(\Q_n\);
Ferrero--Washington enters only in the cyclic cubic case.
\end{remark}

\subsubsection{Freeness over the completed group ring}
\begin{proposition}[{cf.~\cite[\S3.3]{OAPW}, \cite[\S12.2]{OAG}}]\label{prop:free}
Assume \RF{}. For every admissible family,
\begin{equation}\label{eq:free}
 H^1(\Cplus)\otimes_R\mathcal O[G_b]\ \text{is free of rank }2,
 \qquad H^i(\Cplus)\otimes_R\mathcal O[G_b]=0\quad(i\ne1).
\end{equation}
\end{proposition}
\begin{proof}
The residue field of \(\mathcal O[G_b]\) is \(\kLambda\), with
every element of \(G_b\) mapping to one. The quadratic deformation
therefore reduces to \(E[2]\) with its cyclotomic action. Condition
\RF{} and Lemma~\ref{lem:residual-rank} give residual cohomology
\(\kLambda^2\) in degree one and zero elsewhere. A minimal finite free
model of the perfect complex over the local ring \(\mathcal O[G_b]\)
has zero differentials after reduction. It consequently has just one
term, of rank two in degree one. Flatness gives \eqref{eq:free}.
\end{proof}

\begin{remark}\label{rem:RF-scope}
Condition \RF{} depends only on \(E[2]\), hence is invariant under
quadratic twisting. The finite-layer bounds behind it may depend on
the support, but only their generic consequence enters
Proposition~\ref{prop:free}, so the determinant lattice has no
support-dependent denominator.
\end{remark}

\subsection{Beilinson--Kato elements and the determinant}\label{sec:kato}
Proposition~\ref{prop:integral} proves the uniform integral bound
in Theorem~\ref{thm:main}.

Let \(j\colon\operatorname{Spec}\Q\to
\operatorname{Spec}\Z[\tfrac12]\) be the natural map.
Fix an isomorphism \(\iota_f\colon V_f(1)\simeq V_E\), where \(V_f\)
is Kato's representation \cite[\S8.3]{Kato}. His full
cyclotomic Galois group is \(\Gal(\Q(\mu_{2^\infty})/\Q)\);
let \(\jmath\) be complex conjugation in this group. Then
\[
 \Z_2[[\Gal(\Q(\mu_{2^\infty})/\Q)]]\otimes_{\Z_2}\Q_2
 =e_+\LambdaQ\oplus e_-\LambdaQ,\qquad
 e_\pm=\frac{1\pm\jmath}{2}.
\]
After the twist \(V_f(1)\simeq V_E\), we take the component for
the extension \(\Q_\infty/\Q\), the union of the maximal real
subfields \(\Q_n\). The displayed idempotents are defined over
\(\Q_2\). If \(\sigma_a\) is the element corresponding to \(a\in\Z_2^\times\), the Tate twist substitutes \([\sigma_a]\mapsto a^{-1}[\sigma_a]\) in the coefficients (cf.~\cite[\S\S3.3--3.4]{OAPW}, \cite[\S\S12.2--12.3]{OAG}).

\subsubsection{The Euler-system theorem after inverting \texorpdfstring{\(2\)}{2}}\label{subsec:katodiv}

\begin{proposition}\label{prop:kato}
Let \(h\) be admissible and \(V_h=V_E\chi_h\). Then \(H^1_{j_*}(V_h)\) is free of rank one over \(\LambdaQ\), \(H^2_{j_*}(V_h)\) is torsion, and for every height-one prime \(\mathfrak p\) of \(\LambdaQ\),
\begin{equation}\label{eq:Kato}
 \ellp\bigl(H^2_{j_*}(V_h)\bigr)\le\ellp\bigl(H^1_{j_*}(V_h)/Z_h\bigr)+\ellp\bigl(H^2_{\mathrm{Iw}}(\Q_2,V_h)\bigr).
\end{equation}
Here \(H^i_{j_*}\) is rational Iwasawa cohomology of \(j_*V_h\) over \(\Z[\tfrac12]\), \(Z_h\subset H^1_{j_*}(V_h)\) is the \(\LambdaQ\)-module generated by Beilinson--Kato elements in the required sign component, \(H^2_{\mathrm{Iw}}(\Q_2,V_h)\) is local Iwasawa cohomology \cite[12.2]{Kato}, and \(\ellp(M)=\length_{(\LambdaQ)_{\mathfrak p}}M_{\mathfrak p}\).
\end{proposition}

\begin{proof}
We first suppose that \(E\) is non-CM. Fix \(h\), put
\(f_h=f\otimes\chi_h\), and let \(N_h\) be its conductor. Write
\[
 \widetilde\Gamma=\Gal(\Q(\mu_{2^\infty})/\Q),\qquad
 \widetilde\Lambda=\Z_2[[\widetilde\Gamma]],\qquad
 T_{f_h}=V_{\Z_2}(f_h)(1)\subset V_h.
\]
Thus \(T_{f_h}\) is the image in \(V_h\) of the cohomological
lattice \(V_{\Z_2}(f_h)(1)\) of the newform \(f_h\)
\cite[\S8.3]{Kato}. For this fixed \(h\), smoothed Beilinson--Kato elements
\[
 {}_{c,d}z_m(f_h)\in
 H^1\bigl(\Z[\mu_m,1/2],T_{f_h}\bigr)
\]
are the Euler systems of \cite[Example~13.3]{Kato}; the auxiliary
integers, symbol and excluded primes satisfy the conditions there.
Thus the Euler-system theorem is applied to the rank-two
\(\Z_2\)-lattice \(T_{f_h}\), separately for each \(h\).

In the notation of \cite[\S12.2]{Kato}, these systems bound
\[
 \mathbf H^2(T_{f_h})_0
 =\ker\bigl(\mathbf H^2(T_{f_h})
       \longrightarrow\mathbf H^2_{\mathrm{loc}}(T_{f_h})\bigr),
\]
where \(\mathbf H^i\) is ordinary Iwasawa cohomology over the full
cyclotomic tower. We invoke \cite[Theorem~12.5(3)]{Kato}, whose
proof in \S13 uses the Euler-system theorem
\cite[Theorem~13.4(2)]{Kato} at height-one primes
\(\mathfrak q\subset\widetilde\Lambda\) with \(2\notin\mathfrak q\).\footnote{The
rank and torsion assertions of \cite[Theorem~12.4(1)--(2)]{Kato}
and the divisibility of \cite[Theorem~12.5(3)]{Kato} apply at \(p=2\).
The integral refinements in parts~12.4(3) and~12.5(4) assume
\(p\ne2\).}
The hypotheses of Theorem~13.4(2) are generic
non-vanishing of the Euler system, one-dimensional complex-conjugation
eigenspaces, purity, irreducibility of the rational representation,
and the existence of \(\sigma\in G_{\Q(\mu_{2^\infty},\sqrt h)}\) with
\[
 \dim_{\Q_2}\ker(1-\sigma:V_h\longrightarrow V_h)=1.
\]
The first is supplied by \cite[Proposition~13.7]{Kato}; the middle
three hold for the representation of \(f_h\). For the last,
Serre's open image theorem \cite{Serre72} gives an open image in
\(\operatorname{SL}_2(\Z_2)\) for this finite-index subgroup of the
cyclotomic kernel, hence a non-trivial unipotent. Since
\(\chi_h(\sigma)=1\), it has the required action on \(V_h\)
\cite[(12.8.2), \S13.13]{Kato}.

To pass from smoothed systems to the normalized zeta module, Kato
uses their collection as the auxiliary choices vary.
\cite[Theorem~12.6]{Kato} identifies the generated cyclotomic module
with the module of normalized zeta elements up to a finite quotient
for the cohomological lattice. This quotient may depend on \(h\).
The comparison of \(T_{f_h}\) with \(T_2(E^h)\) may also
require a power \(2^{n_h}\) depending on \(h\). Both disappear after inverting \(2\). The pointwise localized
theorem therefore gives the rational divisibility, while the fixed
lattice of Proposition~\ref{prop:class} supplies the uniform
estimate at \((2)\).

Transporting the divisibility by the critical Tate twist and taking the component
for \(\Q_\infty/\Q\) gives \eqref{eq:Kato}. Its local degree-two
term at \(\Q_2\) persists after inverting \(2\); the fixed factor
\(\FE\) of Lemma~\ref{lem:fixedfactor} includes its characteristic
ideal. The rank formulas follow
from \cite[Theorem~12.4(1)--(2)]{Kato}.

If \(E\) is CM, use instead the rational main-conjecture equality
\cite[Theorem~2.6]{BT}, together with its rank formulas
\cite[Theorem~2.3]{BT}. This again gives the required assertion
over \(\LambdaQ\).
\end{proof}

\subsubsection{Local factors}\label{subsec:localfactors}

For an odd prime \(q\) let \(I_q\) be an inertia group, \(\mathrm{Fr}_q\) an arithmetic Frobenius and \(\gamma_q\in\Gamma\) its image, and put
\[
 \mathscr E_{q,h}(t)=\det\bigl(1-\gamma_q^{-1}\mathrm{Fr}_q\mid H^1(I_q,V_h)\bigr),\qquad
 P_q(Z)=1-\frac{a_q}{q}Z+\frac1qZ^2\quad(q\nmid N),
\]
where \(a_q=q+1-\#E(\F_q)\). The first is the determinant of the localization cone at \(q\) in the inverse convention.

\begin{lemma}\label{lem:euler}
Let \(h\) be admissible.
\begin{enumerate}
\item For \(q\in P\), \(\mathscr E_{q,h}=P_q(\chi_h(q)\gamma_q^{-1})\) if \(q\nmid h\), with value \(\#E^h(\F_q)/q\) at \(t=0\), and \(\mathscr E_{q,h}=1\) if \(q\mid h\).
\item For an odd prime \(\ell\mid N\), \(\mathscr E_{\ell,h}=\mathscr E_{\ell,1}\), and \(\mathscr E_{\ell,1}(0)\neq0\).
\end{enumerate}
\end{lemma}

\begin{proof}
If \(q\nmid h\), then \(H^1(I_q,V_h)=V_h(-1)\), on which \(\mathrm{Fr}_q\) has trace \(\chi_h(q)a_q/q\) and determinant \(1/q\); at \(t=0\) this gives \((q+1-\chi_h(q)a_q)/q=\#E^h(\F_q)/q\). If \(q\mid h\), inertia acts on \(V_h\) through a character of order two that is non-trivial, so \(H^1(I_q,V_h)=0\). For \(\ell\mid N\) the character \(\chi_h\) is trivial on \(G_{\Q_\ell}\). Finally \(H^1(I_\ell,V_E)\) has Frobenius coinvariants \(H^2(\Q_\ell,V_E)\), which is dual to \(H^0(\Q_\ell,V_E)=0\); so \(1\) is not an eigenvalue of \(\mathrm{Fr}_\ell\) and \(\mathscr E_{\ell,1}(0)\ne0\).
\end{proof}

\begin{lemma}\label{lem:fixedfactor}
There is \(\FE\in\Lambda\) with \(\FE(0)\ne0\), divisible over \(\LambdaQ\) by the product of the characteristic ideal of \(H^2_{\mathrm{Iw}}(\Q_2,V_E)\) and the elements \(\mathscr E_{\ell,1}\) for odd \(\ell\mid N\).
\end{lemma}

\begin{proof}
Local control gives \(H^2_{\mathrm{Iw}}(\Q_2,V_E)_\Gamma=H^2(\Q_2,V_E)=0\), so the characteristic ideal is prime to \(t\). Take the product of generators of the ideals in question, which do not vanish at \(t=0\) by Lemma~\ref{lem:euler}, and clear a power of \(2\).
\end{proof}

For admissible \(h\), \(V_h\simeq V_E\) as \(G_{\Q_2}\)-modules, so \(\FE\) is also divisible by the characteristic ideal of \(H^2_{\mathrm{Iw}}(\Q_2,V_h)\). Fix \(\FE\) once. Also fix an integral cohomological lattice for \(f\), a symbol \(\xi\in\mathrm{SL}_2(\Z)\) whose Betti vector \(\delta(f,1,\xi)^+\) is non-zero \cite[13.6]{Kato}, and a primitive generator of \(T_E^{\jmath=1}\). These choices constitute the normalization in Theorem~\ref{thm:main}.

\subsubsection{The class attached to the fixed form}

Put \(Q_R=\Frac(\Lambda)[G_b]\) and \(Q_{\mathcal O}=\Frac(\mathcal O)[G_b]\supset Q_R\).
For \(m=\prod_{q\in P}q\) and \(M_n=m2^{n+2}\), the family
character induces the integral coefficient map
\begin{equation}\label{eq:coefficient}
 \Z_2[\Gal(\Q(\mu_{m2^{n+2}})/\Q)]\longrightarrow\Z_2[\Gamma_n\times G_b],\qquad[\sigma]\longmapsto[\sigma|_{\Gamma_n}][\chi_{\mathcal F}(\sigma)].
\end{equation}
Indeed, every \(\chi_{h_x}\) has conductor dividing \(m\).

\begin{lemma}\label{lem:smoothing-units}
Choose \(c,d>1\) with
\(c,d\equiv1\pmod{12Nm}\). The images of the smoothing factors
\(c^2-c[\sigma_c]\) and \(d^2-d[\sigma_d]\) under the coefficient
map \eqref{eq:coefficient} are units in \(\mathcal O[G_b]\).
\end{lemma}
\begin{proof}
The image of \(\sigma_c\) in \(G_b\) is trivial. Its cyclotomic
image is \((1+t)^{\kappa_c}\) with \(\kappa_c\ne0\), since
\(c\in1+4\Z_2\) and \(c>1\). In the residue field of the local
ring \(\mathcal O[G_b]\),
\[
 c^2-c[\sigma_c]\longmapsto
 1-(1+t)^{\kappa_c}\in\kLambda^\times.
\]
The same holds for \(d\). These factors are not units in
\(\Lambda[G_b]\), since their values at \(t=0\) are
\(c(c-1),d(d-1)\in2\Z_2\).
\end{proof}

\begin{proposition}[{cf.~\cite[\S3.3]{OAPW}, \cite[\S12.2]{OAG}}]\label{prop:class}
There is \(c_{\mathrm{int}}\ge0\), depending only on \(E\) and the fixed choices, such that every admissible family has a class \(z_S=\FE z_S^{\mathrm{raw}}\in H^1(\CS)\otimes_RQ_R\), where \(z_S^{\mathrm{raw}}\) is the unsmoothed Beilinson--Kato element for the fixed form pushed to the family, with
\begin{equation}\label{eq:zbound}
 2^{c_{\mathrm{int}}}z_S\in\operatorname{im}\bigl(H^1(\CS)\otimes_R\mathcal O[G_b]\longrightarrow H^1(\CS)\otimes_RQ_{\mathcal O}\bigr).
\end{equation}
\end{proposition}

\begin{proof}
With the fixed comparison \(\iota_f\), put
\[
 T_f=\iota_f\bigl(V_{\Z_2}(f)(1)\bigr),\qquad
 c_{\mathrm{int}}=\min\{c\ge0:2^cT_f\subset T_E\}.
\]
Here \(V_{\Z_2}(f)\) is the image of integral cohomology in the
rational Hecke quotient for \(f\), as defined in
\cite[\S8.3]{Kato}. In particular both \(T_f\) and
\(c_{\mathrm{int}}\) are chosen before \(\mathcal F\).

Choose the smoothing integers \(c,d\) as in Lemma~\ref{lem:smoothing-units}.
The smoothed Beilinson--Kato elements have images
\[
 {}_{c,d}z_{M_n}(f,\xi)\in
 H^1\bigl(\Z[\mu_{M_n},\tfrac1S],T_f\bigr).
\]
Their integrality follows from the Chern character, cup products
and trace maps in \cite[(8.4.3), \S8.7--8.11]{Kato}; all maps
precede projection to the fixed quotient \(V_{\Z_2}(f)\).
The quadratic characters \(\chi_{h_x}\) factor through
\(\Gal(\Q(\mu_m)/\Q)\), since their conductors divide \(m\).
The construction is the composite
\begin{multline*}
 H^1\bigl(\Z[\mu_{M_n},\tfrac1S],T_f\bigr)\simeq H^1\bigl(G_{\Q,S},T_f\otimes\Z_2[\Gal(\Q(\mu_{M_n})/\Q)]\bigr)\\
 \longrightarrow H^1\bigl(G_{\Q,S},T_f\otimes\Z_2[\Gamma_n\times G_b]\bigr)\xrightarrow{\ 2^{c_{\mathrm{int}}}\ }H^1\bigl(G_{\Q,S},T_E\otimes\Z_2[\Gamma_n\times G_b]\bigr),
\end{multline*}
where the isomorphism is Shapiro's lemma and the middle map is
induced by \eqref{eq:coefficient}. The last map is multiplication by \(2^{c_{\mathrm{int}}}\) followed by the inclusion of lattices. Under inverse Shapiro the induced coefficient action is \(\tau[\sigma]=[\tau^{-1}\sigma]\), which is the action on \(\mathbb T\) fixed in Subsection~\ref{sec:positive}. The smoothed elements are compatible under corestriction along the cyclotomic tower, since \(2\in S\). Their inverse limit over \(n\) defines a class in \(H^1(\CS)\).
The coefficient map is integral for every tame conductor.\footnote{For
zeta morphisms in deformation families, compare \cite[\S3]{Nakamura}.}

In Kato's notation \(k=2\) and \(r=r'=1\); the smoothing factor is
\begin{equation}\label{eq:smoothing}
 (c^2-c[\sigma_c])(d^2-d[\sigma_d])
\end{equation}
by \cite[(4.2.4)]{Kato}. Choose \(c,d\) as in
Lemma~\ref{lem:smoothing-units}. These integers satisfy Kato's
coprimality and level conditions. The factors are non-zero elements
of \(\Lambda\), hence invertible in \(Q_R\), and are units in
\(\mathcal O[G_b]\). Dividing gives the unsmoothed class and
preserves the integral bound. Kato's compatibility relations
\cite[(4.2.4), \S8.11]{Kato} identify these quotients for different
\(c,d\); the unsmoothed class and its normalization are therefore
independent of those auxiliary integers. Multiplication by
\(\FE\in\Lambda\) preserves the bound as well. Thus \eqref{eq:zbound} holds with the fixed
\(c_{\mathrm{int}}\), independently of the number and choice of
primes in \(P\). As \(t\) is invertible in \(\mathcal O\), central
values are taken only after Proposition~\ref{prop:integral} places
\(2^{c_{\mathrm{int}}}\mathscr L_{\mathcal F}\) in \(\Lambda[G_b]\).
If \(T_f\subset T_E\), then \(c_{\mathrm{int}}=0\).
\end{proof}

\subsubsection{Comparison with the primitive zeta module}

\begin{lemma}[{cf.~\cite[(3.8)]{OAPW}}]\label{lem:comparison}
Let \(x\) be an index and \(h=h_x\). There are \(a_h\in\Q_2^\times\), a group-like unit \(\varepsilon_h\in\Lambda^\times\) with \(\varepsilon_h(0)=1\), a fixed \(\mathscr E_{\mathrm{bad}}\in\LambdaQ\) with \(\mathscr E_{\mathrm{bad}}(0)\neq0\), and a generator \(z_h\) of \(Z_h\), such that
\begin{equation}\label{eq:comparison}
 z^{\mathrm{raw}}_{S,x}=a_h\,\varepsilon_h(t)\,\mathscr E_{\mathrm{bad}}(t)\prod_{\substack{q\in P\\ q\nmid h}}\mathscr E_{q,h}(t)\ z_h .
\end{equation}
\end{lemma}

\begin{proof}
The \(+\)-part of the Betti realization of \(f\otimes\chi_h\) is a
line. Kato's zeta map \cite[Theorem~12.5(1)--(2)]{Kato} sends any
non-zero vector of it to a generator \(z_h\) of \(Z_h\) over \(\LambdaQ\).
The even twisting isomorphism preserves this line. Its comparison with
the fixed Betti vector gives \(a_h\in\Q_2^\times\).

The character \(\chi_h\) has conductor \(h\mid m\). The \(\chi_h\)-component of the class of tame conductor \(m\) is obtained from that of tame conductor \(h\) by Kato's norm relations \cite[Proposition~8.12]{Kato}, which for \(k=2\), \(r=1\) multiply by \(1-\tfrac{a_q}{q}[\sigma_q]^{-1}+\tfrac1q[\sigma_q]^{-2}\) for each prime \(q\mid m/h\). The map \eqref{eq:coefficient} followed by evaluation at \(x\) sends \([\sigma_q]^{-1}\) to \(\chi_h(q)\gamma_q^{-1}\), so this factor becomes \(\mathscr E_{q,h}\) by Lemma~\ref{lem:euler}. For a character \(\vartheta\) of conductor \(2^n\), the Gauss sums satisfy \(\tau(\chi_h\vartheta)=\chi_h(2^n)\vartheta(h)\tau(\chi_h)\tau(\vartheta)\), and \(\chi_h(2)=1\) by admissibility; the factor \(\vartheta(h)\) is group-like and gives \(\varepsilon_h\). The omitted Euler factors at the odd primes of \(N\) give \(\mathscr E_{\mathrm{bad}}\) \cite[13.9--13.10]{Kato}.

Both sides of \eqref{eq:comparison} lie in
\(H^1(\CS\tensL_{R,x}\LambdaQ)\), the free rank-one
\(\LambdaQ\)-module of Iwasawa cohomology of \(V_h\) with support
\(S\). Outside finitely many finite-order characters \(\rho\),
its specialization is \(H^1(G_{\Q,S},V_h\rho^{-1})\) over the
coefficient field of \(\rho\), of dimension one. Whenever
\(L(f,\chi_h\rho,1)\ne0\), reciprocity gives a non-zero dual
exponential on that line, so it is injective and identifies the two
specialized classes. Rohrlich's theorem supplies infinitely many
such \(\rho\) \cite{Rohrlich}. After clearing denominators, their
difference is an element of \(\Lambda\) vanishing at all these
characters. A non-zero element has only finitely many such zeros by
Weierstrass preparation, proving the identity.
\end{proof}

The factor \(a_h\) is a unit in \(\LambdaQ\) and does not affect divisors away from \((2)\). The integral bound comes instead from the fixed lattice in Proposition~\ref{prop:class}.

\subsubsection{The determinant}

Let \(\eta_\infty\in H^1(\Cplus)\) be the image under \eqref{eq:realseq} of the fixed generator of \(T_E^{\jmath=1}\otimes1\). Over \(Q_R\), Proposition~\ref{prop:kato} gives global cohomology of rank one in degree one and zero in degree two at every character, and \eqref{eq:realseq} adds the real line. Hence \(H^1(\Cplus)\otimes_RQ_R\) is free of rank two, \(H^2(\Cplus)\otimes_RQ_R=0\), and the inverse determinant \(\det_R^{-1}\Cplus\), a free \(R\)-module of rank one, is identified over \(Q_R\) with \(\bigwedge^2H^1(\Cplus)\otimes_RQ_R\).

By Lemma~\ref{lem:comparison} and the decomposition \(\LambdaQ[G_b]=\prod_{x\in X_b}\LambdaQ\), the class \(z_S\) lies in \(H^1(\CS)\otimes_R\LambdaQ[G_b]\). Since \(H^1(\RR,\mathbb T)\) is killed by \(2\), the sequence \eqref{eq:realseq} becomes surjective onto global \(H^1\) after inverting \(2\); choose a lift
\[
 z_S^+\in H^1(\Cplus)\otimes_R\LambdaQ[G_b]
\]
of \(2z_S\). Any two such lifts differ by a multiple of \(\eta_\infty\). The lift has a central specialization at every character.

\begin{definition}[The interpolation element]\label{def:interpolation-element}
Choose an \(R\)-basis \(\mathbf e_S\) of \(\det_R^{-1}\Cplus\).
The \emph{interpolation element} of \(\mathcal F\), relative to this
basis and the fixed zeta normalization, is the unique
\(\mathscr L_{\mathcal F}\in Q_R\) satisfying
\begin{equation}\label{eq:udef}
 z_S^+\wedge\eta_\infty=\mathscr L_{\mathcal F}\,\mathbf e_S.
\end{equation}
\end{definition}

\begin{remark}\label{rem:interpolation-normalization}
The element is independent of the lift \(z_S^+\), since any two
lifts differ by a multiple of \(\eta_\infty\). Changing \(\mathbf e_S\)
multiplies it by an \(R\)-unit. Thus the vanishing and valuation of
\(\mathscr L_{\mathcal F,x}(0)\) are independent of that basis.
\end{remark}

\begin{proposition}[{cf.~\cite[Lemmas~3.3, 3.6]{OAPW}}]\label{prop:integral}
The element satisfies \(\mathscr L_{\mathcal F}\in\LambdaQ[G_b]\). If \RF{} holds, then \(2^{c_{\mathrm{int}}}\mathscr L_{\mathcal F}\in\Lambda[G_b]\).
\end{proposition}

\begin{proof}
Since \(2\) is invertible in \(\LambdaQ\), we have \(\LambdaQ[G_b]=\prod_{x\in X_b}\LambdaQ\), and we may argue one character at a time. Fix \(x\), put \(h=h_x\) and \(C_{S,h}=\CS\tensL_{R,x}\LambdaQ\).
Its groups \(H^i(C_{S,h})\) are Iwasawa cohomology of \(V_h\)
with support \(S\) over \(\Q_\infty/\Q\).
Apply Proposition~\ref{prop:kato} to this twist. The localization triangle from primitive cohomology to support \(S\) has, at each odd \(q\in S\), the cone
\begin{equation}\label{eq:cone}
 \bigl[H^1(I_q,V_h)\otimes\LambdaQ\xrightarrow{\ 1-\gamma_q^{-1}\mathrm{Fr}_q\ }H^1(I_q,V_h)\otimes\LambdaQ\bigr]\quad\text{in degrees }1,2,
\end{equation}
whose differential has determinant \(\mathscr E_{q,h}\neq0\), as \(\gamma_q\) has non-zero cyclotomic exponent. Its degree-one cohomology vanishes and its degree-two cohomology has length \(\operatorname{ord}_{\mathfrak p}(\mathscr E_{q,h})\) at \(\mathfrak p\). Hence the global \(H^1\) line is unchanged, and the length of global \(H^2\) grows by at most \(\sum_q\operatorname{ord}_{\mathfrak p}(\mathscr E_{q,h})\). By \eqref{eq:comparison}, the inactive primes add exactly \(\sum_{q\nmid h}\operatorname{ord}_{\mathfrak p}(\mathscr E_{q,h})\) to the index of \(z_S\), active primes contribute nothing on either side, and the factor \(\FE\) covers the local term in \eqref{eq:Kato} and the cones at the odd primes of \(N\). These are comparisons of divisors over \(\LambdaQ\).
The local term at \(2\) is fixed by admissibility, and the varying
Euler factors at the new primes are matched with their localization
cones.
With \eqref{eq:Kato} this gives
\begin{equation}\label{eq:divisibility}
 \ellp\bigl(H^2(C_{S,h})\bigr)\le\ellp\bigl(H^1(C_{S,h})/\LambdaQ z_{S,x}\bigr)
\end{equation}
at every height-one prime \(\mathfrak p\) of \(\LambdaQ\) (cf.~\cite[(3.9)]{OAPW}). The real line adds no torsion over \(\LambdaQ\), so the divisor of \(\mathscr L_{\mathcal F,x}\) at \(\mathfrak p\) is the right side minus the left side of \eqref{eq:divisibility}, which is non-negative. As \(\LambdaQ\) is regular, \(\mathscr L_{\mathcal F,x}\in\LambdaQ\). Hence \(\mathscr L_{\mathcal F}\in\LambdaQ[G_b]\).

At \((2)\), assume \RF{}, equivalently
\(\mu_\Lambda(\mathcal H^2_{E,+})=0\) by Lemma~\ref{lem:RF-mu}.
Proposition~\ref{prop:free} removes the degree-two cohomology and
identifies the inverse determinant lattice with the exterior square
of a free rank-two module. Let \(c=c_{\mathrm{int}}\). By \eqref{eq:zbound}, \(2^cz_S\) is the image of a class in \(H^1(\CS)\otimes_R\mathcal O[G_b]\). As \(\mathcal O[G_b]\) is flat over \(R\), the sequence \eqref{eq:realseq} remains exact after tensoring with it, and \(H^1(\RR,\mathbb T)\otimes_R\mathcal O[G_b]\) is killed by \(2\); so \(2^{c+1}z_S\) has a lift \(w_S\in H^1(\Cplus)\otimes_R\mathcal O[G_b]\). In \(H^1(\Cplus)\otimes_RQ_{\mathcal O}\) the classes \(w_S\) and \(2^cz_S^+\) both lift \(2^{c+1}z_S\), so they differ by a multiple of \(\eta_\infty\), and
\[
 w_S\wedge\eta_\infty=2^cz_S^+\wedge\eta_\infty=2^c\mathscr L_{\mathcal F}\,\mathbf e_S .
\]
By Proposition~\ref{prop:free}, \(H^1(\Cplus)\otimes_R\mathcal O[G_b]\) is free of rank two and its exterior square is \(\det_R^{-1}\Cplus\otimes_R\mathcal O[G_b]\), with basis \(\mathbf e_S\). Hence \(2^c\mathscr L_{\mathcal F}\in\mathcal O[G_b]\).

Inside \(Q_{\mathcal O}\) one has, coefficientwise, \(\mathcal O[G_b]\cap\LambdaQ[G_b]=\Lambda[G_b]\), since \(2^{-j}f(t)\in\mathcal O\) for \(f\in\Lambda\) means that every coefficient of \(f\) is divisible by \(2^j\). As \(\mathscr L_{\mathcal F}\in\LambdaQ[G_b]\) and \(2^c\mathscr L_{\mathcal F}\in\mathcal O[G_b]\), this gives \(2^c\mathscr L_{\mathcal F}\in\Lambda[G_b]\) (cf.~\cite[(3.11)--(3.12) and Lemma~3.3]{OAPW}, \cite[(12.13)]{OAG}).
\end{proof}

\begin{lemma}\label{lem:central-order}
Put \(A=(\LambdaQ)_{(t)}\). For every admissible family,
\[
 \operatorname{ord}_{t=0}\mathscr L_{\mathcal F,0}
 =\length_A\bigl((\mathcal H_E^1/Z_E)\otimes_{\LambdaQ}A\bigr)
  -\length_A\bigl(\mathcal H_E^2\otimes_{\LambdaQ}A\bigr).
\]
\end{lemma}
\begin{proof}
At the trivial character, \(h_0=1\). The factors \(\FE\),
\(\mathscr E_{\mathrm{bad}}\) and \(\mathscr E_{q,1}\) are units
in \(A\) (Lemmas~\ref{lem:fixedfactor}, \ref{lem:comparison}
and~\ref{lem:euler}). Thus the localization cones in
\eqref{eq:cone} are acyclic over \(A\), and \(z_{S,0}\) generates
the same \(A\)-submodule as the primitive zeta element.
The real-place term adds a free line in degree one. The divisor
calculation in Proposition~\ref{prop:integral}, now over the
discrete valuation ring \(A\), is therefore exactly the difference
of the two displayed lengths.
\end{proof}

\begin{remark}\label{rem:uniformity-separated}
The two estimates have different sources and quantifiers:
\[
 \begin{aligned}
 &\mathscr L_{\mathcal F}\in\LambdaQ[G_b]
 &&\text{for each admissible }\mathcal F,\\
 &2^{c_{\mathrm{int}}}\mathscr L_{\mathcal F}\in\mathcal O[G_b]
 &&\text{for all admissible }\mathcal F,
 \end{aligned}
\]
where the same \(c_{\mathrm{int}}\) works in the second line.
The first follows from the pointwise rational divisibilities of
Proposition~\ref{prop:kato}, after the twist-dependent \(2\)-primary
errors have disappeared. The second follows from the fixed-form construction of
Proposition~\ref{prop:class} and the exact freeness of
Proposition~\ref{prop:free}. Intersecting the two rings gives the
bound over \(\Lambda[G_b]\) with this same exponent for all
admissible families.
\end{remark}

\begin{remark}\label{rem:kato-equivariant}
The coefficient maps in Proposition~\ref{prop:class} are integral
and the Beilinson--Kato elements satisfy norm relations as the
cyclotomic conductor varies. Together these suggest a zeta
morphism for quadratic deformations, with the Euler factors
retained when the ramification support changes.
Nakamura constructs zeta morphisms on universal deformation
spaces \cite{Nakamura}, under different residual
hypotheses.\footnote{These include \(p\ge5\) and absolute
irreducibility on \(G_{\Q(\mu_p)}\)
\cite[Theorem~1.1]{Nakamura}; the cyclic cubic representation here
splits over \(\F_4\).}
Kato's equivariant main conjecture predicts a normalized zeta
element generating the integral determinant line
\cite{Kato93,BF01}. Comparison with \eqref{eq:udef} includes the
real-place condition, the Betti vector and the factor \(\FE\)
\cite[Remark~4.12]{Kato93}.
\end{remark}

\subsection{The central \texorpdfstring{\(L\)}{L}-values}\label{sec:central}

Theorem~\ref{thm:central} computes the central specialization from
Beilinson--Kato reciprocity and Poitou--Tate duality
(cf.~\cite[Proposition~3.7]{OAPW}, \cite[Proposition~12.3]{OAG}).
The global index and the inactive-prime Euler factors cancel;
Corollary~\ref{cor:bsd} identifies the remaining valuation with
the \(2\)-adic BSD defect plus a fixed constant.

Throughout this section \(x\) has Selmer corank zero. Write \(h=h_x\), \(T_h=T_2(E^h)\), \(V_h=T_h\otimes\Q_2\), and take lengths over \(\Z_2\).

\subsubsection{Periods, torsion and local conditions}

\begin{lemma}\label{lem:period}
For admissible \(h\), \(\Delta^{\min}_{E^h}=h^6\Delta^{\min}_E\). If \(\phi_h\colon E^h\to E\) is a twisting isomorphism, defined over \(\Q(\sqrt h)\subset\RR\), then \(\omega_{E^h}=\pm h^{-1/2}\phi_h^*\omega_E\) and \(\Omega^+_{E^h}=\Omega^+_E/\sqrt h\).
\end{lemma}

\begin{proof}
At \(\ell\mid2N\) the curves are isomorphic over \(\Q_\ell\), and \(\sqrt h\in\Z_\ell^\times\), so the minimal discriminants have equal valuation and the differentials correspond up to units. At \(q\mid h\) the twist of good reduction has type \(I_0^*\) and discriminant valuation six; at every other prime the twist is unramified and the valuation is unchanged. The discriminant keeps its sign since \(h>0\). Thus the twisted differential \(h^{-1/2}\phi_h^*\omega_E\) is minimal everywhere. As \(\phi_h\) is defined over \(\RR\), it preserves the invariant homology, which gives the period.
\end{proof}

\begin{lemma}\label{lem:torsion}
Let \(a_E=\dim_{\F_2}E(\Q)[2]\) and \(r(h)=\val(\#E^h(\Q)[2^\infty])\). If \(h\ne1\) is admissible, then \(E^h(\Q)[2^\infty]=E(\Q)[2]\) under the identification \(E^h[2]=E[2]\), and \(r(h)=a_E\).
\end{lemma}

\begin{proof}
Let \(q\mid h\). Inertia at \(q\) acts trivially on \(E[2^\infty]\) and through \(\chi_h\), by \(-1\), on \(E^h[2^\infty]\). A rational point \(Q\) of \(2\)-power order on \(E^h\) therefore satisfies \(Q=-Q\). Quadratic twisting does not change the Galois module \(E[2]\).
\end{proof}

For a place \(v\) let \(K_v\subset H^1(\Q_v,T_h)\) be the image of the \(2\)-adic completion of \(E^h(\Q_v)\) under the Kummer map, and put \(L_2=H^1(\Q_2,T_h)/K_2\). Let \(\epsilon_\infty=\val(\#\pi_0(E(\RR)))\).

\begin{lemma}\label{lem:local}
At odd primes and at the real place \(K_v=H^1(\Q_v,T_h)\). The module \(L_2\) is free of rank one, and \(\length H^1(\RR,T_h)=\length H^2(\RR,T_h)=\epsilon_\infty\). The integers
\begin{equation}\label{eq:fixedlocal}
 \rho_E=\sum_{\ell\mid2N}\length E(\Q_\ell)[2^\infty],\qquad
 \lambda_E\ \text{with}\ \exp_2^*(L_2)=2^{\lambda_E}\Z_2\,\omega_{E^h},
\end{equation}
do not depend on the admissible parameter \(h\).
\end{lemma}

\begin{proof}
The quotient \(H^1(\Q_v,T_h)/K_v\) is the \(2\)-adic Tate module of the local Weil--Ch\^atelet group, which is zero at odd primes and at the real place. At \(2\), local duality identifies \(L_2\) with the dual of the torsion-free part of the completed local points, which has rank one. The Kummer map identifies \(H^1(\RR,T_h)\) with \(\pi_0(E^h(\RR))\otimes\Z_2\), of length \(\epsilon_\infty\) because \(E^h\simeq E\) over \(\RR\); the Herbrand quotient of \(T_h\) is one, which gives the same length in degree two. The local isomorphisms at \(\ell\mid2N\) fix \(\rho_E\). The formal logarithm fixes the comparison between the dual exponential and the differential, and since \(\sqrt h\in\Z_2^\times\) the comparison is unchanged under the local isomorphism at \(2\).
\end{proof}

\subsubsection{Reciprocity and duality}

\begin{lemma}\label{lem:reciprocity}
There is \(\upsilon_E\in\Q_2^\times\), depending only on \(E\) and the fixed choices, such that
\[
 \exp_2^*z_{S,x}(0)=\upsilon_E\,\varepsilon\,\Lalg(E^h,1)\prod_{\substack{q\in P\\ q\nmid h}}\frac{\#E^h(\F_q)}{q}\ \omega_{E^h},\qquad \varepsilon\in\Z_2^\times .
\]
In particular a vanishing central \(L\)-value gives zero dual exponential, and with \(\alpha_E=\val(\upsilon_E)\), for \(L(E^h,1)\ne0\),
\begin{equation}\label{eq:exactreciprocity}
 v_{2,\omega_{E^h}}\bigl(\exp_2^*z_{S,x}(0)\bigr)=\val\bigl(\Lalg(E^h,1)\bigr)+\sum_{\substack{q\in P\\ q\nmid h}}\val\bigl(\#E^h(\F_q)\bigr)+\alpha_E ,
\end{equation}
where \(v_{2,\omega}\) is the valuation of the coefficient of \(\omega\).
\end{lemma}

\begin{proof}
Kato's character-sum formula \cite[Theorems~6.6 and~9.7, 12.5(1)]{Kato}, applied with the fixed Betti vector \(\delta(f,1,\xi)^+\), gives the displayed formula (cf.~\cite[(3.4)]{OAPW}). Its factors are as follows. The period vector, \(\FE(0)\), \(\mathscr E_{\mathrm{bad}}(0)\) and the factor at \(2\) depend only on \(E\). The quadratic Gauss sum has square \(h\), hence is a \(2\)-adic unit, as is the comparison of differentials at \(2\) in Lemma~\ref{lem:period}; together with \(\sqrt h\,L(E^h,1)/\Omega^+_E=\Lalg(E^h,1)\) this accounts for the central \(L\)-value. The character and group-like factors are units, \(\chi_h(2)=1\), and the inactive primes contribute \(\#E^h(\F_q)/q\) by Lemma~\ref{lem:euler}. Since \(\upsilon_E\) comes from the normalization of the reciprocity law, it is defined even when \(L(E,1)=0\).
\end{proof}

\begin{lemma}\label{lem:PT}
Rational global \(H^1(G_{\Q,S},V_h)\) maps isomorphically to \(L_2\otimes\Q_2\), and \(H^2(G_{\Q,S},T_h)\) is finite. Let \(\mathcal J_h\subset L_2\) be the image of \(H^1(G_{\Q,S},T_h)\) and \(i_h=\length(L_2/\mathcal J_h)\). Then \(H^1(G_{\Q,S},T_h)_{\mathrm{tors}}\simeq E^h(\Q)[2^\infty]\), of length \(r(h)\), and
\begin{equation}\label{eq:exactPT}
 \length H^2(G_{\Q,S},T_h)=s_{\Sha}(h)-i_h+\sum_{\ell\in S}\length E^h(\Q_\ell)[2^\infty]+\epsilon_\infty-r(h).
\end{equation}
\end{lemma}

\begin{proof}
Integral Poitou--Tate duality with the Kummer conditions and the
real term retained \cite[Theorem~6.3.4, \S6.9.1, 5.7.1.6]{Nekovar}
gives the exact sequence\footnote{Compare the duality deduction in
\cite[\S14.9, (14.9.1)--(14.9.3)]{Kato}, with the real-place term
included here.}
\begin{multline*}
0\longrightarrow L_2/\mathcal J_h\longrightarrow\Sel_{2^\infty}(E^h)^\vee\longrightarrow H^2(G_{\Q,S},T_h)\\
\longrightarrow\bigoplus_{\ell\in S}H^2(\Q_\ell,T_h)\oplus H^2(\RR,T_h)\longrightarrow E^h(\Q)[2^\infty]^\vee\longrightarrow0,
\end{multline*}
in which the real term is Tate cohomology, equal to ordinary cohomology in positive degrees. The rational sequence gives the first assertions, as the Selmer group is finite. The Selmer group has length \(s_{\Sha}(h)\), local duality identifies \(H^2(\Q_\ell,T_h)\) with the dual of \(E^h(\Q_\ell)[2^\infty]\), and Lemma~\ref{lem:local} gives the real length. Additivity gives \eqref{eq:exactPT}.

The sequence \(0\to T_h\to V_h\to E^h[2^\infty]\to0\), with \(H^0(\Q,V_h)=0\), identifies the torsion of global \(H^1\).
\end{proof}

\subsubsection{The central formula}

\begin{lemma}\label{lem:detvaluation}
The value \(\mathscr L_{\mathcal F,x}(0)\) is the coefficient of \(z^+_{S,x}(0)\wedge\eta_\infty\) in the basis specialized from \(\mathbf e_S\). If \(L(E^h,1)\ne0\), then
\begin{equation}\label{eq:detexact}
 \val\bigl(\mathscr L_{\mathcal F,x}(0)\bigr)=v_{2,\omega_{E^h}}\bigl(\exp_2^*z_{S,x}(0)\bigr)-\lambda_E-i_h+1-\length H^2(G_{\Q,S},T_h)+r(h).
\end{equation}
\end{lemma}

\begin{proof}
Let \(P^+\) be the specialization of \(\Cplus\) at \(x\) and \(t=0\), a perfect \(\Z_2\)-complex computing cohomology with the real-place condition of \(T_h\) (Lemma~\ref{lem:perfect}). Since the \(2^\infty\)-Selmer group of \(E^h\) is finite, Lemma~\ref{lem:PT} shows that rationally \(H^1(P^+)\) has dimension two and \(H^2(P^+)=0\), as at the generic point; this holds whether or not \(L(E^h,1)\) vanishes. Hence the identification of the inverse determinant with \(\bigwedge^2H^1\) specializes, which is the first assertion. Suppose now that \(L(E^h,1)\ne0\).

For a perfect complex \(Q\) over \(\Z_2\) one has \(\det^{-1}Q=\bigotimes_i(\det H^i(Q))^{(-1)^{i+1}}\), the determinant of a finite module being the corresponding fractional ideal; in the rational trivialization, torsion of length \(m\) in degree two multiplies the lattice by \(2^m\) and torsion of length \(m\) in degree one by \(2^{-m}\):
\[
 \text{inverse determinant lattice}=2^{\length H^2-\length H^1_{\mathrm{tors}}}\cdot(\text{lattice of the free part of }H^1),
\]
so the coefficient loses \(\length H^2\) and gains \(\length H^1_{\mathrm{tors}}\).

Let \(P_S\) and \(P_\RR\) be the global and real complexes of \(T_h\). Restriction is an isomorphism in degrees at least three and surjective in degree two (Lemma~\ref{lem:perfect}), so \(P^+\) is the fiber of \(\tau_{\le2}P_S\to\tau_{\le2}P_\RR\), and
\[
 \det{}^{-1}P^+=\det{}^{-1}(\tau_{\le2}P_S)\otimes\bigl(\det{}^{-1}(\tau_{\le2}P_\RR)\bigr)^{-1}.
\]
In \(\tau_{\le2}P_\RR\) the torsion groups \(H^1\) and \(H^2\) have the same length \(\epsilon_\infty\) and cancel, while \(H^0=T_h^{\jmath=1}\) is free of rank one with generator mapping to \(\eta_\infty\). In \(\tau_{\le2}P_S\), \(H^1\) is free of rank one plus torsion of length \(r(h)\), and \(H^2\) is finite. Hence \(\val(\mathscr L_{\mathcal F,x}(0))\) is the valuation of the coefficient of \(2z_{S,x}(0)\) in \(H^1(G_{\Q,S},T_h)/\mathrm{tors}\), plus \(r(h)-\length H^2(G_{\Q,S},T_h)\). By Lemma~\ref{lem:PT} that free line maps isomorphically onto \(\mathcal J_h=2^{i_h}L_2\), and \(\exp_2^*(L_2)=2^{\lambda_E}\Z_2\omega_{E^h}\); so the coefficient of \(z_{S,x}(0)\) has valuation \(v_{2,\omega_{E^h}}(\exp_2^*z_{S,x}(0))-\lambda_E-i_h\), and the factor \(2\) adds one.
\end{proof}

The constants \(\alpha_E,\lambda_E,\rho_E,\epsilon_\infty\)
are fixed under admissible twisting. The central \(L\)-value,
\(s_{\Sha}(h)\) and \(\tau(h)\) vary; Lemma~\ref{lem:torsion}
also makes \(r(h)\) fixed for \(h\ne1\).

\begin{theorem}\label{thm:central}
At an index of Selmer corank zero, \(\mathscr L_{\mathcal F,x}(0)=0\) if and only if \(L(E^{h_x},1)=0\). If \(L(E^{h_x},1)\ne0\), then
\begin{equation}\label{eq:central-r}
 \val\bigl(\mathscr L_{\mathcal F,x}(0)\bigr)=D(h_x)+2r(h_x)+\alpha_E+1-\lambda_E-\rho_E-\epsilon_\infty .
\end{equation}
\end{theorem}

\begin{proof}
Suppose \(L(E^h,1)\ne0\), and put
\(n_q=\val(\#E^h(\F_q))\) for each inactive prime.
At an active prime, inertia invariants of \(E^h[2^\infty]\) are
\(E[2]\), whose Frobenius invariants have length \(t_q\).
At an inactive prime, reduction identifies the local torsion with
\(E^h(\F_q)[2^\infty]\). Thus
\[
 \sum_{\ell\in S}\length E^h(\Q_\ell)[2^\infty]
 =\rho_E+\tau(h)+\sum_{\substack{q\in P\\q\nmid h}}n_q.
\]
Substituting \eqref{eq:exactreciprocity} and \eqref{eq:exactPT}
into \eqref{eq:detexact} gives
\begin{align*}
 \val(\mathscr L_{\mathcal F,x}(0))
 &=\val(\Lalg(E^h,1))+\sum_{q\nmid h}n_q+\alpha_E
     -\lambda_E-i_h+1\\
 &\quad-\left(s_{\Sha}(h)-i_h+\rho_E+\tau(h)
          +\sum_{q\nmid h}n_q+\epsilon_\infty-r(h)\right)+r(h)\\
 &=D(h)+2r(h)+\alpha_E+1-\lambda_E-\rho_E-\epsilon_\infty.
\end{align*}
Here both sums run over inactive primes in \(P\).
The two occurrences of \(r(h)\) come from duality and the torsion
of global \(H^1\). Every new-prime contribution has either cancelled
or entered \(c(h)\); the remaining constant is fixed.

If \(L(E^h,1)=0\), Lemma~\ref{lem:reciprocity} gives \(\exp_2^*z_{S,x}(0)=0\). The specialized class then satisfies the finite condition at \(2\), and at the other places it does so automatically by Lemma~\ref{lem:local}; so it lies in the rational Selmer group, which is zero. Hence \(z^+_{S,x}(0)\) is rationally a multiple of \(\eta_\infty\), its wedge with \(\eta_\infty\) vanishes, and \(\mathscr L_{\mathcal F,x}(0)=0\) by Lemma~\ref{lem:detvaluation}. Conversely, a non-zero central \(L\)-value gives a finite valuation in \eqref{eq:detexact}.
\end{proof}

\begin{proof}[Proof of Theorem~\ref{thm:main}]
Proposition~\ref{prop:integral} gives \(\mathscr L_{\mathcal F}\in\LambdaQ[G_b]\) and \(2^{C_E}\mathscr L_{\mathcal F}\in\Lambda[G_b]\) with \(C_E=c_{\mathrm{int}}\). Theorem~\ref{thm:central} gives the equivalence \(\mathscr L_{\mathcal F,x}(0)=0\Leftrightarrow L(E^{h_x},1)=0\), and with Lemma~\ref{lem:torsion} it gives \eqref{eq:exactcentral} with
\[
 \beta_E=2a_E+\alpha_E+1-\lambda_E-\rho_E-\epsilon_\infty .
\]
All these integers depend only on \(E\) and the fixed choices.
\end{proof}

\begin{definition}\label{def:bsd}
Let \(\mathcal T_E=\sum_{\ell\mid2N}\val(c_\ell(E))\).
For admissible \(h\) with \(L(E^h,1)\ne0\), its
\(2\)-adic BSD defect is
\[
 \bsd(h)=\val\left(\frac{L(E^h,1)}{\Omega_{\mathrm{real}}(E^h)}\right)
 -s_{\Sha}(h)-\sum_\ell\val(c_\ell(E^h))+2r(h),
 \qquad \Omega_{\mathrm{real}}(E^h)=2^{\epsilon_\infty}\Omega^+_{E^h}.
\]
Here \(\Omega_{\mathrm{real}}\) is the full real N\'eron period;
\(\Omega^+\) is half of it when \(E(\RR)\) has two components,
and equals it otherwise.
\end{definition}

\begin{corollary}\label{cor:bsd}
Assume \RF{} and put
\(\kappa_E=\alpha_E+1-\lambda_E-\rho_E+\mathcal T_E\).
For every admissible family and every index with
\(L(E^{h_x},1)\ne0\), including \(h_x=1\),
\begin{equation}\label{eq:bsd-determinant}
 \val(\mathscr L_{\mathcal F,x}(0))=\bsd(h_x)+\kappa_E.
\end{equation}
Consequently every admissible \(h\) with \(L(E^h,1)\ne0\) satisfies
\(\bsd(h)\ge-C_E-\kappa_E\).
\end{corollary}
\begin{proof}
One has \(\bsd(h)=D(h)+2r(h)-\mathcal T_E-\epsilon_\infty\).
Theorem~\ref{thm:central} gives \eqref{eq:bsd-determinant}.
Apply the integral bound in a family containing \(h\); for
\(h=1\), take \(b=1\) and \(P=\varnothing\).
\end{proof}
\section{Waldspurger coefficients and the minimum}\label{sec:symbol}
Proposition~\ref{prop:wald} supplies the half-integral-weight form and
Lemma~\ref{lem:mindepth} proves that its normalized coefficient valuations
attain a minimum. The reduced coefficient \(A(U)\) detects this minimum.
Theorem~\ref{thm:families} is the quadratic twist family construction
proved in the following three sections.

\subsection{Three cases and the auxiliary set}
The conclusion of Theorem~\ref{thm:A} is invariant under
isogeny.\footnote{An isogeny identifies rational Tate modules and
hence \(L\)-functions. Its maps on \(2^\infty\)-Selmer groups,
together with those of the dual isogeny, preserve corank.}
After replacing \(E\) by an isogenous curve, the rational
two-torsion case reduces to \CaseI{} or \CaseII{} below.
Abelian two-division fields give the following three cases:
\begin{enumerate}
\item[(I)]\phantomsection\label{case:I} \(E[2]\subset E(\Q)\);
\item[(II)]\phantomsection\label{case:II} \(E(\Q)[2]\simeq\Z/2\), and no curve isogenous to \(E\) has full rational two-torsion;
\item[(III)]\phantomsection\label{case:III} \(\Q(E[2])/\Q\) is cyclic cubic.
\end{enumerate}
In case~\CaseII{}, \(F=\Q(E[2])\) is a quadratic field ramified only at primes dividing \(2N\).

In case~\CaseI{}, let \(\Sigma\) be the primes dividing \(2N\).
In cases~\CaseII{} and~\CaseIII{}, put \(F=\Q(E[2])\) and enlarge
this set so that the classes of primes above \(\Sigma\) generate
\(\Cl(F)\). Thus the ring \(\mathcal O_\Sigma\) of
\(\Sigma\)-integers of \(F\) has trivial class group. Make this
fixed choice before constructing \(g\).

\begin{definition}\label{def:admissibleset}
For a finite set \(U\) of primes outside \(\Sigma\) put \(h(U)=\prod_{p\in U}p^*\). The set \(U\) is \emph{admissible} if \(h(U)>0\) and \(h(U)\in(\Q_\ell^\times)^2\) for every \(\ell\in\Sigma\). The \emph{weight} of a prime \(p\notin\Sigma\) is \(w(p)=t_p/2\), and that of \(U\) is \(w(U)=\sum_{p\in U}w(p)\).
\end{definition}

The set \(U\) is admissible exactly when \(h(U)\) is \(\Sigma\)-admissible in Definition~\ref{def:intro-family}. The interpolation theorem therefore applies with its original constants. For admissible \(U\), one has \(h(U)=\prod_{p\in U}p\) and \(w(U)=\tau(h(U))/2\).

When a set \(U\) occurs as a Fourier index, it denotes the positive product \(\prod_{p\in U}p\), even when \(U\) is not admissible. Thus \(a(U)\) means \(a(\prod_{p\in U}p)\), while \(h(U)=\prod_{p\in U}p^*\) is the signed twisting parameter. In case~\CaseI{} every \(w(p)\) is \(1\); in case~\CaseII{}, \(w(p)=1\) if \(p\) splits in \(F\) and \(w(p)=\tfrac12\) if \(p\) is inert in \(F\). In case~\CaseIII{}, \(w(p)=1\) at split primes and
\(w(p)=0\) at inert primes.

\begin{lemma}\label{lem:admissible-sign}
Let \(U\) be admissible. Then \(w(U)\in\Z\), and \(E^{h(U)}\) has the same root number as \(E\).
\end{lemma}

\begin{proof}
Let \(\chi_\Delta\) be the sign of the permutation action of
\(G_\Q\) on the three non-zero points of \(E[2]\). It is a
quadratic character, possibly trivial, ramified only in \(\Sigma\).
For \(p\notin\Sigma\), the three possible Frobenius types give
\(\chi_\Delta(p)=(-1)^{t_p}\). Admissibility and the product
formula give \(\prod_{p\in U}\chi_\Delta(p)=1\). Thus
\(\sum_{p\in U}t_p\) is even and \(w(U)\in\Z\). The root number of \(E^h\), \(h=h(U)\),
has the same local factors as \(E\) at \(\Sigma\) and infinity,
factor one at the other primes not dividing \(h\), and factor
\(\chi_{h,q}(-1)\) at \(q\mid h\). The product of these last
factors is one, again by the product formula.
\end{proof}

\subsection{Waldspurger's coefficients}\label{subsec:waldspurger}
Put \(q=e^{2\pi iz}\) and \(\theta(z)=\sum_{m\in\Z}q^{m^2}\).
We use Shimura's normalization of half-integral weight
\cite{Shimura73}, with multiplier \(j(\gamma,z)=\theta(\gamma z)/\theta(z)\). As \(j(-I,z)=1\), a non-zero form of weight \(\tfrac32\) and character \(\psi\) has \(\psi\) even, and then \(\eta=\psi\chi_{-4}\) is odd. For such a form \(g'\) and a prime \(p\) not dividing the level, Shimura's formula \cite[Theorem~1.7]{Shimura73} for the Hecke operator \(T(p^2)\) reads
\begin{equation}\label{eq:Tp2}
 (T(p^2)g')[n]=g'[p^2n]+\psi(p)\Bigl(\frac{-n}p\Bigr)g'[n]+p\,\psi(p)^2\,g'[n/p^2],
\end{equation}
with \(g'[n/p^2]=0\) if \(p^2\nmid n\). Waldspurger's correspondence supplies a weight-\(3/2\) eigenform with
non-zero coefficient at any prescribed non-vanishing admissible twist
\cite[Th\'eor\`eme~1]{Wald} (see also \cite{Wald80}).
At each \(\ell\in\Sigma\), the prescribed index is a local unit
square, so the required squareclass is the trivial class; at
infinity it is positive and the component is holomorphic of weight
\(3/2\). The weight-two form to which we apply the correspondence
is \(f\otimes\eta_1\). Choose \(\eta_1\), with conductor supported
on \(\Sigma\), as a product of an even character of sufficiently
large \(2\)-power conductor and odd local characters at the finite
places where the local root number of \(f\) is \(-1\), as in
\cite[\S2]{OS98b}. At odd \(\ell\) these may be taken of conductor
\(\ell\) or \(\ell^2\). The archimedean root number is \(-1\),
so global root number \(+1\) makes \(\eta_1\) odd and
\(\psi_1=\eta_1\chi_{-4}\) even. These local parities and the
dyadic ramification verify Waldspurger's local conditions, and hence
give a local vector with non-zero coefficient in the prescribed
squareclass. Outside
\(\Sigma\), choose the spherical vector: its local factor is one
at valuations zero and one \cite[\S I.3, pp.~377--378]{Wald}.
Thus the resulting half-integral-weight form has level supported
on \(\Sigma\). We use the following local-square-class form of
\cite[Corollaire~2]{Wald}. If \(g'\) has character \(\psi\), is a \(T(p^2)\)-eigenform outside the level, and corresponds
to a weight-two newform \(\varphi\), and \(\eta=\psi\chi_{-4}\), then
\begin{equation}\label{eq:wald-standard}
 a_{g'}(D)^2\sqrt{D_0}\,L(\varphi\otimes\eta^{-1}\chi_{D_0},1)
 =\psi(D/D_0)a_{g'}(D_0)^2\sqrt D\,L(\varphi\otimes\eta^{-1}\chi_D,1)
\end{equation}
for positive square-free \(D,D_0\), prime to the level, with
\(D/D_0\) a square at every prime of the level.

\begin{proposition}\label{prop:wald}
Suppose that the root number of \(E\) is \(+1\). There exist a character \(\psi\) of \(2\)-power order, with \(\eta=\psi\chi_{-4}\), a cusp form \(g=\sum_{n\ge1}a(n)q^n\) of weight \(\tfrac32\), character \(\psi\) and level supported on \(\Sigma\), and an admissible \(U_0\) with \(a(n_0)=1\) for \(n_0=h(U_0)\), with the following properties.
\begin{enumerate}
\item \(a(n)=0\) unless \(n\in(\Z_\ell^\times)^2\) for every \(\ell\in\Sigma\).
\item For every prime \(p\notin\Sigma\) and every \(n\ge1\),
\[
 a(p^2n)=\eta(p)\Bigl[\Bigl(a_p-\Bigl(\frac np\Bigr)\Bigr)a(n)-p\,\eta(p)\,a(n/p^2)\Bigr],
\]
where \((\frac np)\) is the Legendre symbol, zero if \(p\mid n\), and \(a(n/p^2)=0\) if \(p^2\nmid n\).
\item For admissible \(U\) and \(n=h(U)\), \(a(n)^2=\psi(n/n_0)\,\Lalg(E^n,1)/\Lalg(E^{n_0},1)\). In particular \(a(n)\ne0\) if and only if \(L(E^n,1)\ne0\), and then \(2\val(a(n))=\val(\Lalg(E^n,1))+c_g\) with \(c_g=-\val(\Lalg(E^{n_0},1))\).
\item For square-free \(n\), \(a(n)\) is a rational number times a value of \(\eta\). All coefficients of \(g\) lie in \(\Q(\mu_{2^k})\), where \(2^k\) is the order of \(\eta\).
\end{enumerate}
\end{proposition}

\begin{proof}
\emph{An initial twist.} The theorem of Friedberg and Hoffstein \cite{FH}, applied to the twists of \(f\) by quadratic characters trivial at the places of \(\Sigma\) and at infinity, which have root number \(+1\) by Lemma~\ref{lem:admissible-sign}, gives such a character \(\chi_{n_0}\) with \(L(E^{n_0},1)\ne0\), where \(n_0>0\) is square-free and a square at every \(\ell\in\Sigma\). As \(n_0\equiv1\pmod8\), we have \(n_0=h(U_0)\) for the set \(U_0\) of its prime factors, which is admissible (cf.~\cite[Proposition~7.1]{OAPW}).

\emph{Character and support.} Use the characters \(\psi_1,\eta_1\)
chosen above. Waldspurger's correspondence, with \(D_0=n_0\), gives
\(g_1\) of level \(L_1\) supported on \(\Sigma\) and character
\(\psi_1\). Write \(\psi_1=\psi\mu^{-2}\) with \(\psi\) of \(2\)-power order and \(\mu\) of odd order, which is possible because squaring is bijective on characters of odd order. The twist of \(g_1\) by \(\mu\) has character \(\psi\), level supported on \(\Sigma\), \(T(p^2)\)-eigenvalues \(\mu(p)^2\eta_1(p)a_p=\eta(p)a_p\), and coefficient at \(n_0\) multiplied by \(\mu(n_0)\neq0\) \cite{Shimura73}.

For odd \(\ell\in\Sigma\), with \((U_\ell g')[n]=g'[\ell n]\) and \((V_\ell g')[n]=g'[n/\ell]\), the operator \(g'\mapsto\tfrac12\bigl(g'-V_\ell U_\ell g'+g'\otimes(\cdot/\ell)\bigr)\) keeps the coefficients at \(n\) with \((\frac n\ell)=1\) and kills the others, and \(g'\mapsto\tfrac14\sum_{\chi\bmod8}g'\otimes\chi\) keeps those with \(n\equiv1\pmod8\). These operators preserve the weight, the character \(\psi\) (a quadratic twist multiplies it by a square, and \(V_\ell\), \(U_\ell\) multiply it by the same quadratic character) and the \(T(p^2)\)-eigenvalues for \(p\notin\Sigma\); they raise the level only at \(\ell\), and they keep the coefficient at \(n_0\), which is a unit square at every \(\ell\in\Sigma\). Applying them at every \(\ell\in\Sigma\) and normalizing at \(n_0\) gives \(g\) with (1), of level \(L\) supported on \(\Sigma\). By \eqref{eq:Tp2} with \(T(p^2)g=\eta(p)a_pg\),
\[
 a(p^2n)=\Bigl(\eta(p)a_p-\psi(p)\Bigl(\frac{-n}p\Bigr)\Bigr)a(n)-p\,\psi(p)^2a(n/p^2),
\]
and since \(\psi(p)=\eta(p)\chi_{-4}(p)\) and \(\chi_{-4}(p)(\frac{-n}p)=(\frac np)\), this is (2).

\emph{Central values.} Let \(f_\eta=f\otimes\eta\), and let \(\mathcal V\)
be the space of weight-\(3/2\) cusp forms at the level of \(g\), with
character \(\psi\) and eigenvalues \(a_{f_\eta}(p)\) for \(T(p^2)\) outside
\(\Sigma\). For admissible \(D\), both \(D\) and \(n_0\) are unit
squares at the level primes. Formula~\eqref{eq:wald-standard}, together
with \(f_\eta\otimes\eta^{-1}\chi_D=f\otimes\chi_D\), gives
\begin{equation}\label{eq:waldcomp}
 a(D)^2\sqrt{n_0}\,L(E^{n_0},1)
 =\psi(D/n_0)\sqrt D\,L(E^D,1).
\end{equation}
Lemma~\ref{lem:period} proves (3).

\emph{Rationality of the coefficients.} We show that square-free
coefficients are rational numbers times roots of unity. By (3) the square of every square-free coefficient in the support is algebraic, so these coefficients are algebraic, and by (2) so are all coefficients. Let \(D\) be square-free in the support of \(g\). Formula~\eqref{eq:wald-standard}, applied to every \(g'\in\mathcal V\), gives \(a_{g'}(D)^2=c\,a_{g'}(n_0)^2\) with \(c\) independent of \(g'\). The product of the linear forms \(a_{g'}(D)-\sqrt c\,a_{g'}(n_0)\) and \(a_{g'}(D)+\sqrt c\,a_{g'}(n_0)\) vanishes on \(\mathcal V\), and a vector space is not the union of two proper subspaces, so one of them vanishes identically. Thus the ratio \(a_{g'}(D)/a_{g'}(n_0)\) is the same for all \(g'\in\mathcal V\) with \(a_{g'}(n_0)\neq0\).

Let \(\sigma\in\operatorname{Aut}(\C)\). It acts on the values of \(\psi\) and \(\eta\) by \(\zeta\mapsto\zeta^e\) for an odd integer \(e\); put \(\vartheta_\sigma=\eta^{(e-1)/2}\). The form \(g^\sigma\), obtained by applying \(\sigma\) to the coefficients, is a cusp form of weight \(\tfrac32\) with character \(\psi^e\) \cite[\S1]{Shimura81} and \(T(p^2)\)-eigenvalues \(\eta(p)^ea_p\). The form \(g\otimes\vartheta_\sigma\) has character \(\psi\vartheta_\sigma^2=\psi\eta^{e-1}=\psi^e\), because \(\psi/\eta=\chi_{-4}\) and \(e-1\) is even, and the same eigenvalues \(\vartheta_\sigma(p)^2\eta(p)a_p\). At a common level supported on \(\Sigma\) both lie in one space of the type of \(\mathcal V\), and comparing their ratios gives \(\sigma(a(D))=\vartheta_\sigma(D/n_0)\,a(D)\).

As \(D/n_0\) is a unit square at the primes dividing the conductor of \(\eta\), \(\eta(D/n_0)=z^2\) for a value \(z\) of \(\eta\); then \(\sigma(z)=z^e\) and \(\vartheta_\sigma(D/n_0)=z^{e-1}\), so \(\sigma(a(D)/z)=a(D)/z\) for every \(\sigma\), and \(a(D)/z\in\Q\). The other coefficients in the support are obtained from the square-free ones by (2), which involves only values of \(\eta\) and rational integers; this gives the last assertion.
\end{proof}

Let \(\mathcal O_g\) be the ring of integers of
\(\Q_2(\mu_{2^k})\); its residue field is \(\F_2\).

\begin{remark}\label{rem:integral-coefficient-valuations}
Proposition~\ref{prop:wald}(4) and the shell recurrence in
Subsection~\ref{subsec:forms} express square-free coefficients and
shell coefficients of the divided translates as rational numbers
times roots of unity. Their normalized valuations are integers.
In particular, if \(r\in\Q_2\), \(\zeta\in\mu_{2^\infty}\)
and \(r\zeta\in\mathcal O_g\), then
\[
 \overline{r\zeta}=0\quad\Longrightarrow\quad r\zeta\in2\mathcal O_g.
\]
This implication is used when dividing the coefficient forms by \(2\).
\end{remark}

\subsection{The symbol}\label{subsec:symbol}
Assume \RF{} throughout this subsection. Return to the fixed curve \(E\) and the form \(g\) of Proposition~\ref{prop:wald}. The contribution of the new Tamagawa factors is \(w(U)\). Subtracting it from the coefficient valuation allows comparison of sets with different numbers of prime factors.

\begin{definition}\label{def:normval}
The \emph{normalized valuation} of an admissible set \(U\) with \(a(h(U))\ne0\) is \(\val(a(h(U)))-w(U)\).
\end{definition}

\begin{lemma}\label{lem:mindepth}
Suppose that the root number of \(E\) is \(+1\). For a finite set \(U\) of primes outside \(\Sigma\), \(a(U)=0\) unless \(U\) is admissible. For admissible \(U\) with \(a(h(U))\ne0\) and \(n=h(U)\),
\begin{equation}\label{eq:normval}
 2\bigl(\val(a(n))-w(U)\bigr)=D(n)+s_{\Sha}(n)+c_g .
\end{equation}
The normalized valuations are integers, and for \(U\ne\varnothing\) they are at least \((c_g-C_E-\beta_E)/2\).
\end{lemma}

\begin{proof}
If \(a(U)\neq0\), then \(U\) lies in the support, so \(\prod_{p\in U}p\equiv1\pmod8\); as \(h(U)\equiv1\pmod4\) and has the same absolute value, \(h(U)=\prod_{p\in U}p\), and \(U\) is admissible. For admissible \(U\) with \(a(n)\ne0\), Proposition~\ref{prop:wald}(3) gives \(L(E^n,1)\ne0\), so \(E^n\) has Selmer corank zero, and \eqref{eq:normval} is Proposition~\ref{prop:wald}(3) with \(\tau(n)=2w(U)\). By Proposition~\ref{prop:wald}(4) and Lemma~\ref{lem:admissible-sign} its left side is an even integer. For \(U\ne\varnothing\), Corollary~\ref{cor:sha} gives \(D(n)\ge-C_E-\beta_E\).
\end{proof}

The lower bound, together with the value at \(1\) when \(a(1)\ne0\),
makes the following minimum well defined.

\begin{definition}\label{def:symbol}
Suppose that the root number of \(E\) is \(+1\). The \emph{minimal normalized valuation} \(\delta\) is the minimum
of \(\val(a(h(U)))-w(U)\) over admissible sets \(U\) with
\(a(h(U))\ne0\). This set of integers is bounded below by
Lemma~\ref{lem:mindepth} and the preceding observation, and is
nonempty because \(a(n_0)=1\). The \emph{symbol} of a finite set
\(U\) of primes outside \(\Sigma\) is
\[
 A(U)=\begin{cases}\text{the residue of }2^{-\delta-w(U)}a(h(U))\in\mathcal O_g,&U\text{ admissible},\\ 0,&U\text{ not admissible}.\end{cases}
\]
\end{definition}

\begin{lemma}\label{lem:depth}
Suppose that the root number of \(E\) is \(+1\). For admissible \(U\), \(A(U)=1\) exactly when \(\val(a(h(U)))=\delta+w(U)\), and then
\[
 L(E^{h(U)},1)\ne0,\qquad D(h(U))\le C_A:=2\delta-c_g .
\]
Some admissible \(U\) has \(A(U)=1\).
\end{lemma}

\begin{proof}
For admissible \(U\), \(w(U)\) is an integer, so \(2^{-\delta-w(U)}a(h(U))\) is a rational number times a root of unity, by Proposition~\ref{prop:wald}(4); it has non-negative valuation by the choice of \(\delta\) and lies in \(\mathcal O_g\), and it has residue \(1\) exactly when it is a unit. In that case \(L(E^{h(U)},1)\ne0\), and \eqref{eq:normval} gives \(D(h(U))=2\delta-c_g-s_{\Sha}(h(U))\le C_A\). A set attaining the minimum \(\delta\) has symbol \(1\).
\end{proof}

\begin{lemma}\label{lem:selmerbound}
Suppose that the root number of \(E\) is \(+1\), and let \(U\ne\varnothing\) be admissible with \(A(U)=1\). Then
\[
 s_{\Sha}(h(U))\le C_A+C_E+\beta_E,\qquad\dim_{\F_2}\Sel_2(E^{h(U)}/\Q)\le\dim_{\F_2}E(\Q)[2]+C_A+C_E+\beta_E .
\]
\end{lemma}

\begin{proof}
By Lemma~\ref{lem:depth}, \(L(E^{h(U)},1)\ne0\), so \(E^{h(U)}\) has rank zero and finite Shafarevich--Tate group, and \eqref{eq:normval} gives \(s_{\Sha}(h(U))=C_A-D(h(U))\). Corollary~\ref{cor:sha} gives \(D(h(U))\ge-C_E-\beta_E\). The \(2\)-Selmer group is an extension of \(\Sha(E^{h(U)}/\Q)[2]\), of dimension at most \(s_{\Sha}(h(U))\), by \(E^{h(U)}(\Q)/2\), whose dimension is that of \(E^{h(U)}(\Q)[2]=E(\Q)[2]\).
\end{proof}

\subsection{Families at the minimum}

\begin{theorem}\label{thm:families}
Suppose that \(E\) is in case \CaseI{}, \CaseII{}, or \CaseIII{}, has root number \(+1\), and \(A(\varnothing)=0\). For every \(b\ge1\) there are a finite set \(P\) of primes outside \(\Sigma\) and non-zero ramification forms \(\alpha_q\in(\F_2^b)^*\), \(q\in P\), such that for every \(x\ne0\) the set \(U_x=\{q\in P:\alpha_q(x)=1\}\) is admissible and \(A(U_x)=1\).
\end{theorem}

By Lemma~\ref{lem:depth}, \(A(\varnothing)=0\) says that \(\val(a(1))>\delta\), and \(A(U_x)=1\) says that \(a(h_x)\ne0\) and \(\val(a(h_x))-w(U_x)=\delta\); so Theorem~\ref{thm:families} is Theorem~\ref{thm:E}. Section~\ref{sec:rules} proves it in case~\CaseI{},
Section~\ref{sec:residual} in case~\CaseII{}, and
Section~\ref{sec:cyclic-cubic} in case~\CaseIII{}. As \(A(\varnothing)=0\), every \(U_x\) with \(x\neq0\) is non-empty.

Theorem~\ref{thm:C} follows from Theorem~\ref{thm:main} and Proposition~\ref{prop:RF-cases}; Proposition~\ref{prop:wald} and Lemma~\ref{lem:mindepth} prove Theorem~\ref{thm:D}. It remains to prove Theorem~\ref{thm:families}.

\section{Coefficient congruences and simultaneous non-vanishing}\label{sec:rules}
The main result is Theorem~\ref{thm:families} in case~\CaseI{}.
We first prove the contraction and isolation rules of
Theorem~\ref{thm:rules}, then use auxiliary primes to obtain simultaneous
non-vanishing. Throughout this section \(E[2]\subset E(\Q)\), so
\(w(U)=|U|\), and \(\Sigma\) consists of the primes dividing \(2N\).

A prime is fresh relative to \(U\) if it lies outside \(\Sigma\cup U\); fresh primes used together are taken distinct.

\subsection{Labels and graph data}
Let \(\mathcal W_\Sigma\subset\Q^\times/\Q^{\times2}\) be the subgroup generated by \(-1\) and the primes of \(\Sigma\), and \(\mathcal W_\Sigma^\vee=\operatorname{Hom}(\mathcal W_\Sigma,\F_2)\). The \emph{label} of a prime \(p\notin\Sigma\) is \(\lambda(p)\in\mathcal W_\Sigma^\vee\), \(w\mapsto[(\frac wp)=-1]\), and for distinct primes \(p,q\notin\Sigma\) put \(\varepsilon(p,q)=[(\frac qp)=-1]\in\F_2\). Quadratic reciprocity reads
\begin{equation}\label{eq:recip-bits}
 \varepsilon(p,q)+\varepsilon(q,p)=\lambda(p)(-1)\,\lambda(q)(-1).
\end{equation}

\begin{definition}\label{def:datum}
A \emph{graph datum} \(\Dat=(V,\lambda,\varepsilon)\) consists of a finite set \(V\), a map \(\lambda\colon V\to\mathcal W_\Sigma^\vee\), and \(\varepsilon\colon\{(p,q)\in V^2:p\ne q\}\to\F_2\) satisfying \eqref{eq:recip-bits}. Its \emph{total label} is \(\lambda(\Dat)=\sum_{p\in V}\lambda(p)\). A finite set \(U\) of primes outside \(\Sigma\) has the datum \(\Dat(U)\) given by Legendre symbols, and \(U\) \emph{realizes} \(\Dat\) if \(\Dat\simeq\Dat(U)\).
\end{definition}

A datum records, for each vertex, its residue symbols against \(-1\) and the primes of \(\Sigma\), and for each pair of vertices one independent residue symbol.

\begin{lemma}\label{lem:filter}
A finite set \(U\) of primes outside \(\Sigma\) is admissible if and only if \(\lambda(\Dat(U))=0\).
\end{lemma}

\begin{proof}
Write \(h=h(U)\). The sign of \(h\) is \(\prod_{p\in U}(\frac{-1}p)\), so \(h>0\) exactly when \(\lambda(\Dat(U))(-1)=0\). Each \(p^*\) is \(1\) modulo \(4\), and \(p^*\equiv1\pmod8\) exactly when \((\frac2p)=1\); so \(h\in(\Q_2^\times)^2\), that is \(h\equiv1\pmod8\), exactly when \(\lambda(\Dat(U))(2)=0\). For odd \(\ell\in\Sigma\), reciprocity gives \((\frac{p^*}\ell)=(\frac\ell p)\), so \((\frac h\ell)=\prod_{p\in U}(\frac\ell p)\).
\end{proof}

\begin{lemma}\label{lem:realize}
Let \(U\) be a finite set of primes outside \(\Sigma\), \(Z\) a finite set of primes, and \(\Dat\) a graph datum on \(U\sqcup V'\) restricting to \(\Dat(U)\) on \(U\). There are distinct primes \(q_v\notin\Sigma\cup U\cup Z\), \(v\in V'\), such that \(U\cup\{q_v\}\) realizes \(\Dat\).
\end{lemma}

\begin{proof}
Order \(V'\) and choose the primes successively. When \(q_v\) is chosen, its label and the bits \(\varepsilon(v,r)\), for \(r\) in \(U\) or among the earlier vertices, prescribe the values at \(q_v\) of the quadratic characters attached to \(-1\), to the primes of \(\Sigma\) and to the earlier primes. These classes are independent in \(\Q^\times/\Q^{\times2}\), so the prescriptions define a non-empty union of residue classes, and Dirichlet's theorem supplies infinitely many primes in it. The reverse bits \(\varepsilon(r,v)\) are then determined by \eqref{eq:recip-bits}, consistently with \(\Dat\).
\end{proof}

\begin{definition}\label{def:operations}
For a datum \(\Dat\) and distinct \(p,q\in V\), the \emph{flip} \(\Dat^{pq}\) is the datum with \(\varepsilon(p,q)\) and \(\varepsilon(q,p)\) both changed, and the \emph{contraction} \(\Dat/pq\) is the datum on \((V\smallsetminus\{p,q\})\sqcup\{c\}\) with
\[
 \lambda(c)=\lambda(p)+\lambda(q),\qquad\varepsilon(c,r)=\varepsilon(p,r)+\varepsilon(q,r),\qquad\varepsilon(r,c)=\varepsilon(r,p)+\varepsilon(r,q),
\]
the other data unchanged. A vertex \(p\) is \emph{isolated} if \(\lambda(p)=0\) and \(\varepsilon(p,r)=0\) for every \(r\ne p\).
\end{definition}

Both satisfy \eqref{eq:recip-bits}, and \(\lambda(\Dat^{pq})=\lambda(\Dat/pq)=\lambda(\Dat)\). If \(c\) is realized by a prime, its residue symbols are the products of those of \(p\) and \(q\): its Frobenius in the maximal elementary abelian \(2\)-extension unramified outside \(\Sigma\cup V\smallsetminus\{p,q\}\) is the product of theirs. A vertex \(p\) is isolated exactly when \(p\) splits completely in \(\Q(\sqrt{\mathcal W_\Sigma},\sqrt r:r\in V\smallsetminus\{p\})\).

\begin{theorem}\label{thm:rules}
The symbol \(A\) of Definition~\ref{def:symbol} has the following properties for admissible sets \(U\).
\begin{enumerate}
\item \emph{Graph dependence.} \(A(U)\) depends only on \(\Dat(U)\).
\item \emph{Contraction.} For distinct \(p,q\in U\), \(A(\Dat(U))+A(\Dat(U)^{pq})=A(\Dat(U)/pq)\).
\item \emph{Isolation.} If some \(p\in U\) is isolated in \(\Dat(U)\), then \(A(U)=0\).
\end{enumerate}
By \textup{(1)} and Lemma~\ref{lem:realize}, \(A\) is a function on graph data of total label zero, and \textup{(2)} is an identity between such functions.
\end{theorem}

In the contraction rule, the merged vertex has the products of the two old residue symbols (cf.~\cite[Lemmas~7.6--7.8 and Proposition~7.9]{OAPW}). Multiplication by \(\theta\) and the use of Galois representations to read coefficients modulo \(2\) go back to Ono and Skinner \cite{OS98a,OS98b}.

Throughout the proof, \(\mathcal O_g\) is as in Definition~\ref{def:symbol}, \(\mathfrak m\) is its maximal ideal, and \(\bar x\in\F_2\) is the residue of \(x\in\mathcal O_g\). A root of unity of \(2\)-power order has residue \(1\).

\subsection{Forms attached to the symbol}\label{subsec:forms}
The integral forms \(H_U\) below have fresh-prime coefficients
\(A(Uq)\). Unary theta corrections at the indices \(Jk^2\),
\(J\subseteq U\), permit one further division of \(B_U\) by \(2\).
Multiplication by \(\theta\equiv1\pmod2\) gives a weight-two form
with the same reduced coefficients. Its Galois traces factor through
quadratic residue data (Lemma~\ref{lem:factor}); the local inertia
identity gives contraction (Lemma~\ref{lem:inertia}).

\subsubsection{The Hecke recurrence}\label{subsec:hecke}
Fix \(p\notin\Sigma\) and an integer \(m\) prime to \(\Sigma\) with \(v_p(m)\le1\). Put \(\kappa=(\frac mp)\in\{0,\pm1\}\), and define \(b_0=1\), \(b_1=a_p-\kappa\), \(b_{s+1}=a_pb_s-pb_{s-1}\). Iterating Proposition~\ref{prop:wald}(2), and noting that \((\frac np)=0\) once \(p\mid n\), gives
\begin{equation}\label{eq:shell}
 a(p^{2s}m)=\eta(p)^sb_s\,a(m)\qquad(s\ge0).
\end{equation}

\begin{lemma}\label{lem:recurrence}
Suppose \(a_p\) is even. If \(p\mid m\), then \(b_s\equiv s+1\pmod2\), so every \(b_s-b_{s-1}\) \((s\ge1)\) is odd. If \(p\nmid m\), then every \(b_s\) is odd and every \(\Delta_s=b_{s+1}-b_s\) is even; put \(\nu_p(m)=\Delta_0/2=(a_p-\kappa-1)/2\), whose parity changes when \(\kappa\) changes sign. If moreover \(a_p\equiv p+1\pmod4\), then \(\Delta_s\equiv\Delta_0\pmod4\) for all \(s\).
\end{lemma}

\begin{proof}
The first case follows from \(b_0=1\), \(b_1=a_p\) even and \(b_{s+1}\equiv b_{s-1}\pmod2\). In the second, \(b_1=a_p-\kappa\) is odd and \(b_{s+1}\equiv b_{s-1}\pmod2\), so all \(b_s\) are odd; replacing \(\kappa\) by \(-\kappa\) changes \(\Delta_0/2\) by \(\kappa=\pm1\). The differences satisfy \(\Delta_{s+1}=a_p\Delta_s-p\Delta_{s-1}\), so \(\Delta_{s+1}\equiv-\Delta_{s-1}\equiv\Delta_{s-1}\pmod4\). Finally \(\Delta_1-\Delta_0=a_p(a_p-\kappa-2)+2\kappa+1-p\); if \(a_p\equiv p+1\pmod4\), this is congruent modulo \(4\) to \((p+1)(p-1-\kappa)+2\kappa+1-p=(p-1)(p-\kappa)\), a product of two even numbers.
\end{proof}

In case~\CaseI{}, \(E[2]\subset E(\F_p)\) gives \(4\mid p+1-a_p\), so both hypotheses hold at every \(p\notin\Sigma\).

\subsubsection{Divided operators}
Lemma~\ref{lem:HU} gives the coefficients used in the residue rules.

For a prime \(p\notin\Sigma\) let \((U_pg')[n]=g'[pn]\) and \((V_pg')[n]=g'[n/p]\), with \(g'[n/p]=0\) if \(p\nmid n\). For a finite set \(U\) of primes outside \(\Sigma\) put
\[
 B_U=2^{-\delta-|U|}\prod_{p\in U}\bigl(U_p-\eta(p)V_p\bigr)g .
\]
The operators \(V_p\) and \(U_p\) multiply the character of a weight-\(\tfrac32\) form by \((\frac p\cdot)\) and raise the level by \(p\) \cite[Proposition~1.3]{Shimura73}, so \(B_U\) is a weight-\(\tfrac32\) form of character \(\psi\prod_{p\in U}(\frac p\cdot)\) and level supported on \(\Sigma\cup U\), with exponent one at each \(p\in U\).

The coefficients of \(B_U\) vanish at indices divisible by
\(\Sigma\). Write any other index as \(n=Jrk^2\), where
\(J\subseteq U\) and \(r\) are square-free, with \(r\) prime to
\(U\). The shell recurrence gives
\begin{equation}\label{eq:factor}
 B_U[n]=2^{-\delta-|U|}a(D)\prod_\ell\mu_\ell(n),\qquad
 D=(U\smallsetminus J)r.
\end{equation}
At \(p\in U\smallsetminus J\), with \(v_p(n)=2s\), the multiplier is
\(\eta(p)^s(b_s-b_{s-1})\), using the recurrence based at an index
divisible exactly once by \(p\), and \(b_{-1}=0\). It is odd after
removing the root of unity. At \(p\in J\), with \(v_p(n)=2s+1\),
the multiplier is \(\eta(p)^{s+1}\Delta_s\), using the recurrence
based at an index prime to \(p\); it is even. At \(\ell\mid k\)
outside \(U\), the multiplier is \(\eta(\ell)^sb_s\), integral and
odd if \(\ell\nmid r\). The Legendre symbols in these recurrences
depend only on \(D\), the common square-free part of the indices
before applying the operators.

Thus \(B_U[n]=0\) unless \(D\) is in the support of \(g\), in which case \(D=h(U_D)\) for an admissible set \(U_D\) (Lemma~\ref{lem:mindepth}), and the minimality of \(\delta\) applies to \(a(D)\).

\begin{lemma}\label{lem:BU}
Let \(n=Jrk^2\) be as above.
\begin{enumerate}
\item \(B_U[n]\in\mathcal O_g\), and \(\val(B_U[n])\ge w(r)\), the number of primes dividing \(r\).
\item If \(r=1\), then \(B_U[Jk^2]\equiv\eta(k)B_U[J]\pmod{2\mathcal O_g}\).
\item For distinct fresh primes \(q,q'\): \(B_U[q]=2^{-\delta-|U|}a(Uq)\), and \(\val(B_U[qq'])\ge2\).
\item For \(p\in U\), \(U'=U\smallsetminus\{p\}\) and a fresh prime \(q'\): \(B_U[pq']=2\eta(p)\,\nu_p(U'q')\cdot2^{-\delta-|U|}a(U'q')\).
\end{enumerate}
\end{lemma}

\begin{proof}
In \eqref{eq:factor} the \(|J|\) multipliers at \(p\in J\) are divisible by \(2\), and the minimality of \(\delta\) gives \(\val(a(D))\ge\delta+w(D)=\delta+|U|-|J|+w(r)\). This proves (1), and in particular
\begin{equation}\label{eq:depthJ}
 \val\bigl(2^{-\delta-|U|}a(U\smallsetminus J)\bigr)\ge-|J|.
\end{equation}
For (2), \(r=1\) and \(D=U\smallsetminus J\). Write each multiplier at \(Jk^2\) as \(\eta(\cdot)^s\) times the corresponding multiplier at \(J\), plus an error. For \(\ell\mid k\) outside \(U\) the multipliers are \(\eta(\ell)^sb_s\) and \(1\), with \(b_s\) odd, so the error lies in \(2\mathcal O_g\). For \(p\in U\smallsetminus J\) they are \(\eta(p)^s(b_s-b_{s-1})\) and \(1\), so again the error lies in \(2\mathcal O_g\). For \(p\in J\) they are \(\eta(p)^{s+1}\Delta_s\) and \(\eta(p)\Delta_0\in2\mathcal O_g\), and the error lies in \(4\mathcal O_g\) by Lemma~\ref{lem:recurrence}. Expanding the product, \(B_U[Jk^2]-\eta(k)B_U[J]\) is \(2^{-\delta-|U|}a(D)\) times a sum of terms in \(2^{|J|+1}\mathcal O_g\), and \eqref{eq:depthJ} gives (2).

Assertion (3) is the case \(J=\varnothing\), \(k=1\), with \(w(r)=2\) for \(qq'\). For (4), \(n=pq'\) has \(J=\{p\}\), \(r=q'\), \(k=1\), \(D=U'q'\), and the multiplier at \(p\) is \(\eta(p)\Delta_0=2\eta(p)\nu_p(U'q')\).
\end{proof}

Recall that \(\psi\) is even and \(\eta=\psi\chi_{-4}\) is odd. For square-free \(J\) prime to \(\Sigma\) put
\[
 \Theta_J=\sum_{k\ge1,\ (k,\Sigma)=1}\eta(k)\,k\,q^{Jk^2},
\]
a weight-\(\tfrac32\) cusp form of character \(\eta\chi_{-4}(\frac J\cdot)=\psi(\frac J\cdot)\) and level supported on \(\Sigma\cup J\), with exponent one at the primes of \(J\) \cite[Proposition~2.2]{Shimura73}, \cite{SerreStark}. Its coefficients are supported on the indices \(Jk^2\), and \(\Theta_J[Jk^2]\equiv\eta(k)\pmod2\) as \(k\) is odd. Since \(\mathcal O_g\) is ramified over \(\Z_2\), dividing by \(2\) requires numerators in \(2\mathcal O_g\), not merely in \(\mathfrak m\); Lemma~\ref{lem:BU}(2) supplies this for
\[
 H_U=\tfrac12\Bigl(B_U-\sum_{J\subseteq U}B_U[J]\,\Theta_J\Bigr).
\]

\begin{lemma}\label{lem:HU}
The form \(H_U\) has coefficients in \(\mathcal O_g\). For distinct fresh primes \(q,q'\), \(p\in U\) and \(U'=U\smallsetminus\{p\}\),
\[
 \overline{H_U[q]}=A(U\cup\{q\}),\qquad\overline{H_U[qq']}=0,\qquad\overline{H_U[pq']}=\overline{\nu_p(U'q')}\,A(U'\cup\{q'\}).
\]
\end{lemma}

\begin{proof}
At indices with \(r\ne1\), \(B_U\) is even by Lemma~\ref{lem:BU}(1) and the \(\Theta_J\) vanish. At \(Jk^2\), Lemma~\ref{lem:BU}(2) and \(\Theta_J[Jk^2]\equiv\eta(k)\) show that \(B_U[Jk^2]-B_U[J]\Theta_J[Jk^2]\) is even, and the other \(\Theta_{J'}\) vanish there. The indices \(q\), \(qq'\), \(pq'\) are not of the form \(Jk^2\) with \(J\subseteq U\), so \(H_U=\tfrac12B_U\) there. By Lemma~\ref{lem:BU}(3), \(H_U[q]=2^{-\delta-|U\cup\{q\}|}a(Uq)\), whose residue is \(A(U\cup\{q\})\) by definition; the other two formulas follow in the same way from Lemma~\ref{lem:BU}(3)--(4) and \(\overline{\eta(p)}=1\).
\end{proof}

\subsubsection{Weight two and traces of Frobenius}\label{subsec:trace-operators}
With Shimura's normalization of half-integral weight by the theta multiplier, \(f_U=\theta H_U\) is a weight-two modular form on \(\Gamma_1(N_U)\), with \(N_U\) supported on \(\Sigma\cup U\) and exponent one at each \(p\in U\), and \(\theta f\) has, in weight two, the character of \(f\). Since \(\theta\equiv1\pmod2\) coefficientwise, \(f_U[n]\equiv H_U[n]\pmod{2\mathcal O_g}\).

Let \(M_U\) be the \(\mathcal O_g\)-module spanned by the \(q\)-expansions of \(f_U\), of the forms \(\theta\Theta_J\) \((J\subseteq U)\), and of all their images under the Hecke operators \(T_\ell\), \(\ell\nmid N_U\). The character \(\chi_B=\psi(\frac U\cdot)\) of \(\theta B_U\) and the characters \(\chi_J=\psi(\frac J\cdot)\) of the \(\theta\Theta_J\) have \(2\)-power order, and \(\chi_B/\chi_J=(\frac{UJ}\cdot)\) is quadratic. Hence
\[
 \langle u\rangle f_U=\chi_B(u)f_U+\tfrac12\sum_JB_U[J]\bigl(\chi_B(u)-\chi_J(u)\bigr)\theta\Theta_J,
\]
where each \(\chi_B(u)-\chi_J(u)\) is \(0\) or twice a root of unity; so \(M_U\) is stable under the diamond and Hecke operators (cf.~\cite[Lemma~7.7]{OAPW}). It is a finite \(\mathcal O_g\)-module, since finitely many coefficients determine a form of this level, and it is closed.

Let \(\mathbf T_U\subset\operatorname{End}_{\mathcal O_g}(M_U)\) be the subalgebra generated by the \(T_\ell\) and \(\langle\ell\rangle\), \(\ell\nmid N_U\). The generators of \(M_U\) are cusp forms: \(g\), and hence \(B_U\), is cuspidal, the \(\Theta_J\) are cusp forms, and multiplication by \(\theta\) preserves vanishing at the cusps. Over \(\overline\Q_2\) the operators of \(\mathbf T_U\) act semisimply, and \(M_U\otimes\overline\Q_2\) is a sum of cuspidal Hecke eigensystems \(j\), each with a multiplicity space.

To each eigensystem, Eichler--Shimura and Deligne attach a two-dimensional representation \(\rho_j\) of \(G_{\Q,\Sigma\cup U}\), with trace \(T_\ell\) and determinant \(\ell\langle\ell\rangle\) at \(\mathrm{Frob}_\ell\) \cite{Deligne}. Define \(t(g),d(g)\in\mathbf T_U\otimes\overline\Q_2\) as the operators acting on eigensystem \(j\) by the scalars \(\operatorname{tr}\rho_j(g)\) and \(\det\rho_j(g)\).

Then \(t(\mathrm{Frob}_\ell)=T_\ell\) and \(d(\mathrm{Frob}_\ell)=\ell\langle\ell\rangle\), so by continuity, Chebotarev density and closedness of \(\mathbf T_U\), the maps \(t,d\colon G_{\Q,\Sigma\cup U}\to\mathbf T_U\) are continuous and satisfy, as on each eigensystem,
\begin{equation}\label{eq:pseudo}
 t(g)t(u)=t(gu)+d(u)\,t(gu^{-1}),\qquad t(1)=2 .
\end{equation}
For \(g\in G_{\Q,\Sigma\cup U}\) put \(\Phi_U(g)=\overline{(t(g)f_U)[1]}\in\F_2\). Since \((T_\ell f)[1]=f[\ell]\) and \((T_\ell T_{\ell'}f)[1]=f[\ell\ell']\) in weight two for distinct \(\ell,\ell'\nmid N_U\), Chebotarev density gives, for primes \(q\ne q'\) not dividing \(N_U\),
\begin{equation}\label{eq:reading}
 \Phi_U(\mathrm{Frob}_q)=\overline{f_U[q]},\qquad\overline{\bigl(t(\mathrm{Frob}_q)t(\mathrm{Frob}_{q'})f_U\bigr)[1]}=\overline{f_U[qq']}.
\end{equation}
Let \(K_U=\Q(\sqrt w,\sqrt p:w\in\mathcal W_\Sigma,\ p\in U)\), the maximal elementary abelian \(2\)-extension of \(\Q\) unramified outside \(\Sigma\cup U\) and infinity.

The vanishing of the two-fresh-prime coefficients is the crucial point. By \eqref{eq:pseudo}, it gives \(\Phi_U(gu)=\Phi_U(gu^{-1})\), hence invariance under squares. This makes \(\Phi_U\) factor through the multiquadratic extension \(K_U\).

\begin{lemma}\label{lem:factor}
The function \(\Phi_U\) factors through \(\Gal(K_U/\Q)\), and \(\Phi_U(1)=0\).
\end{lemma}

\begin{proof}
First, \(\overline{(d(u)t(v)f_U)[1]}=\overline{(t(v)f_U)[1]}\) for all \(u,v\). Indeed \(d(u)\) acts on \(f_U\) as the cyclotomic character, which is \(\equiv1\pmod2\), times a diamond operator, which by the formula above acts as a root of unity plus an integral combination of the \(\theta\Theta_J\). By Chebotarev density it suffices to take \(v=\mathrm{Frob}_q\) with \(q\) fresh, and \((\theta\Theta_J)[q]\equiv\Theta_J[q]=0\pmod2\) because \(q\) is not of the form \(Jk^2\).

Next, for all \(g,u\), \eqref{eq:reading}, Lemma~\ref{lem:HU} and \(f_U\equiv H_U\) give \(\overline{(t(g)t(u)f_U)[1]}=0\). Applying \eqref{eq:pseudo} to \(f_U\) at the coefficient \(1\), we obtain \(\Phi_U(gu)=\Phi_U(gu^{-1})\), hence \(\Phi_U(gu^2)=\Phi_U(g)\).

Thus \(\Phi_U\) is constant on the cosets of the closed subgroup generated by squares, and factors through the maximal elementary abelian \(2\)-quotient of \(G_{\Q,\Sigma\cup U}\), which is \(\Gal(K_U/\Q)\). Finally \(\Phi_U(1)=\overline{2f_U[1]}=0\).
\end{proof}

\subsubsection{The inertia identity}\label{subsec:inertia}
Fix \(p\in U\). Let \(\tau_p\) be an element of inertia at \(p\) acting by \(-1\) on \(\sqrt p\), and \(\varphi_p\) a Frobenius lift in the same decomposition group.

\begin{lemma}\label{lem:inertia}
For every \(g\in G_{\Q,\Sigma\cup U}\) and \(g'\in M_U\otimes\overline\Q_2\),
\[
 U_p^{(2)}\bigl(t(g\tau_p\varphi_p)g'-t(g\varphi_p)g'\bigr)=d(\varphi_p)\bigl(t(g\tau_p)g'-t(g)g'\bigr),
\]
in the space of weight-two forms of level \(N_U\), where \(U_p^{(2)}\) is the weight-two operator \(g'\mapsto\sum g'[pn]q^n\).
\end{lemma}

\begin{proof}
Decompose \(g'\) into eigensystems. Each occurs in the Hecke span of \(\theta B_U\) or of some \(\theta\Theta_J\) with \(p\in J\), whose characters are ramified at \(p\) with \(p\)-part the Legendre symbol and whose levels have exponent one at \(p\), or in the Hecke span of some \(\theta\Theta_J\) with \(p\notin J\), whose level is prime to \(p\).

An eigensystem of level prime to \(p\) has \(\rho_j\) unramified at \(p\), so both differences vanish.

Every other eigensystem has conductor exponent one at \(p\) and a central character whose restriction to inertia at \(p\) is the ramified quadratic character. A Steinberg constituent with this conductor would have unramified central character, and a supercuspidal one would have conductor exponent at least two; so the local representation is a principal series \(\pi(\mu_1,\mu_2)\) with \(\mu_1\) unramified and \(\mu_2\) ramified.

By local--global compatibility at \(p\ne2\) \cite[Th\'eor\`eme~(A)]{Carayol}, with the normalization of Frobenius, the \(U_p\)-eigenvalue of the newform is the Frobenius eigenvalue \(\lambda\) on the unramified line. Since the \(p\)-exponents of \(N_U\)
and the newform conductor are both one, its oldforms are
\(B_df_j(z)=f_j(dz)\) with \(p\nmid d\).
The operator \(U_p\) commutes with these \(B_d\), so acts by
\(\lambda\) on the whole eigensystem at level \(N_U\).

In a basis adapted to the two lines, \(\rho_j(\tau_p)=\mathrm{diag}(1,-1)\) and \(\rho_j(\varphi_p)=\mathrm{diag}(\lambda,\mu)\); if \(z\) is the lower-right entry of \(\rho_j(g)\), the left side is \(\lambda(-2z\mu)\) and the right side is \(\lambda\mu(-2z)\).
\end{proof}

\subsection{Proof of Theorem~\ref{thm:rules}}\label{subsec:proofrules}
By Lemmas~\ref{lem:HU} and~\ref{lem:factor}, the symbol can be read at any one prime through its Frobenius.

\subsubsection{The last vertex, and isolation}
Let \(U'\) be a finite set of primes outside \(\Sigma\), and \(q\notin\Sigma\cup U'\) a prime with \(\lambda(q)=\lambda(\Dat(U'))\), so that \(U'\cup\{q\}\) is admissible. By Lemma~\ref{lem:HU} and \eqref{eq:reading}, \(A(U'\cup\{q\})=\Phi_{U'}(\mathrm{Frob}_q)\), and by Lemma~\ref{lem:factor} this depends only on the Frobenius of \(q\) in \(K_{U'}\), that is, on \(\lambda(q)\) and the bits \(\varepsilon(q,r)\), \(r\in U'\). Every vertex of an admissible set may play the role of \(q\). An isolated vertex has trivial Frobenius in \(K_{U'}\), and \(\Phi_{U'}(1)=0\); this is the isolation rule.

\subsubsection{Graph dependence}
Let \(U,U'\) realize the same labelled datum. By Lemma~\ref{lem:realize}
choose a disjoint auxiliary realization \(U''\) whose pair bits against
both old realizations match the corresponding bits of the datum.
Choose its primes successively; the prescriptions involve distinct
prime radicals, so reciprocity is the only compatibility condition.
Replacing the primes of \(U\) by those of \(U''\), one at a time,
preserves the Frobenius of the replaced prime relative to all remaining
primes. Thus it preserves \(A\) by Lemma~\ref{lem:factor}. The same
replacement from \(U'\) gives \(A(U)=A(U'')=A(U')\). If \(U\) and
\(U'\) overlap, first replace \(U\) by a disjoint realization in the
same way, then apply this construction.

\subsubsection{Contraction}
By graph dependence and Lemma~\ref{lem:realize}, it suffices to treat a datum \(\Dat(U\cup\{q\})\) with \(U\cup\{q\}\) admissible, \(p\in U\) and \(U'=U\smallsetminus\{p\}\). Put \(g=\mathrm{Frob}_q\). The reductions of \(t\) and \(d\) on \(M_U/2M_U\) are locally constant, so they factor through a finite quotient of \(G_{\Q,\Sigma\cup U}\), which we take to dominate \(\Gal(K_U/\Q)\). Chebotarev density gives fresh primes \(q_1,q'_0,q'_1\) whose Frobenius elements in that quotient are \(g\tau_p\), \(g\varphi_p\), \(g\tau_p\varphi_p\); the prime \(q\) itself has Frobenius \(g\). Apply Lemma~\ref{lem:inertia} to \(f_U\) and take the coefficient at \(1\). Since \((U_p^{(2)}g')[1]=g'[p]\) and \((T_{q'}f)[p]=f[pq']\) for \(q'\nmid pN_U\), the left side reduces to \(\overline{f_U[pq'_1]}-\overline{f_U[pq'_0]}\), and by the first step of the proof of Lemma~\ref{lem:factor} the right side reduces to \(\overline{f_U[q_1]}-\overline{f_U[q]}\). Lemma~\ref{lem:HU} turns this into
\[
 A(U\cup\{q_1\})+A(U\cup\{q\})=\overline{\nu_p(U'q'_1)}\,A(U'\cup\{q'_1\})+\overline{\nu_p(U'q'_0)}\,A(U'\cup\{q'_0\}).
\]
The element \(\tau_p\) acts trivially on \(K_{U'}\), so \(q'_0\) and \(q'_1\) have the same Frobenius there and \(A(U'\cup\{q'_1\})=A(U'\cup\{q'_0\})\). On the other hand \(\tau_p\) changes the sign of \(\sqrt{p^*}\), so \((\frac{p^*}{q'_1})=-(\frac{p^*}{q'_0})\), and by reciprocity \((\frac{q'_1}p)=-(\frac{q'_0}p)\); thus \(\kappa=(\frac{U'q'}p)\) changes sign between \(U'q'_0\) and \(U'q'_1\), and \(\overline{\nu_p(U'q'_1)}+\overline{\nu_p(U'q'_0)}=1\) by Lemma~\ref{lem:recurrence}. It remains to identify the two resulting data.

The Frobenius of \(q_1\) differs from that of \(q\) only on \(\sqrt p\), so \(U\cup\{q_1\}\) realizes \(\Dat(U\cup\{q\})^{pq}\). In \(K_{U'}\) the Frobenius of \(q'_0\) is \(\mathrm{Frob}_q\mathrm{Frob}_p\), so the residue symbols of \(q'_0\) are the products of those of \(q\) and \(p\), and \(U'\cup\{q'_0\}\) realizes \(\Dat(U\cup\{q\})/pq\). This is the contraction rule and completes the proof of Theorem~\ref{thm:rules}.\hfill\qedsymbol\par

\subsection{Auxiliary graphs and simultaneous non-vanishing}\label{sec:graphs}

Finite differences contract forests, and auxiliary vertices realize
these differences on the terminal polynomials. Selectors lower the
degree at a chosen index while preserving the leading terms at all
others. Iteration gives \(A=1\) simultaneously.

Throughout this subsection \(A\) is any \(\F_2\)-valued function on graph data of total label zero satisfying Theorem~\ref{thm:rules}.

\subsubsection{Forests}
For a datum \(\Dat\) and a pair \(e=\{p,q\}\) of vertices put \(\Delta_eA(\Dat)=A(\Dat)+A(\Dat^{pq})\). These difference operators commute, and \(\Delta_e^2=0\). For a set \(F\) of pairs write \(\Delta_F=\prod_{e\in F}\Delta_e\), and if \(F\) is a forest let \(\Dat/F\) be the datum obtained by contracting each of its trees to one vertex, with labels and bits summed as in Section~\ref{sec:rules}.

\begin{lemma}\label{lem:forest}
Let \(F\) be a set of pairs of vertices of \(\Dat\). If \(F\) is a forest, then \(\Delta_FA(\Dat)=A(\Dat/F)\). If \(F\) contains a cycle, then \(\Delta_FA(\Dat)=0\).
\end{lemma}

\begin{proof}
For a forest, apply the contraction rule successively: each remaining
edge still joins distinct contracted vertices. If \(F\) contains a
cycle, choose one edge \(e\) and a path \(F'\subseteq F\setminus\{e\}\)
joining its endpoints. Then \(\Dat/F'=\Dat^e/F'\), since the bit of
\(e\) becomes internal to one contracted vertex. Thus
\(\Delta_e\Delta_{F'}A=0\), and hence \(\Delta_FA=0\).
\end{proof}

\begin{definition}\label{def:terminal}
Fix an ordered set \(T\) of vertices whose labels sum to zero. The \emph{terminal polynomial} \(\mathscr P_T\) is the function \((\varepsilon(p,q))_{p<q}\mapsto A(\Dat)\) of the independent bits of the data \(\Dat\) on \(T\) with these labels, viewed as a multilinear polynomial over \(\F_2\). The ideal \(\mathfrak d\) is generated by the \(\Delta_e\), for pairs \(e\) in \(T\), in the commutative algebra of operators they generate.
\end{definition}

On multilinear polynomials \(\Delta_e\) is differentiation in the bit of \(e\); hence \(\mathfrak d^r\) lowers degree by at least \(r\). A forest spanning the blocks of a partition \(\pi\) acts by contraction to \(\Dat/\pi\).

\subsubsection{Networks}
An auxiliary vertex of label zero joined to two terminals \(a,b\)
acts by \(\Delta_{ab}\): expansion in its two incident bits gives
\(0+A(\Dat_T)+A(\Dat_T)+A(\Dat_T/ab)\).
The next lemma gives the corresponding formula for an arbitrary graph.

\begin{definition}\label{def:network}
Let \(\Dat_T\) be a datum on a set \(T\) of \emph{terminals}, \(I\) a finite set of \emph{auxiliary vertices} disjoint from \(T\), and \(E_0\) a set of pairs in \((I\times I)\cup(I\times T)\). The \emph{network datum} \(\Dat_T\oplus_{E_0}I\) on \(T\sqcup I\) extends \(\Dat_T\) by \(\lambda=0\) on \(I\), by \(\varepsilon=1\) on the pairs of \(E_0\), and by \(0\) on all other pairs meeting \(I\); \eqref{eq:recip-bits} holds since \(\lambda(i)(-1)=0\). Its \emph{reduced Laplacian} is the matrix \(Q\in M_I(\F_2)\) with \(Q_{ij}=[\{i,j\}\in E_0]\) for \(i\ne j\) and \(Q_{ii}\) the number of pairs of \(E_0\) containing \(i\): the Laplacian of the graph \((I\sqcup T,E_0)\), reduced modulo \(2\), with the rows and columns of the terminals removed. The \emph{incidence vector} of a terminal \(a\) is \(\mathbf c_a\in\F_2^I\), with \((\mathbf c_a)_i=[\{i,a\}\in E_0]\).
\end{definition}

\begin{lemma}\label{lem:network}
There is an operator \(O\), a sum of products \(\Delta_F\) over forests \(F\) of pairs in \(T\), such that \(A(\Dat_T\oplus_{E_0}I)=O\,A(\Dat_T)\) for every datum \(\Dat_T\) on \(T\) with the given labels, and
\begin{equation}\label{eq:leading}
 O\equiv\det Q+\sum_{\{a,b\}\subseteq T}\mathbf c_a^{\mathsf t}\operatorname{adj}(Q)\,\mathbf c_b\,\Delta_{ab}\pmod{\mathfrak d^2}.
\end{equation}
If \(I=I_1\sqcup I_2\) and \(E_0=E_1\sqcup E_2\) with \(E_j\) meeting only \(I_j\), the operator for \(E_0\) is the product of those for \(E_1\) and \(E_2\).
\end{lemma}

\begin{proof}
Start with all bits meeting \(I\) zero and expand in the bits of
\(E_0\). By Lemma~\ref{lem:forest}, only forests contribute.
A component without terminals, including an untouched auxiliary vertex, contracts to an isolated vertex and
contributes zero. Every other forest contracts the terminal partition
\(\pi\) given by its components, so contributes \(\Delta_{F_\pi}\),
where \(F_\pi\) spans its blocks. This defines \(O\).
Its term has order \(|T|-\#\pi\) in \(\mathfrak d\).

The constant term counts forests with exactly one terminal in each
component, and is \(\det Q\) by the matrix-tree theorem. The
coefficient of \(\Delta_{ab}\) counts those with just one two-terminal
component, joining \(a,b\). By the all-minors formula \cite{Chaiken},
this count is
\[
 \det\begin{pmatrix}Q&\mathbf c_b\\\mathbf c_a^{\mathsf t}&0\end{pmatrix}
 =\mathbf c_a^{\mathsf t}\operatorname{adj}(Q)\mathbf c_b
 \quad\text{in }\F_2.
\]
This proves \eqref{eq:leading}.

For disjoint auxiliary graphs, expand both operators. An acyclic
union gives the product of their terminal contractions. A cycle in
the union gives a cycle of terminal differences or a repeated
difference; its product vanishes by Lemma~\ref{lem:forest} or
\(\Delta_e^2=0\). Components without terminals vanish in either
expansion. Thus the operators compose as asserted.
\end{proof}

\subsubsection{Selectors}
A selector acts at one chosen family index. At that index it takes a first difference in a chosen terminal bit; at every other index it preserves the leading homogeneous part. 

\begin{definition}\label{def:familyconfig}
A \emph{family configuration} is a finite set of vertices, each with a non-zero \emph{ramification form} \(\alpha\in(\F_2^b)^*\); the vertices present at \(x\) are those with \(\alpha(x)=1\), so that none is present at \(x=0\). A network on a family configuration of terminals has auxiliary vertices with non-zero ramification forms; at an index \(x\), its present part is a network on the terminals present at \(x\), with reduced Laplacian \(Q_x\) and operator \(O_x\) on the terminal polynomial at \(x\) (Lemma~\ref{lem:network}). For \(x^*\ne0\) and distinct terminals \(a_0,a_1\) present at \(x^*\), a \emph{selector} for \((x^*;a_0,a_1)\) is such a network with
\[
 O_{x^*}\equiv\Delta_{a_0a_1}\pmod{\mathfrak d^2},\qquad O_x\equiv1\pmod{\mathfrak d}\quad(x\ne x^*).
\]
\end{definition}

An operator in \(1+\mathfrak d\) preserves the leading homogeneous
part. We seek an auxiliary matrix invertible on every proper
subset of its pairs, with exactly one kernel direction when all
pairs are present.

\begin{figure}[htbp]
\centering
\begin{tikzpicture}[every node/.style={font=\small},
  vertex/.style={circle,fill=black,inner sep=1.6pt}]
 \node[vertex,label=above:\(i\)] (i) at (0,1.2) {};
 \node[vertex,label=left:\(a\)] (a) at (-1,0) {};
 \node[vertex,label=right:\(b\)] (b) at (1,0) {};
 \draw (a)--(i)--(b);
 \draw[dashed] (a)--(b);
 \node at (0,-.55) {\(O=\Delta_{ab}\)};
 \begin{scope}[xshift=5cm,yshift=.5cm]
  \foreach \name/\angle/\lab in {u1/150/u_1,v2/90/v_2,u2/30/u_2,v3/-30/v_3,u3/-90/u_3,v1/-150/v_1}
    \node[vertex,label={\angle:\(\lab\)}] (\name) at (\angle:1cm) {};
  \draw (u1)--(v2)--(u2)--(v3)--(u3)--(v1)--(u1);
 \end{scope}
\end{tikzpicture}
\caption{Left: a label-zero auxiliary vertex realizes a difference in
the dashed terminal bit. Right: the internal graph of the selector for
\(b=3\); \(u_i,v_i\) have form \(\ell_i\). Terminal attachments,
which determine the diagonal entries of \(Q\), are omitted.}
\label{fig:graphs}
\end{figure}
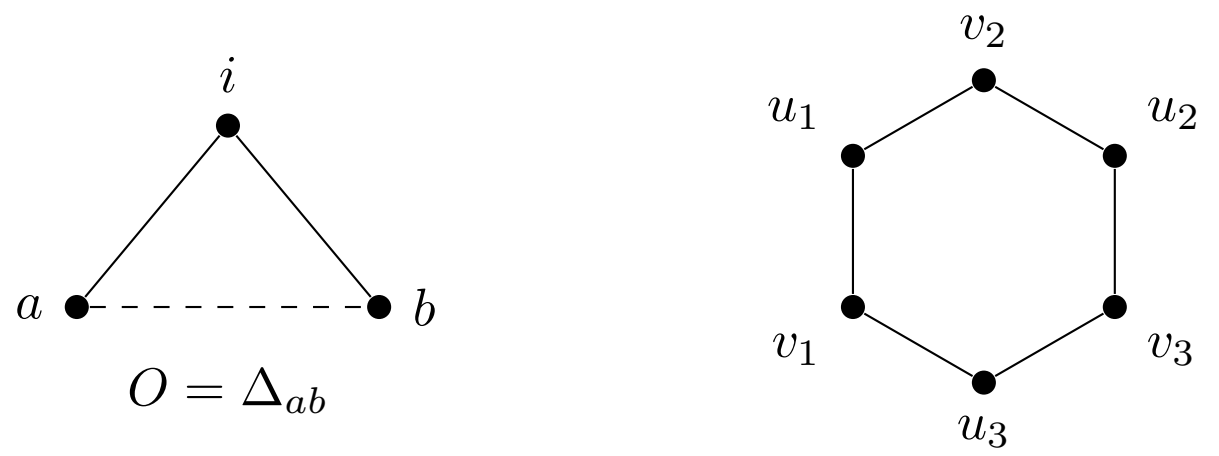

\begin{lemma}\label{lem:selector}
Suppose that a family configuration of terminals contains a vertex of every non-zero ramification form. Let \(x^*\ne0\), and let \(a_0,a_1\) be distinct terminals present at \(x^*\). There is a selector for \((x^*;a_0,a_1)\) on \(2b\) auxiliary vertices of label zero.
\end{lemma}

\begin{proof}
Choose a basis \(\ell_1,\dots,\ell_b\) of \((\F_2^b)^*\) with
\(\ell_i(x^*)=1\) for every \(i\) (cf.~\cite[Lemma~8.6(1)]{OAPW}). Take auxiliary vertices \(u_i,v_i\)
of ramification form \(\ell_i\), and put
\[
 Q=\begin{pmatrix}0&\mathsf T\\\mathsf T^{\mathsf t}&\mathsf E\end{pmatrix},
 \qquad \mathsf T=1+\mathsf P,
\]
in the order \((u_1,\dots,u_b,v_1,\dots,v_b)\), where \(\mathsf P\)
is a cyclic permutation matrix and \(\mathsf E_{11}=1\) is the only
non-zero entry of \(\mathsf E\). Use the off-diagonal entries of \(Q\)
as the internal pair bits. For each \(i\), attach a terminal of form
\(\ell_i\) through the incidence vector \(Qe_i\), where
\(e_i=e_{u_i}+e_{v_i}\). Attach \(\mathbf c=e_{u_1}\) to both
\(a_0,a_1\). If their forms \(\alpha,\beta\) differ, attach
\(\mathbf c\) also to a terminal of form \(\alpha+\beta\).
A terminal with several roles receives the sum of its incidence vectors.

These last attachments contribute zero to the total incidence at every
index, since \(\alpha+\beta+(\alpha+\beta)=0\). Let \(I_x\) be the
present auxiliary vertices. As
\(\mathbf1_{I_x}=\sum_i\ell_i(x)e_i\), the diagonal entry of the reduced
Laplacian at a present vertex \(v\) is
\[
 \sum_{w\in I_x\smallsetminus\{v\}}Q_{vw}
 +(Q\mathbf1_{I_x})_v=Q_{vv}.
\]
Thus this Laplacian is the principal submatrix \(Q_J\), where
\(J=\{i:\ell_i(x)=1\}\) specifies the present pairs \(u_i,v_i\).

For proper \(J\), the restriction \(\mathsf P_J\) is nilpotent, so
\(\mathsf T_J\), and hence \(Q_J\), is invertible. At \(x^*\) all
pairs are present. Here
\(\ker\mathsf T=\ker\mathsf T^{\mathsf t}=\F_2\mathbf1\).
If \((u,v)\in\ker Q\), then \(v=s\mathbf1\) and
\(\mathsf T^{\mathsf t}u=s\mathsf E\mathbf1\). The left side has
coordinate sum zero, whereas the right side has sum \(s\). Thus
\(s=0\), and \(\ker Q=\F_2r\), with \(r=\sum_i e_{u_i}\).
Consequently \(\det Q=0\) and \(\operatorname{adj}Q=rr^{\mathsf t}\).
Now \(r^{\mathsf t}Qe_i=0\) and \(r^{\mathsf t}\mathbf c=1\);
the balancing terminal is absent at \(x^*\). Formula~\eqref{eq:leading}
therefore gives the required first difference at \(x^*\), and constant
term \(1\) elsewhere. For \(b=1\), \(\mathsf T=0\) and the same
calculation applies.
\end{proof}

\subsubsection{Simultaneous non-vanishing}
The following construction makes all the non-zero terminal polynomials equal to one, by successively decreasing the sum of their degrees.

\begin{theorem}\label{thm:address}
Let \(\mathcal T\) be a family configuration of primes outside \(\Sigma\), containing a vertex of every non-zero ramification form, with \(\lambda(\mathcal T_x)=0\) for every \(x\), where \(\mathcal T_x\) is the set of primes present at \(x\). Suppose that for each \(x\neq0\) the terminal polynomial \(\mathscr P_x\) of \(\Dat(\mathcal T_x)\) is non-zero. Then there are finitely many further primes outside \(\Sigma\), of label zero and with non-zero ramification forms, such that the enlarged configuration \(\mathcal U\) satisfies \(A(\mathcal U_x)=1\) for every \(x\neq0\).
\end{theorem}

\begin{proof}
Add selectors one at a time, on disjoint sets of auxiliary vertices with no pairs between different selectors, keeping the terminals \(\mathcal T\). By Lemma~\ref{lem:network}, after any number of steps the symbol at \(x\ne0\), as a function of the terminal bits, is \(R_x=O_x\mathscr P_x\), where \(O_x\) is the product of the operators at \(x\) of the selectors added so far. Initially \(R_x=\mathscr P_x\ne0\). If some \(R_{x^*}\) has positive degree \(d\), choose a variable \(z\), the bit of a pair \(\{a_0,a_1\}\) of terminals present at \(x^*\), occurring in a monomial of degree \(d\) of \(R_{x^*}\); then \(\Delta_zR_{x^*}\) has degree exactly \(d-1\). Adding a selector of Lemma~\ref{lem:selector} for \(x^*\) and \(a_0,a_1\) acts by
\[
\begin{aligned}
 R_{x^*}&\longmapsto\Delta_zR_{x^*}+(\text{terms of degree at most }d-2),\\
 R_x&\longmapsto R_x+(\text{terms of lower degree})\qquad(x\ne x^*).
\end{aligned}
\]
The new \(R_{x^*}\) has degree \(d-1\), and every other \(R_x\) keeps its degree and leading part. Thus the nonnegative integer
\[
 \sum_{x\ne0}\deg R_x
\]
decreases by one at each step, while every \(R_x\) remains non-zero. After finitely many steps all \(R_x\) are the constant \(1\).

 By Lemma~\ref{lem:realize} the prescribed bits are realized by primes of label zero, without changing the terminal data, and then \(A(\mathcal U_x)=R_x=1\). The auxiliary vertices have label zero, so \(\lambda(\mathcal U_x)=\lambda(\mathcal T_x)=0\).
\end{proof}

\subsection{The initial configuration in case~\CaseI{}}\label{sec:seeds}

\begin{proof}[Proof of Theorem~\ref{thm:families} in case~\CaseI{}]
By Lemma~\ref{lem:depth}, choose an admissible \(U\) with \(A(U)=1\).
It is nonempty since \(A(\varnothing)=0\). Write
\(\varsigma=\Dat(U)\), with vertex labels \(\mu_1,\dots,\mu_r\).
For every non-zero \(\ell\in(\F_2^b)^*\), take \(r\) terminals of
these labels, all with ramification form \(\ell\); realize their data
by Lemma~\ref{lem:realize}. Each such group has total label zero, so
\(\lambda(\mathcal T_x)=0\) for every \(x\).

Fix \(x\ne0\) and one present group \(p_1,\dots,p_r\). Let \(F\)
be a spanning tree on
\(C=\mathcal T_x\smallsetminus\{p_1,\dots,p_{r-1}\}\).
By Lemma~\ref{lem:forest}, \(\Delta_F\mathscr P_x\) is the symbol of an
\(r\)-vertex datum whose labels are \(\mu_1,\dots,\mu_r\).
Its independent bits are sums over disjoint nonempty sets of the
remaining bits. They can therefore be assigned to recover
\(\varsigma\), giving \(\Delta_F\mathscr P_x=1\) and hence \(\mathscr P_x\ne0\).

Theorem~\ref{thm:address} now supplies an enlarged configuration
\(\mathcal U\) with \(A(\mathcal U_x)=1\) for every \(x\ne0\).
Its auxiliary vertices have label zero, so each \(\mathcal U_x\)
remains admissible by Lemma~\ref{lem:filter}. Take \(P\) to be its
prime set and \(\alpha_q\) the ramification form of \(q\).
\end{proof}
\section{The quadratic two-division field}\label{sec:residual}

This section proves Theorem~\ref{thm:families} in case~\CaseII{}.
The graph rules of Proposition~\ref{prop:qgraph}, including removal
of pairs of inert primes, give the simultaneous construction in
Lemma~\ref{lem:qpacket}.

\subsection{Split, inert and doubly inert primes}\label{subsec:qprimes}
Throughout, \(E\) is in case~\CaseII{}: \(E(\Q)[2]\) is generated by a point \(P_0\) of order two, and no curve isogenous to \(E\) has full rational two-torsion. Put
\[
 E'=E/\langle P_0\rangle,\qquad F=\Q(E[2]),\qquad F'=\Q(E'[2]).
\]
At a prime inert in both \(F\) and \(F'\), one has
\(\#E(\F_p)\equiv2\pmod4\). This permits removal of inert-prime
pairs without changing the symbol.

The fields \(F,F'\) are quadratic and ramified only at primes dividing \(2N\). Indeed, \(E'\) has rational two-torsion but not full rational two-torsion, by the isogeny hypothesis. Write \(d,d'\in\mathcal W_\Sigma\) for their discriminant classes.

Take \(\Sigma\) as fixed in Section~\ref{sec:symbol}, so
\(\mathcal O_\Sigma\) has trivial class group.\footnote{The
label criterion and prime realization of Lemmas~\ref{lem:filter}
and~\ref{lem:realize} allow any finite \(\Sigma\) containing the
primes dividing \(2N\).}

\begin{definition}\label{def:doublyinert}
A prime \(p\notin\Sigma\) is \emph{split} or \emph{inert} according to its behavior in \(F\), and \emph{doubly inert} if it is inert in both \(F\) and \(F'\). A rational label \(\mu\) is \emph{doubly inert} if \(\mu(d)=\mu(d')=1\).
\end{definition}

In terms of labels, \(p\) is inert exactly when \(\lambda(p)(d)=1\), and doubly inert exactly when its label is doubly inert.
At split primes \(t_p=2\), \(w(p)=1\), and \(a_p\equiv p+1\pmod4\); at inert primes \(t_p=1\) and \(w(p)=1/2\).
An admissible set contains an even number of inert primes by Lemma~\ref{lem:admissible-sign}.

\begin{lemma}\label{lem:good}
If \(p\) is doubly inert, then \(\#E(\F_p)\equiv2\pmod4\).
\end{lemma}

\begin{proof}
As \(p\) is inert in \(F\), \(E(\F_p)[2]=\langle P_0\rangle\). If \(4\mid\#E(\F_p)\), the \(2\)-primary part of \(E(\F_p)\) is cyclic of order at least four and contains \(P_1\) with \(2P_1=P_0\). The image of \(P_1\) in \(E'(\F_p)\) has order two and does not lie in the kernel \(E[2]/\langle P_0\rangle\) of the dual isogeny, since \(P_1\notin E[2]\). Together with the rational point of that kernel it gives full two-torsion in \(E'(\F_p)\), so \(p\) splits in \(F'\), a contradiction.
\end{proof}

For these primes see also \cite[Lemma~2.1(ii)]{Zhai25} (cf.~\cite[Proposition~7.9(4)]{OAPW}).

\subsection{Coefficient forms and the rules for one new prime}\label{subsec:qcoeff}

\subsubsection{Divided operators}
The weight of a set can now be half-integral. We distinguish the parity of its inert primes so that the divided translates still use an integral power of \(2\).

\begin{definition}\label{def:matching}
A finite set \(U\) of primes outside \(\Sigma\) is \emph{matching} if it contains an even number of inert primes, and \emph{opposite} otherwise. Put \(e_U=0\) or \(\tfrac12\) accordingly, so that \(n_U=\delta+w(U)+e_U\) is an integer.
\end{definition}

The correction \(e_U\) makes \(w(U)+e_U\) integral. Define
\[
 B_U=2^{-n_U}\prod_{p\in U}\bigl(U_p-\eta(p)V_p\bigr)g,\qquad
 H_U=\begin{cases}\tfrac12\Bigl(B_U-\sum_{J}B_U[J]\,\Theta_J\Bigr)&U\text{ matching},\\[2pt] B_U&U\text{ opposite},\end{cases}
\]
with \(\Theta_J\) as in Section~\ref{sec:rules} and \(J\) running over the sets of split primes of \(U\). In case~\CaseI{} these are the forms of Section~\ref{sec:rules}.

\begin{lemma}\label{lem:BUq}
Write \(n=Jrk^2\) as in Section~\ref{sec:rules}. Then \(B_U[n]\) is a rational number times a root of unity of \(2\)-power order,
\[
 \val(B_U[n])\ \ge\ w(r)+\tfrac12\#\{p\in J\ \text{inert}\}-e_U,
\]
and \(B_U\) and \(H_U\) have coefficients in \(\mathcal O_g\). For distinct fresh primes \(q,q'\), \(p\in U\) and \(U'=U\smallsetminus\{p\}\), modulo \(\mathfrak m\):
\begin{enumerate}
\item \(H_U[q]\equiv A(U\cup\{q\})\);
\item \(H_U[pq']\equiv \nu_p(U'q')\,A(U'\cup\{q'\})\), except that \(H_U[pq']\equiv0\) when \(U\) is opposite and \(p\) is inert;
\item \(H_U[qq']\equiv A(U\cup\{q,q'\})\) if \(U\) is matching and \(q,q'\) are inert, and \(H_U[qq']\equiv0\) otherwise.
\end{enumerate}
\end{lemma}

\begin{proof}
Use \eqref{eq:factor} with exponent \(n_U\). Each multiplier is an
integer times a value of \(\eta\), and every \(p\in J\) contributes
a factor \(2\). The underlying square-free part is
\(D=(U\setminus J)r\), so minimality gives
\[
 \val(B_U[n])\ge |J|-w(J)+w(r)-e_U.
\]
The valuation is an integer by Proposition~\ref{prop:wald}(4), and
the right side is at least \(-1/2\). Thus \(B_U\) is integral.
For matching \(U\), it is even unless \(r=1\) and \(J\) consists of
split primes. At those exceptional indices the congruence
\(\Delta_s\equiv\Delta_0\pmod4\) at split primes gives
\(B_U[Jk^2]\equiv\eta(k)B_U[J]\pmod{2\mathcal O_g}\), exactly as in
Lemma~\ref{lem:BU}. Subtracting the displayed unary theta series
therefore makes \(H_U\) integral.

At the indices in (1)--(3), the unary corrections vanish, and
\[
 B_U[q]=2^{-n_U}a(Uq),\quad
 B_U[qq']=2^{-n_U}a(Uqq'),\quad
 B_U[pq']=2\eta(p)\nu_p(U'q')2^{-n_U}a(U'q').
\]
An admissible set has integral weight. Comparing these exponents with
\(\delta+w(Uq)\), \(\delta+w(Uqq')\), and \(\delta+w(U'q')\) gives
the three formulas. In (3), a non-zero residue can occur only when \(U\) is matching and both new primes are inert. In (2), opposite \(U\)
and inert \(p\) leave an additional factor \(2\).
\end{proof}

\subsubsection{Traces}
For matching \(U\) put \(f_U=\theta H_U\), and for opposite \(U\) put \(f_U=\theta B_U\). Define \(N_U\), \(M_U\), \(\mathbf T_U\), \(t\), \(d\) and \(\Phi_U\) as in Section~\ref{sec:rules}; for opposite \(U\) the diamond operators act on \(f_U\) by a character of \(2\)-power order. Let \(\Kt_U\) be the maximal elementary abelian \(2\)-extension of \(F\) unramified outside the places above \(\Sigma\cup U\) and infinity. It is Galois over \(\Q\) and contains \(K_U\).

\begin{lemma}\label{lem:factorq}
If \(U\) is matching, \(\Phi_U\) factors through \(\Gal(\Kt_U/\Q)\) and vanishes outside \(\Gal(\Kt_U/F)\). If \(U\) is opposite, \(\Phi_U\) factors through \(\Gal(K_U/\Q)\) and vanishes on \(\Gal(K_U/F)\).
\end{lemma}

\begin{proof}
The first step of the proof of Lemma~\ref{lem:factor} applies, and so does \eqref{eq:reading}. Let \(G_{F,\Sigma\cup U}\) be the kernel of \(G_{\Q,\Sigma\cup U}\to\Gal(F/\Q)\). By Lemma~\ref{lem:BUq} and Chebotarev density, \(\overline{(t(g)t(u)f_U)[1]}=0\) for all \(g,u\) if \(U\) is opposite, and for all \(g\) and all \(u\in G_{F,\Sigma\cup U}\) if \(U\) is matching, since then a split prime and any second prime give a vanishing coefficient. The identity \eqref{eq:pseudo} gives \(\Phi_U(gu^2)=\Phi_U(g)\) for all \(g\) and all \(u\in G_{\Q,\Sigma\cup U}\), respectively all \(u\in G_{F,\Sigma\cup U}\). The closed subgroup generated by these squares is normal in \(G_{\Q,\Sigma\cup U}\), with quotient \(\Gal(K_U/\Q)\), respectively \(\Gal(\Kt_U/\Q)\). The vanishing assertions are the readings of Lemma~\ref{lem:BUq} at primes \(q\) for which \(U\cup\{q\}\) has an odd number of inert primes, hence is not admissible.
\end{proof}

\subsubsection{The rules for one new prime}
For a prime \(p\notin\Sigma\) let \(\tau_p\) be an element of inertia at a place above \(p\) acting by \(-1\) on \(\sqrt p\); it lies in \(G_F\), as \(F\) is unramified at \(p\). Let \(\varphi_p\) be a lift of the rational arithmetic Frobenius at \(p\) in the decomposition group of the same place; \(\varphi_p\in G_F\) exactly when \(p\) is split. For a set \(U\) and a prime \(q\), the \emph{Frobenius of \(q\) relative to \(U\)} is its Frobenius class in \(\Gal(\Kt_U/\Q)\) if \(q\) is split and in \(\Gal(K_U/\Q)\) if \(q\) is inert.

\begin{proposition}\label{prop:qrules}
Let \(U\) be a finite set of primes outside \(\Sigma\) and \(q\notin\Sigma\cup U\) a prime.
\begin{enumerate}
\item \(A(U\cup\{q\})\) depends only on the Frobenius of \(q\) relative to \(U\).
\item If \(q\) is split and its Frobenius in \(\Kt_U\) is trivial, then \(A(U\cup\{q\})=0\).
\item Let \(p\in U\), \(U'=U\smallsetminus\{p\}\), let \(g\) be a Frobenius element of \(q\) in the group attached to \(q\) relative to \(U\), let \(q_1\) be a prime with Frobenius \(g\tau_p\) relative to \(U\), and \(q'\) a prime with Frobenius \(g\varphi_p\) relative to \(U'\). If \(p\) and \(q\) are not both inert, then
\[
 A(U\cup\{q_1\})+A(U\cup\{q\})=A(U'\cup\{q'\}).
\]
If \(p\) and \(q\) are both inert, then \(A(U\cup\{q_1\})=A(U\cup\{q\})\).
\item If \(U\) is matching and \(q,q'\) are inert primes with the same Frobenius in \(K_U\), then \(A(U\cup\{q,q'\})=A(U)\cdot\overline{\#E(\F_q)/2}\) in \(\F_2\). In particular \(A(U\cup\{q,q'\})=A(U)\) when \(q\) is doubly inert.
\end{enumerate}
\end{proposition}

\begin{proof}
If \(U\cup\{q\}\) is admissible, Lemma~\ref{lem:BUq} and \eqref{eq:reading} give \(A(U\cup\{q\})=\Phi_U(\mathrm{Frob}_q)\); otherwise both sides of each assertion vanish. Assertion~(1) follows from Lemma~\ref{lem:factorq}, since \(\Phi_U\) is a trace and hence a class function, and (2) follows from \(\Phi_U(1)=0\).

For (3), we follow the contraction step of Subsection~\ref{subsec:proofrules}. Lemma~\ref{lem:inertia} holds for \(M_U\), since its proof uses only that the characters of \(\theta B_U\) and of the \(\theta\Theta_J\) with \(p\in J\) are ramified at \(p\), with \(p\)-part the Legendre symbol, which remains true here. Choose a finite quotient through which \(t\) and \(d\) act on \(M_U/2M_U\), dominating the groups attached to \(q\) relative to \(U\) and \(U'\). Chebotarev supplies fresh primes realizing \(g\tau_p\), \(g\varphi_p\) and \(g\tau_p\varphi_p\); use \(q\) itself for \(g\). The same calculation gives
\[
 A(U\cup\{q_1\})+A(U\cup\{q\})=\overline{\nu_p(U'q'_1)}\,A(U'\cup\{q'_1\})+\overline{\nu_p(U'q'_0)}\,A(U'\cup\{q'_0\})
\]
by Lemma~\ref{lem:BUq}, where the right side is zero when \(U\) is opposite and \(p\) is inert, that is, when \(p\) and \(q\) are both inert. Otherwise \(q'_0\) and \(q'_1\) have the same Frobenius relative to \(U'\), because \(\tau_p\) acts trivially on \(\Kt_{U'}\), while \(\kappa=(\frac{U'q'}p)\) changes sign between them, because \(\tau_p\) changes the sign of \(\sqrt{p^*}\). Thus the two values of \(\overline{\nu_p}\) sum to \(1\) by Lemma~\ref{lem:recurrence}, and (1) gives the assertion.

For (4), let \(U\) be matching and \(q\) inert, and put \(b=B_U[1]=2^{-n_U}a(U)\), whose residue is \(A(U)\), and \(\epsilon=\chi_{-4}(q)\). Since \(t(\mathrm{Frob}_q)=T_q\) and \(T_q^2=T_{q^2}+q\langle q\rangle\) in weight two,
\[
 \bigl(t(\mathrm{Frob}_q)^2f_U\bigr)[1]=f_U[q^2]+q\,(\langle q\rangle f_U)[1].
\]
With Shimura's normalization the substitution \(z\mapsto tz\) multiplies the character by \(d\mapsto(\frac td)\) \cite[Proposition~1.3]{Shimura73}, while \(U_{p^2}\) preserves it; so \(B_U\) has character \(\psi_U=\psi\prod_{p\in U}(\frac p\cdot)\), and \(\psi_U(q)=\eta(q)\epsilon(\frac Uq)\) because \(\psi/\eta=\chi_{-4}\). The series \(\Theta_\varnothing\) has character \(\psi\), and \(\theta f\) has, in weight two, the character of \(f\). By Proposition~\ref{prop:wald}(2), \(B_U[q^2]=\eta(q)(a_q-(\frac Uq))b\), and \(\Theta_\varnothing[q^2]=\eta(q)q\). Since \(f_U\equiv H_U\pmod2\) and only \(J=\varnothing\) contributes at the indices \(q^2\) and \(1\), modulo \(2\),
\begin{align*}
 f_U[q^2]+q\,(\langle q\rangle f_U)[1]&\equiv\frac{b\,\eta(q)}2\Bigl(a_q-\bigl(\tfrac Uq\bigr)-q\Bigr)+\frac{b\,\eta(q)}2\,q\Bigl(\epsilon\bigl(\tfrac Uq\bigr)-\epsilon\Bigr)\\
 &=\frac{b\,\eta(q)}2\Bigl(a_q-\bigl(\tfrac Uq\bigr)(1-q\epsilon)-q(1+\epsilon)\Bigr).
\end{align*}
If \(q\equiv3\pmod4\), then \(\epsilon=-1\), \(1-q\epsilon=1+q\equiv0\pmod4\), and the bracket divided by two is congruent to \(a_q/2\) modulo \(2\). If \(q\equiv1\pmod4\), then \(1-q\epsilon\equiv0\pmod4\) and the bracket divided by two is congruent to \(a_q/2-q\equiv a_q/2+1\). In both cases it is congruent to \(\#E(\F_q)/2=(q+1-a_q)/2\). On the other hand, for a fresh prime \(q''\) with the Frobenius of \(q\) in a finite quotient through which \(t\) acts on \(M_U/2M_U\) and which dominates \(\Gal(K_U/\Q)\), the left side is congruent to \(f_U[qq'']\), whose residue is \(A(U\cup\{q,q''\})\) by Lemma~\ref{lem:BUq}.

Finally \(A(U\cup\{q,q'\})=A(U\cup\{q,q''\})\): relative to the opposite set \(U\cup\{q\}\), the inert primes \(q'\) and \(q''\) have the same Frobenius in \(K_U\) and may differ only on \(\sqrt q\), and that change does not affect the symbol by~(3).
\end{proof}

\subsection{Graph data over \texorpdfstring{\(F\)}{F}}\label{subsec:qgraph}

\subsubsection{Generators}
The local squareclasses of prime generators determine the conditions
at \(\Sigma\), and their pair residues determine the Frobenius data.
The following ray-class realization will also be used in the cubic case.

\begin{lemma}\label{lem:ray-realization}
Let \(F/\Q\) be finite Galois and let \(\Sigma\) contain \(2\)
and its ramified primes.
Prescribe local squareclasses at the places above \(\Sigma\) and
at the real places, and unit squareclasses at finitely many other
prime ideals. Infinitely many degree-one prime ideals
\(\mathfrak q\) admit a generator of
\(\mathfrak q\mathcal O_\Sigma\) with these classes.
\end{lemma}
\begin{proof}
Weak approximation gives \(a\in F^\times\) with the prescribed
classes, a unit at the other specified primes. Choose a modulus
\(\mathfrak f\) supported on these primes and those above
\(\Sigma\), including the real places, so deep that its principal
ray subgroup consists of local squares. Let \(\mathfrak b\) be
the part of \((a)\) outside \(\Sigma\). It is prime to
\(\mathfrak f\). For a prime \(\mathfrak q\) in its ray class,
\[
 \mathfrak q=(\beta)\mathfrak b,\qquad
 \beta\equiv1\pmod{\mathfrak f},\qquad
 \pi=\beta a,
\]
with \(\beta\) positive at the real places. Then
\(\pi\mathcal O_\Sigma=\mathfrak q\mathcal O_\Sigma\), with
all the required classes. Chebotarev supplies infinitely many such
prime ideals. Those of degree greater than one have density zero;
degree one implies complete rational splitting because \(F/\Q\)
is Galois \cite[Chapters~V--VI]{MilneCFT}.
\end{proof}

For the quadratic field \(F\) of this section, let \(\sigma\) be its non-trivial automorphism and put
\[
 \mathcal Y=\prod_vF_v^\times/(F_v^\times)^2,
\]
the product over the places of \(F\) above \(\Sigma\) and the real places. For a prime \(\mathfrak q\) of \(F\) and \(x\in F^\times\) a unit at \(\mathfrak q\), \([x/\mathfrak q]\in\F_2\) is the additive quadratic residue symbol, and \((x,y)_v\in\F_2\) is the additive Hilbert symbol.

\begin{definition}\label{def:type}
For a split prime \(p\notin\Sigma\), a \emph{chosen place} is a place \(\mathfrak p\) of \(F\) above \(p\), and a \emph{generator} is a generator \(\pi_p\) of \(\mathfrak p\mathcal O_\Sigma\); for an inert prime \(p\) put \(\mathfrak p=p\mathcal O_F\) and \(\pi_p=p\). The \emph{type} \(\typ(p)\in\mathcal Y\) is the class of \(\pi_p\), that is, the collection of its local square classes at the places of \(\mathcal Y\).
\end{definition}

The chosen generator has valuation one at \(\mathfrak p\) and zero valuation at all other primes outside \(\Sigma\). Its type includes the choice of generator; this choice is fixed in a realization.

\begin{lemma}\label{lem:kummer}
For a finite set \(U\) of primes outside \(\Sigma\), \(\Kt_U\)
is generated over \(F\) by the square roots of
\(\mathcal O_\Sigma^\times\), of \(\pi_p\) for \(p\in U\),
and of \(\sigma\pi_p\) for split \(p\in U\). The prime
generators are independent modulo
\(F(\sqrt{\mathcal O_\Sigma^\times})\).
\end{lemma}

\begin{proof}
Kummer theory identifies \(\Kt_U\) with the classes of even
valuation outside \(\Sigma\cup U\). Divide by prime generators
to remove the valuations at \(U\). The remaining principal ideal
is a square in \(\mathcal O_\Sigma\), so its class is a unit
times a square, since this ring has trivial class group.
Independence follows from the valuations above \(U\).
\end{proof}

\begin{definition}\label{def:qdatum}
Let \(U\) be a finite set of primes outside \(\Sigma\), with chosen places and generators. The \emph{forward data} of \(q\in U\) are:
\begin{enumerate}
\item if \(q\) is split, the symbols \([x/\mathfrak q]\) for \(x\in\mathcal O_\Sigma^\times\), for \(x=\pi_p,\sigma\pi_p\) with \(p\in U\smallsetminus\{q\}\) split, and for \(x=p\) with \(p\in U\) inert;
\item if \(q\) is inert, the Legendre symbols \((\frac wq)\), \(w\in\mathcal W_\Sigma\), and \((\frac rq)\), \(r\in U\smallsetminus\{q\}\).
\end{enumerate}
The \emph{datum} \(\Dat(U)\) consists of the splitting of each \(q\in U\) in \(F\), the types of the split primes, and the forward data of all \(q\in U\). Two sets \emph{realize the same datum} if a bijection preserves splitting in \(F\), types and forward data, with the chosen generators corresponding.
\end{definition}

For an old split prime \(p\), write
\(\pi_p\sigma\pi_p=p\,u_p\), with \(u_p\) a rational
\(\Sigma\)-unit. Its rational data at a split \(q\) are therefore
\begin{equation}\label{eq:quadratic-norm-bit}
 [p/\mathfrak q]=[\pi_p/\mathfrak q]
 +[\sigma\pi_p/\mathfrak q]+[u_p/\mathfrak q].
\end{equation}
The last term is fixed by the type of \(q\). Thus these rational
data, used also in mixed contractions, are determined by the forward
data and the types.

By Lemma~\ref{lem:kummer}, the forward data of \(q\) determine its Frobenius relative to \(U\smallsetminus\{q\}\) at the chosen place. For a split vertex, Hilbert reciprocity gives \([x/\mathfrak q]=\sum_v(x,\pi_q)_v\) for \(x\in\mathcal O_\Sigma^\times\); these unit residues are therefore fixed by the type. For an inert vertex, the data at \(\mathcal W_\Sigma\) form its rational label.

The reverse pair data are determined by reciprocity. For split \(p,q\),
\[
 [\pi_q/\mathfrak p]+[\pi_p/\mathfrak q]=\sum_v(\pi_p,\pi_q)_v,\qquad[\sigma\pi_q/\mathfrak p]+[\sigma\pi_p/\mathfrak q]=\sum_v(\pi_p,\sigma\pi_q)_v
\]
where the sums run over the places of \(\mathcal Y\). If one prime is inert, use rational reciprocity and \([\pi_p/q\mathcal O_F]=\bigl[(\frac{N_{F/\Q}\pi_p}q)=-1\bigr]\). The rational label of a split prime is also determined by its type.

A split pair carries two oriented variables
\(x^+=[\pi_p/\mathfrak q]\) and \(x^-=[\sigma\pi_p/\mathfrak q]\).
Every other pair carries one rational bit. A change contracts the pair
unless both vertices are inert, in which case the symbol is unchanged.
Replacing the chosen place and generator of \(q\) by \(\sigma(\mathfrak q)\) and \(\sigma\pi_q\) changes the type to \(\sigma\typ(q)\), interchanges \(x^+\) and \(x^-\) on the pairs at \(q\), and does not change the symbol.

\begin{lemma}\label{lem:realizeq}
Let \(U\) be a finite set of primes outside \(\Sigma\) and \(Z\) a finite set of primes. Let new vertices be prescribed one at a time, each split or inert, with a type if it is split and a rational label \(\mu\) with \(\mu(d)=1\) if it is inert, and with forward data relative to \(U\) and the earlier new vertices, the data of a split vertex at \(\mathcal O_\Sigma^\times\) being those determined by its type. Then there are distinct primes outside \(\Sigma\cup U\cup Z\), with chosen places and generators for the split ones, realizing these data; the vertices of \(U\) keep their choices.
\end{lemma}

\begin{proof}
Choose new vertices successively. For an inert vertex, its rational
label and pair bits prescribe quadratic characters at distinct prime
radicals; Dirichlet's theorem realizes them.

For a split vertex, prescribe its type at the support. At each old
split prime prescribe the two residues of its generator at the
conjugate places; Hilbert reciprocity converts these into the required
forward residues. At an old inert prime, prescribe the residue of
the generator there; its norm gives the required rational bit.
Lemma~\ref{lem:ray-realization} realizes all these conditions by a
degree-one prime, avoiding the finite excluded set. The rational
bits at old split primes are then \eqref{eq:quadratic-norm-bit}.
The data at \(\Sigma\)-units and the reverse data follow from
reciprocity, so no further independent conditions are imposed.
\end{proof}

\subsubsection{Graph rules}
For a datum \(\Dat\) and a variable \(x\) on a pair \(\{p,q\}\), let \(\Dat^x\) be the datum with \(x\) and its reciprocity partner changed. If \(x=x^+\) and \(p,q\) are split, let \(\Dat/x\) be the datum in which \(p,q\) are replaced by a split vertex \(c\) of type \(\typ(p)+\typ(q)\), whose forward data are the sums of those of \(p\) and \(q\), and at which every other vertex has the sums of its data at \(p\) and \(q\); this is the datum of a vertex with generator \(\pi_p\pi_q\). For \(x=x^-\), reorient \(q\) first. If one of \(p,q\) is inert, \(\Dat/x\) replaces them by an inert vertex with the sum of their rational labels and of their rational data; the type of the split vertex is forgotten, because only rational Frobenius data enter at an inert vertex. 

\begin{proposition}\label{prop:qgraph}
Let \(\Dat\) be a datum of total rational label zero.
\begin{enumerate}
\item \(A\) depends only on \(\Dat\).
\item For a variable \(x\) on a pair of vertices not both inert, \(A(\Dat)+A(\Dat^x)=A(\Dat/x)\). For a pair of inert vertices, \(A(\Dat^x)=A(\Dat)\).
\item If a split vertex has all forward data zero, then \(A(\Dat)=0\).
\item If \(\Dat\) has an even number of inert vertices and \(q,q'\notin\Dat\) are doubly inert vertices with equal forward data relative to \(\Dat\), then \(A(\Dat\cup\{q,q'\})=A(\Dat)\), whatever the data between \(q\) and \(q'\).
\end{enumerate}
\end{proposition}

\begin{proof}
For (1), the successive replacement proof of Subsection~\ref{subsec:proofrules} applies, with Proposition~\ref{prop:qrules}(1) for the last vertex and Lemma~\ref{lem:realizeq} for realization: the bridging primes are prescribed with the splitting and the types of the original vertices and with forward data relative to both realizations, which are first made disjoint, and reciprocity makes each intermediate set a realization.

For (2), realize \(\Dat\) by a set \(U\cup\{q\}\) with \(q\) as the last vertex and \(p\in U\). The variable \(x^+=[\pi_p/\mathfrak q]\) is changed by \(\tau_p\), taken at the chosen place of \(p\), which acts by \(-1\) on \(\sqrt{\pi_p}\) and trivially on the other generators of \(\Kt_U\); in the other cases \(\tau_p\) changes \(\sqrt p\) and the rational bit. Proposition~\ref{prop:qrules}(3) gives the identity with a prime \(q'\) whose Frobenius relative to \(U'=U\smallsetminus\{p\}\) is that of \(q\) times \(\varphi_p\), so that its forward data are the sums. By Proposition~\ref{prop:qrules}(1), the value \(A(U'\cup\{q'\})\) depends only on this Frobenius.

If \(p\) and \(q\) are split, the type \(\typ(p)+\typ(q)\) gives the required values on the square roots of \(\Sigma\)-units: for \(u\in\mathcal O_\Sigma^\times\), bilinearity of the Hilbert symbols gives \(\sum_v(u,\pi_p\pi_q)_v=[u/\mathfrak p]+[u/\mathfrak q]\). Lemma~\ref{lem:realizeq} therefore realizes this Frobenius by a prime \(q'\) of type \(\typ(p)+\typ(q)\), and then \(\Dat(U'\cup\{q'\})=\Dat/x\) by reciprocity. If one of \(p,q\) is inert, \(q'\) is inert and only its rational data enter.

Assertions (3) and (4) are Proposition~\ref{prop:qrules}(2) and~(4), with Lemma~\ref{lem:good}; the data between \(q\) and \(q'\) are immaterial by~(2).
\end{proof}

\begin{lemma}\label{lem:oforest}
Let \(\Dat\) be a datum and \(S\) a set of variables, viewed as edges of a multigraph on the vertices.
\begin{enumerate}
\item If \(S\) contains a cycle, including two variables on one pair, then \(\Delta_SA(\Dat)=0\).
\item If \(S\) is a forest in which some tree contains two inert vertices, then \(\Delta_SA(\Dat)=0\).
\item Otherwise \(\Delta_SA(\Dat)=A(\Dat/S)\), where \(\Dat/S\) contracts each tree along its variables, the orientation of each split vertex being dictated by the path to a root.
\end{enumerate}
\end{lemma}

\begin{proof}
Assertion (3) follows by induction as in Lemma~\ref{lem:forest}; the image in \(\Dat/x\) of a variable sharing an endpoint with \(x\) is again a variable, since the data of the merged vertex are sums. For (2), contract along the tree: at the first merger of two parts each containing an inert vertex, both parts are inert vertices of the contracted datum, and the difference vanishes by Proposition~\ref{prop:qgraph}(2). For (1), first realize the two oriented variables on a split
pair \(\{p,q\}\), with \(q\) last. After the first inertia difference,
the coefficient form is attached to the set with \(p\) omitted,
so its level is prime to the entire rational prime \(p\).
Its Galois representations are therefore unramified at both places
above \(p\).
The second oriented difference is zero. A repeated difference in
the same Boolean variable is also zero. For a longer cycle,
contract a path joining the endpoints of a closing variable, using
(2) if that path already gives zero. The closing difference then
reduces to one of these cases.
\end{proof}

\subsection{Networks and packets}\label{subsec:qnetworks}

\subsubsection{Networks}
An \emph{auxiliary vertex} is a split vertex of type \(0\). Networks are defined as in Subsection~\ref{sec:graphs}, with these conventions: on a pair of auxiliary vertices the chosen variable is \(x^+\) for a fixed ordering; on a pair formed by an auxiliary vertex \(i\) and a split terminal \(a\), the selected variable is \([\pi_i/\mathfrak a]\) for alignment \(\epsilon_a=+\), and \([\sigma\pi_i/\mathfrak a]\) for \(\epsilon_a=-\), regardless of the ordering of the pair; on a pair with an inert terminal it is the rational bit; all other variables at auxiliary vertices are zero. In the datum with all chosen variables zero, an auxiliary vertex has all forward data zero, by type \(0\), and the data of the other vertices at it are zero by reciprocity, since its generator is a local square at every place of \(\mathcal Y\). Its rational label is zero for the same reason.

\begin{lemma}\label{lem:qnetwork}
Lemmas~\ref{lem:network} and~\ref{lem:selector} and Theorem~\ref{thm:address} hold for the data of this section, with auxiliary vertices of type \(0\) and the conventions above, and with the following changes. For terminals \(a\ne b\) let \(\Delta^\epsilon_{ab}\) be the difference along the variable on \(\{a,b\}\) of relative orientation \(\epsilon_a\epsilon_b\) if \(a\) and \(b\) are split, along its rational bit if exactly one of them is inert, and \(0\) if both are inert. In Lemma~\ref{lem:network}, \(\Delta_{ab}\) is replaced by \(\Delta^\epsilon_{ab}\), and \(\mathfrak d\) is generated by the differences along the terminal variables. In Lemma~\ref{lem:selector}, \(a_0\) and \(a_1\) are not both inert, and for a variable \(z\) on \(\{a_0,a_1\}\) the alignments are chosen so that \(\Delta^\epsilon_{a_0a_1}=\Delta_z\); then \(O_{x^*}\equiv\Delta_z\pmod{\mathfrak d^2}\).
\end{lemma}

\begin{proof}
The chosen variables give label-zero auxiliary vertices with all
forward data zero in the initial datum. Expand as in
Lemma~\ref{lem:network}, using Lemma~\ref{lem:oforest} in place of
Lemma~\ref{lem:forest}. Cycles, components without terminals, and
components with two inert terminals vanish. Every other component
contracts its terminals, with orientation the product of the alignments
along its paths. The constant and linear terms are therefore the same
matrix-tree counts, with \(\Delta_{ab}\) replaced by
\(\Delta^\epsilon_{ab}\). The same expansion proves composition;
the Laplacian calculation in Lemma~\ref{lem:selector} is unchanged.

The terminal polynomial is independent of inert--inert bits, so a
highest-degree monomial has a variable with permitted endpoints.
Thus the degree descent of Theorem~\ref{thm:address} applies.
Lemma~\ref{lem:realizeq} realizes the resulting data.
\end{proof}

\subsubsection{Packets}
To recover a datum with symbol one, we first contract suitable split vertices into the unwanted inert vertices, then remove pairs with equal labels and equal forward data.

A rational label \(\mu\) is that of a split prime if \(\mu(d)=0\), and every such label is the label of some split prime, the label of a split prime being determined by its type. For each rational label needed for the additional split vertices below, choose the type of one split prime with that label; the split vertices of \(\varsigma\) and their copies keep their original types. If \(F\ne F'\), fix a label \(\xi\) with \(\xi(d)=0\) and \(\xi(d')=1\), the label of a prime split in \(F\) and inert in \(F'\). If \(F=F'\), every inert prime is doubly inert, and no correction is needed below.

\begin{lemma}\label{lem:qpacket}
Let \(\varsigma\) be a datum of total rational label zero with \(A(\varsigma)=1\), with split vertices \(v_1,\dots,v_r\) and inert vertices \(u_1,\dots,u_s\), \(s\ge1\), of rational labels \(\mu_1,\dots,\mu_s\). There is a finite list \(\mathcal P\) of vertices, each split with a type or inert with a rational label, containing those of \(\varsigma\), with total rational label zero, such that the following holds. Let \(\mathcal T\) consist, for each non-zero \(\ell\in(\F_2^b)^*\), of two packets of primes realizing \(\mathcal P\) with ramification form \(\ell\), chosen by Lemma~\ref{lem:realizeq} with arbitrary forward data. Then every \(\mathcal T_x\) has rational label zero, and \(\mathscr P_x\ne0\) for every \(x\ne0\).
\end{lemma}

\begin{proof}
\emph{Correcting the labels.}
Evaluating the total label at \(d\) shows that \(s\) is even, so \(s\ge2\). Put \(\mu'_j=\mu_j+\xi\) if \(u_j\) is not doubly inert, and \(\mu'_j=\mu_j\) otherwise; each \(\mu'_j\) is the label of a doubly inert prime. Let \(m(\lambda)=\#\{j:\mu'_j=\lambda\}\), let \(O\) be the set of labels \(\lambda\) with \(m(\lambda)\) odd, and fix a doubly inert label \(\lambda_*\). Let \(\mathcal P\) contain the vertices of \(\varsigma\) and the following
split vertices: a vertex \(d_j\) of label \(\xi\) for each \(u_j\) not
doubly inert; a vertex \(z_\lambda\) of label \(\lambda_*+\lambda\) for every
\(\lambda\in O\); and a last vertex of label \(\sum_j\mu_j\).

Since \(s\) is even and \(\sum_\lambda m(\lambda)=s\), the set \(O\) has even cardinality, and \(\sum_{\lambda\in O}\lambda=\sum_j\mu'_j\). The labels of the \(d_j\) and \(z_\lambda\) therefore sum to \(\sum_j\mu_j\), and the last label makes the total zero, as \(\varsigma\) has total label zero. All added labels vanish at \(d\), the last because \(s\) is even, so they are labels of split primes.

\emph{Choosing the clusters.}
Fix \(x\ne0\). The number of packets present at \(x\) is
\(k=2\cdot2^{b-1}=2^b\). Doubling the packets makes \(k\) even
also for \(b=1\), so the unwanted inert clusters can be paired.
Fix two packets, \(B_0\) and \(B_1\). Form clusters as follows. The split vertices \(v_1,\dots,v_r\) of \(B_0\) are singletons, and so are \(u_1,\dots,u_{s-1}\) of \(B_0\). In every packet \(B\ne B_0\) present at \(x\), each \(u_j\) that is not doubly inert forms a cluster with the \(d_j\) of \(B\), of label \(\mu'_j\); the doubly inert \(u_j\) of \(B\) are singletons. For each \(\lambda\in O\) choose \(j\) with \(\mu'_j=\lambda\) and add \(z_\lambda\) of \(B_1\) to the cluster of \(u_j\) in \(B_1\), which then has label \(\lambda_*\); these moves use distinct clusters because the labels \(\lambda\) are distinct. Every remaining split vertex joins the cluster of \(u_s\) of \(B_0\). Each cluster contains at most one inert vertex, and the clusters outside \(B_0\) are inert with doubly inert labels. For \(\lambda\ne\lambda_*\) the number of clusters of label \(\lambda\) outside \(B_0\) is
\[
 (k-1)\,m(\lambda)-\mathbf 1_{\{\lambda\in O\}},
\]
which is even, as \(k-1\) is odd. The total number of clusters outside \(B_0\) is \((k-1)s\), which is even, so the number of clusters of label \(\lambda_*\) is even too. Pair the clusters of equal labels.

\emph{Recovering the seed.}
Let \(S\) consist of the rational variables joining each inert vertex to the split vertices of its cluster, a union of stars. By Lemma~\ref{lem:oforest}(3), \(\Delta_S\mathscr P_x=A(\Dat(\mathcal T_x)/S)\), a function of the variables outside \(S\). The contracted datum consists of the split vertices of \(\varsigma\), with their types; inert vertices of labels \(\mu_1,\dots,\mu_{s-1}\); a last inert vertex, whose label is determined by the total to be \(\mu_s\); and the pairs. An inert vertex carries only a rational label, so the types of the split vertices absorbed into it are forgotten. The variables of the contracted datum are sums of the variables outside \(S\) over disjoint sets of pairs, so they can be assigned freely. Assign them so that the vertices coming from \(B_0\) carry the data of \(\varsigma\), and the two members of each pair have equal forward data relative to all other vertices. Removing the pairs one at a time by Proposition~\ref{prop:qgraph}(4) gives \(A(\Dat(\mathcal T_x)/S)=A(\varsigma)=1\) at this assignment. Hence \(\mathscr P_x\ne0\).
\end{proof}

\begin{proof}[Proof of Theorem~\ref{thm:families} in case~\CaseII{}]
By Lemma~\ref{lem:depth} there is an admissible \(U\) with \(A(U)=1\); let \(\varsigma=\Dat(U)\). If \(\varsigma\) has no inert vertex, adjoin two doubly inert primes with equal forward data relative to \(\varsigma\); by Proposition~\ref{prop:qgraph}(4) the symbol is unchanged. With this \(\varsigma\), Lemmas~\ref{lem:qpacket} and~\ref{lem:qnetwork} give a family configuration \(\mathcal U\), whose auxiliary vertices have type \(0\), such that \(A(\mathcal U_x)=1\) for every \(x\ne0\). Auxiliary vertices have rational label zero, so every \(\mathcal U_x\) has total rational label zero and is admissible by Lemma~\ref{lem:filter}. Take \(P\) to be the set of primes of \(\mathcal U\) and \(\alpha_q\) their ramification forms.
\end{proof}

\Needspace{7\baselineskip}
\begin{remark}
The pair-removal identity is what allows the construction to start from an arbitrary non-zero seed, including one with inert vertices. Contraction with a split prime preserves the number of inert vertices, while changing an inert--inert bit leaves the symbol unchanged. Removing the unwanted inert vertices therefore uses the separate coefficient identity of Proposition~\ref{prop:qgraph}(4).
\end{remark}

\section{The cyclic cubic case}\label{sec:cyclic-cubic}

Theorem~\ref{thm:cubic-E} proves Theorem~\ref{thm:families} in
case~\CaseIII{}. The principal steps are the least-weight transfer
of Proposition~\ref{prop:cubic-minimal-transfer} and the disjoint selection
of Lemma~\ref{lem:cubic-disjoint}.

Suppose \(F=\Q(E[2])\) is cyclic cubic and \(E\) has root number
\(+1\). Take \(\Sigma\) as fixed in Section~\ref{sec:symbol}, so
\(\mathcal O_\Sigma\) has trivial class group. Write \(V=E[2]\). Primes outside \(\Sigma\)
are either split or inert in \(F\). Their weights are one and
zero, respectively. For an admissible set \(T\),
\[
 c(h(T))=4^{\,\#\{p\in T:p\text{ split}\}},\qquad
 w(T)=\#\{p\in T:p\text{ split}\}.
\]
Let \(A\) be the normalized symbol of Section~\ref{sec:symbol},
\[
 A(T)=\overline{2^{-\delta-w(T)}a(T)}\in\F_2.
\]
It is one exactly when \(T\) attains the minimum \(\delta\), and
is zero off the admissible support. Bars denote reduction modulo
the maximal ideal \(\mathfrak m\) of \(\mathcal O_g\).

\begin{theorem}\label{thm:cubic-E}
If \(A(\varnothing)=0\), then for every \(b\ge1\) there is an
admissible quadratic twist family with \(A(U_x)=1\) at every \(x\ne0\).
\end{theorem}

Inert primes have weight zero and carry no local condition of their
own; the symbol depends on them only through their combined quadratic
label (Proposition~\ref{prop:cubic-W0}), and they serve to adjust the
local conditions at split primes. Let \(m\) be the least number of
split primes in a set attaining the minimum. At this weight the symbol
depends on the split primes only through their Frobenius in a fixed
finite extension of \(F\), so a configuration attaining the minimum can
be copied with prescribed residue symbols between distinct primes
(Proposition~\ref{prop:cubic-minimal-transfer}). At each split prime
the symbol satisfies the same three-term relation as the determinant
of the matrix of local Selmer conditions
(Lemma~\ref{lem:cubic-state-relation}). The central comparison turns
attainment of the minimum into a Selmer bound, which limits the number
of singular blocks of this matrix (Proposition~\ref{prop:cubic-nullity}).
Selectors and a maximality argument over disjoint blocks attain the
minimum at every non-zero index, and inert primes realize the required
local conditions in one family (Lemma~\ref{lem:cubic-disjoint},
Proposition~\ref{prop:cubic-calibration}; cf.~\cite[Sections~5.5--5.6,
8--9]{OAG}).

\subsection{Integral translates and their trace operators}

In this subsection \(U\) is a finite set of split primes. It need not
be admissible. Use the half-integral-weight
operators \((U_pf)[n]=f[pn]\) and
\((V_pf)[n]=f[n/p]\). Define
\begin{equation}\label{eq:cubic-B}
 B_U=2^{-\delta-|U|}\prod_{p\in U}
          (U_p-\eta(p)V_p)g,
 \qquad f_U=\theta B_U.
\end{equation}
Here \(\eta=\psi\chi_{-4}\), as in Subsection~\ref{subsec:waldspurger},
and the weight-two character is
\(\psi_U=\psi\prod_{p\in U}(\frac p\cdot)\).
For a coefficient outside the support of \(g\), the valuation is
understood to be infinite.

\begin{lemma}\label{lem:cubic-B-integral}
The form \(B_U\) is integral. Write an index prime to \(\Sigma\) as
\(n=Ir k^2\), where \(I\subseteq U\), \(r\) is square-free and prime
to \(U\). Then
\[
 v_2(B_U[n])\ge w(r).
\]
For a square-free product \(r\) of inert primes,
\[
 \overline{B_U[r]}=A(U\cup\{q:q\mid r\}).
\]
\end{lemma}
\begin{proof}
The shell recurrence of Subsection~\ref{subsec:hecke} applies at every
\(p\in U\), since \(a_p\) is even and
\(a_p\equiv p+1\pmod4\). At a prime in \(I\) the difference operator
contributes an even integer after its root-of-unity phase is removed;
at every other prime its multiplier is integral. The underlying
square-free coefficient is that of \(D=(U\setminus I)r\). Minimality
gives \(v_2(a(D))\ge\delta+|U|-|I|+w(r)\), and the factors at \(I\)
supply \(|I|\) additional powers of two. Division in
\eqref{eq:cubic-B} proves the estimate. Indices divisible by
\(\Sigma\) have zero coefficient. The rationality assertion of
Proposition~\ref{prop:wald}(4) makes these integral valuations.
At \(r\) prime to \(U\), every term containing \(V_p\)
vanishes, so \(B_U[r]=2^{-\delta-|U|}a(Ur)\). This proves the last
assertion, including \(r=1\).
\end{proof}

The weight-two form \(f_U\) is cuspidal, has the single character
\(\psi_U\), and satisfies \(f_U\equiv B_U\pmod2\).
Take its Hecke span over the fraction field of \(\mathcal O_g\),
and intersect with the integral \(q\)-expansions. Denote this lattice
by \(L_U\). A finite coefficient map embeds \(L_U\) in a finite free
\(\mathcal O_g\)-module, so \(L_U\) is finite free. It is Hecke
stable. It is also saturated: if every coefficient of \(g'\in L_U\)
is divisible by a uniformizer, then dividing \(g'\) by that
uniformizer gives another element of \(L_U\). Thus coefficientwise
reduction detects equality in \(L_U/\mathfrak mL_U\).

Let \(\mathscr A_U\) be the image of the Hecke algebra modulo
\(\mathfrak m\) on the cyclic module generated by
\(\overline f_U\) in \(L_U/\mathfrak mL_U\). The construction
of Galois trace operators in Subsection~\ref{subsec:trace-operators} applies to
this lattice. Denote the resulting trace by \(T_U(\sigma)\). The
nebentypus has \(2\)-power order, so its determinant reduces to one.
Thus \(\mathscr A_U\) is a finite commutative \(\F_2\)-algebra and
\begin{equation}\label{eq:cubic-trace-identity}
 T_U(a)T_U(b)=T_U(ab)+T_U(ab^{-1}),\qquad T_U(1)=0.
\end{equation}
The algebra need not be reduced. If the cyclic module is zero, every
symbol tested on it is zero and the assertions below are immediate.

\begin{lemma}\label{lem:cubic-trace-types}
On this cyclic module,
\[
 T_U(\sigma)=0\quad(\sigma|_F=1),\qquad
 T_U(\sigma)^2=1\quad(\sigma|_F\ne1).
\]
\end{lemma}
\begin{proof}
Fix a coefficient index \(n\). Chebotarev permits each trace test to
be represented by a prime \(q\) outside \(\Sigma\cup U\) and prime
to \(n\), with the same conjugacy class in a finite quotient through which
the reduced traces and \(V\) factor. If \(q\) is split, the
coefficient of \(T_qf_U\) at \(n\) is, modulo the maximal ideal,
\(B_U[qn]\). It vanishes by Lemma~\ref{lem:cubic-B-integral}, since
the fresh split factor adds weight one. This holds at every index.
Saturation and commutation with the Hecke operators give the first
operator identity on the cyclic module.

For an inert prime \(q\), the weight-two Hecke formula gives
\[
 (T_q^2f_U)[n]=f_U[q^2n]+q\psi_U(q)f_U[n]
 \quad(q\nmid n).
\]
The good-prime recurrence for \(B_U\) is
\[
 B_U[q^2n]=\eta(q)\left(a_q-
       \Bigl(\frac{n\prod_{p\in U}p}{q}\Bigr)\right)B_U[n].
\]
Indeed, the \(T(q^2)\)-eigenvalue of \(B_U\) remains \(\eta(q)a_q\),
while the middle term uses the character \(\psi_U\).
For inert \(q\), both \(a_q\) and the displayed quadratic symbol
are odd. The first term of the weight-two formula therefore
vanishes after reduction. The second term reduces to \(f_U[n]\). Choosing \(q\) separately
for each coefficient index proves \(T_U(\sigma)^2\overline f_U=\overline f_U\); by
commutation this holds on the entire cyclic module.
\end{proof}

\subsection{Aggregation of inert primes}

\begin{lemma}\label{lem:cubic-character-algebra}
Let a finite group \(\Gamma\) map onto \(C_3\), and let
\(T:\Gamma\to\mathscr A\) be a conjugacy-invariant function to a
commutative \(\F_2\)-algebra satisfying
\eqref{eq:cubic-trace-identity}. Suppose \(T\) is zero on the kernel
and \(T(g)^2=1\) off the kernel. Then there is a character
\(\alpha:\Gamma\to\mathscr A^\times\), with \(\alpha^2=1\), such
that \(T(g)=\alpha(g)\) off the kernel.
\end{lemma}
\begin{proof}
For \(u\) in the kernel, the trace identity gives
\(T(gu^2)=T(g)\). The subgroup generated by these squares is
characteristic in the kernel, hence normal in \(\Gamma\); the
function descends to its quotient. Its kernel \(H\) is now
elementary abelian. If \(g\) lifts a generator of \(C_3\), then
\(g^3\in H\), so \(g^6=1\) and \(\varrho=g^4\) has order three
with the same image. The trace identity
first gives \(T(\varrho)=T(\varrho^2)\), and then
\(1=T(\varrho)T(\varrho^2)=T(\varrho)\). For \(a\in H\), set
\[
 \alpha(a)=T(a\varrho)=T(a\varrho^2).
\]
The equality follows by applying the identity to \(a,\varrho\).
Write \(\varrho(b)=\varrho b\varrho^{-1}\). Conjugacy by
\(b\in H\) makes this function invariant under
multiplication by \(b\varrho(b)^{-1}\), so it factors through
\(H/[H,\varrho]\). Therefore
\[
 \alpha(a)\alpha(b)
 =T(a\varrho)T(b\varrho)
 =T(a\varrho(b)\varrho^2)+T(ab^{-1})
 =\alpha(ab).
\]
It is invariant under \(\varrho\), and its values square to one. Extend it
by \(\alpha(\varrho)=1\), and pull it back to the original group.
These identities hold in \(\mathscr A\), including its nilpotents.
\end{proof}

The character can contain information that disappears on passing
to residue fields. For example, in \(\F_2[\epsilon]/(\epsilon^2)\)
the expression \(1+\epsilon\chi(g)\) is an order-two character
whenever \(\chi\) is an additive quadratic character. All its values
reduce to one modulo \(\epsilon\), while the full algebra still
records \(\chi\).

Let \(K_U\) be the maximal elementary abelian \(2\)-extension of
\(\Q\) unramified outside \(\Sigma\cup U\) and infinity. The quadratic
label of a set \(\mathcal R\) of fresh primes is the product of its Frobenius
elements in \(\Gal(K_U/\Q)\); the empty product is the identity.

\begin{proposition}\label{prop:cubic-W0}
Fix a set \(U\) of split primes. If \(\mathcal R\) and \(\mathcal R'\) are finite sets
of inert primes outside \(\Sigma\cup U\) with the same
quadratic label in \(K_U\), then
\[
 A(U\cup \mathcal R)=A(U\cup \mathcal R').
\]
In particular, suppose \(U\) and \(U\cup \mathcal R\) are admissible, and put
\(r^*=\prod_{q\in \mathcal R}q^*\). Then \(A(U\cup \mathcal R)\) depends only on \(U\)
and the symbols \((r^*/p)\), \(p\in U\). If these symbols are all one,
\[
 A(U\cup \mathcal R)=A(U).
\]
For an admissible set \(\mathcal R\) consisting entirely of inert primes,
\(A(\mathcal R)=A(\varnothing)\).
\end{proposition}
\begin{proof}
Apply Lemma~\ref{lem:cubic-character-algebra} to a common finite
quotient for the trace operators of Lemma~\ref{lem:cubic-trace-types}
and the action on \(V\). It gives an order-two character
\(\alpha_U\) into \(\mathscr A_U^\times\), unramified outside
\(\Sigma\cup U\), which agrees with the traces at cubic elements.
It therefore factors through \(K_U\).

For distinct fresh primes, the product of the weight-two Hecke
operators reads the coefficient at their square-free product. Thus,
with \(r=\prod_{q\in \mathcal R}q\),
\[
 A(U\cup \mathcal R)
 =\left[\left(\prod_{q\in \mathcal R}T_U(\operatorname{Frob}_q)\right)
                   \overline f_U\right][1]
 =\left[\alpha_U\left(\prod_{q\in \mathcal R}
                     \operatorname{Frob}_q\right)\overline f_U\right][1].
\]
The empty product gives \(\overline{B_U[1]}=A(U)\). This proves the
first assertion.

In the last assertions, \(r^*\) is positive and locally square at
\(\Sigma\), since both sets are admissible. Positivity fixes the coordinate at \(-1\), and
\(r^*\equiv1\pmod8\) fixes that at \(2\). Reciprocity fixes the
coordinates at the odd primes of \(\Sigma\). The remaining
coordinates of \(K_U\) are
exactly the Legendre symbols against \(p\in U\). If they too are
trivial, the character value is the identity operator.
\end{proof}

Thus a collection of inert primes contributes only its combined
quadratic label. When \(A(\varnothing)=0\), every set attaining
the minimum contains a split prime.

\subsection{Local conditions and the Selmer bound}

Write the Weil pairing additively, in a basis of \(V\) with matrix
\(J=\left(\begin{smallmatrix}0&1\\1&0\end{smallmatrix}\right)\),
and lift this basis to \(E[4]\). At a split prime \(p\), choose
a Frobenius lift \(\varphi_p\) with
\(\chi_{p^*}(\varphi_p)=1\). Its action on \(E[4]\) is
\(I+2D_p\), with \(D_p\) a binary matrix. For a parameter \(h\)
containing \(p^*\), define its unramified twisting bit by
\[
 \kappa_p(h)=\frac{1-\chi_{h/p^*}(\varphi_p)}2\pmod2.
\]

\begin{lemma}\label{lem:cubic-local-kummer}
At a split prime dividing \(h\), the local Kummer image of
\(E^h\), parametrized by its inertia coordinate \(x\in V\), has
unramified coordinate
\begin{equation}\label{eq:cubic-local-kummer}
 (D_p+\kappa_p(h)I)x.
\end{equation}
\end{lemma}
\begin{proof}
Halve the rational two-torsion point \(x\) in \(E[4]\).
Ramified twisting gives Kummer inertia \(x\); Frobenius gives
\(D_px\), with an additional \(x\) when its unramified twisting
character is non-trivial. These classes fill the Kummer image,
which has dimension two.
\end{proof}

Changing \(\kappa_p(h)\) gives two transverse Kummer planes: their
graphs differ by \(I\). At an inert prime, Frobenius has no
invariant on \(V\), so local duality and the Euler characteristic
give \(H^1(\Q_p,V)=0\) (cf.~\cite[Lemmas~2.9--2.10]{MazurRubin}).
Thus inert primes change local conditions at split primes while
contributing no condition at their own places.

There is also a fixed bound at the minimum. Since \(E(\Q)[2]=0\),
Lemma~\ref{lem:selmerbound} gives
\begin{equation}\label{eq:cubic-selmer-bound}
 \dim_{\F_2}\Sel_2(E^{h(U)}/\Q)\le M_E,
 \qquad M_E=C_A+C_E+\beta_E,
 \quad A(U)=1,
\end{equation}
for nonempty admissible \(U\). The existence of such a set makes
\(M_E\ge0\).

\subsection{Global classes and local matrices}

Let \(F_\Sigma^{\rm el}\) be the maximal elementary abelian
\(2\)-extension of \(F\) unramified outside the places above
\(\Sigma\) and infinity. It is finite and Galois over \(\Q\).
A split prime, with a chosen place of \(F\), has an oriented
Frobenius label in \(\Gal(F_\Sigma^{\rm el}/F)\).

Choose a generator \(\sigma\) of \(\Gal(F/\Q)\) whose action on \(V\) is
\(\mathsf R=\left(\begin{smallmatrix}0&1\\1&1\end{smallmatrix}\right)\).
Thus \(\mathsf R^2+\mathsf R=I\) and \(\F_2[\mathsf R]\simeq\F_4\).
The matrices \(D_p\) of the local Kummer conditions need not commute
with \(\mathsf R\); these planes need not be \(\F_4\)-lines. The rotation
algebra \(\F_2[\mathsf R]\) describes the off-diagonal blocks.

For an oriented split prime \(p\), the choice of \(\Sigma\) gives
\(\pi_p\in F^\times\) with divisor \(\mathfrak p\) away from
\(\Sigma\), where \(\mathfrak p\) is its chosen place. Set
\begin{equation}\label{eq:cubic-global-class}
 C_p(x)=\operatorname{cor}_{F/\Q}(\chi_{\pi_p}x)
 \in H^1(\Q,V),\qquad x\in V.
\end{equation}
Here \(\chi_{\pi_p}\in H^1(F,\F_2)\) is the additive Kummer
class, and \(V\) is constant over \(F\). Local restriction--corestriction
shows that the inertia coordinate of \(C_p(x)\) at \(p\) is
\(x\), and it is zero at every other prime outside \(\Sigma\).
Consequently the classes for distinct \(p\) are independent.
Restriction of \eqref{eq:cubic-global-class} to \(G_F\) is
\(\sum_{i=0}^2\chi_{\sigma^i\pi_p}\mathsf R^ix\).
Thus, for \(q\ne p\) split and oriented by \(\mathfrak q\),
\begin{equation}\label{eq:cubic-global-evaluation}
 C_p(x)(\varphi_q)=
 \sum_{i=0}^2[\sigma^i\pi_p/\mathfrak q]\,\mathsf R^ix.
\end{equation}
The brackets are additive quadratic residue symbols.

Define the raw Selmer matrix by
\begin{equation}\label{eq:cubic-raw-matrix}
 B_{qp}(y,x)=\langle C_p(x)(\varphi_q),y\rangle
 -\mathbf1_{p=q}\langle(D_p+\kappa_p(h)I)x,y\rangle.
\end{equation}
Its rows impose precisely the local conditions at the active split
primes. At an active inert prime the local cohomology is zero.

The \emph{type} \(\mathfrak h_p\) of \(p\) records the local squareclasses of
\(\pi_p\) at the places of \(F\) above \(\Sigma\), including
all real signs. This is finite information. It determines the
oriented label in \(F_\Sigma^{\rm el}\): for a quadratic radical
unramified away from \(\Sigma\), Hilbert reciprocity computes its
Frobenius at \(\mathfrak p\) from its pairings with \(\pi_p\)
at those support places. It also determines \(D_p\), since
\(F(E[4])/F\) is elementary abelian and unramified outside
\(\Sigma\). The local maps \(C_p\) at \(\Sigma\) also depend only on the
type. We retain this local information, rather than only the
Frobenius label, because the Selmer comparison requires these maps to
repeat. All types used below are fixed before the family dimension
is chosen.

\begin{proposition}\label{prop:cubic-matrix-realization}
Fix finitely many types occurring at oriented split primes. Any
finite array of these types admits arbitrarily fresh realizations
with the following properties. Cross residues are prescribed freely;
self residues are retained as they occur, modulo the correction below.
\begin{enumerate}
\item Off-diagonal blocks can be prescribed in \(J\F_2[\mathsf R]\), in
one direction for each pair; reciprocity determines the reverse
blocks. For that direction, the rational Legendre bit \((p/q)\)
is independently prescribable; its reverse is fixed by quadratic
reciprocity.
\item For each type, the self block is fixed modulo \(\F_2J\).
There is a primary type, trivial at the support places, whose self
block lies in \(\F_2J\).
\end{enumerate}
The construction extends any fixed finite array and permits the
exclusion of finitely many primes.
\end{proposition}
\begin{proof}
For two distinct split primes, \eqref{eq:cubic-global-evaluation} gives
\begin{equation}\label{eq:cubic-pair-bits}
 B_{qp}=J(e_0+e_1\mathsf R+e_2\mathsf R^2),
\end{equation}
where \(e_i=[\sigma^i\pi_p/\mathfrak q]\). Write \(N_{F/\Q}\pi_p=p\,u_{\mathfrak h_p}\) up to a rational
square, where the squareclass \(u_{\mathfrak h_p}\) is supported
on \(\Sigma\) and infinity. Its class is fixed by the type.
The norm identity gives, in additive Legendre notation,
\[
 e_0+e_1+e_2=[p/q]+[u_{\mathfrak h_p}/q].
\]
The last bit is fixed by the two types. If \(\omega\) denotes the
image of \(\mathsf R\) in \(\F_4\), the map
\[
 (e_0,e_1,e_2)\longmapsto
 (e_0+e_1\omega+e_2\omega^2,\ e_0+e_1+e_2)
\]
is an isomorphism \(\F_2^3\simeq\F_4\oplus\F_2\): the first
component has kernel \(\langle(1,1,1)\rangle\), on which the
second is non-zero. Thus the off-diagonal block and the rational
Legendre bit can be prescribed independently.

Apply Lemma~\ref{lem:ray-realization}, prescribing the type at
the support and the three residue squareclasses at the conjugate
places above every old split prime. The preceding isomorphism
prescribes these three bits from the desired block and norm bit.
Hilbert reciprocity fixes the reverse residues. The lemma gives a
generator with the specified local classes and a single fresh prime
divisor outside \(\Sigma\).

The three forward bits at each old prime are therefore free.
Hilbert reciprocity fixes the reverse bits, with correction only
at the support. Preserving the local type also preserves every
fixed Frobenius label just described. Proceeding one new prime at
a time proves (1).

For the self block the same restriction formula gives the
unramified coordinate \(d_0I+e_1\mathsf R+e_2\mathsf R^2\).
Here \(d_0=\chi_{\pi_p}(\varphi_p)\) uses the chosen Frobenius
lift, while \(e_i=[\sigma^i\pi_p/\mathfrak p]\) for \(i=1,2\). Hilbert reciprocity for
\((\pi_p,\sigma\pi_p)\) gives
\[
 e_1+e_2=s_{\mathfrak h},
\]
where \(s_{\mathfrak h}\) is determined by the support type.
The only off-support contributions are at
\(\mathfrak p\) and \(\sigma\mathfrak p\).
Since \(\mathsf R+\mathsf R^2=1\),
\[
 d_0I+e_1\mathsf R+e_2\mathsf R^2=(d_0+e_1)I+s_{\mathfrak h}\mathsf R^2.
\]
Consequently, writing \(D_{\mathfrak h}=D_p\) for this type,
\[
 B_{pp}\equiv J(s_{\mathfrak h}\mathsf R^2-D_{\mathfrak h})
       \pmod{\F_2J}.
\]
This is the asserted fixed self block. Changing the unit twist
changes precisely the remaining scalar multiple of \(J\).

For the primary type, take a generator locally square at every
support place and positive at every real place. Principal ray
classes supply arbitrarily fresh such primes. Its fixed-support
Frobenius is trivial, so \(D_p=0\), and reciprocity gives
\(e_1=e_2\). Its self block is thus a scalar multiple of \(J\).
\end{proof}

\begin{remark}\label{rem:cubic-spin}
The individual self residues
\((\pi_p/\sigma\mathfrak p)\) and \((\pi_p/\sigma^2\mathfrak p)\), with
the generator fixed, are spin data; compare \cite{FIMR}.
Reciprocity fixes their product; flipping both signs adds \(J\)
to the self block. Proposition
\ref{prop:cubic-matrix-realization} prescribes cross residues
between distinct rational primes and accepts the self residue
as it occurs. In Proposition~\ref{prop:cubic-calibration}, the
offset \(s_p\) is known after the split primes have been chosen.
An inert-prime bit then compensates for it.

\end{remark}

\subsection{Deletion and transfer at the least non-zero weight}

Write \(K_\Sigma=K_\varnothing\), and use additive notation for
quadratic Galois groups. For a prime \(q\notin\Sigma\), let
\(\lambda_\Sigma(q)\) be its Frobenius in
\(\Gal(K_\Sigma/\Q)\). With signed prime discriminants, a finite
prime set \(T\) is admissible exactly when the sum of these labels is
zero. The coordinates at \(-1\) and \(2\) impose positivity
and squareness at \(2\); quadratic reciprocity identifies the
coordinate at an odd \(\ell\in\Sigma\) with \((h(T)/\ell)\).
For a split set \(U\) and
\(\zeta\in\Gal(K_U/\Q)\), define \(A_U(\zeta)\) to be
\(A(U\cup\{r\})\), where \(r\) is an inert prime with quadratic
label \(\zeta\). Such a prime exists: \(F\cap K_U=\Q\), so its
order-three Frobenius can be prescribed independently. The value is
well defined by Proposition~\ref{prop:cubic-W0}. Several inert primes,
or none, give the same value when their product label is \(\zeta\).
Admissibility means that the restriction of \(\zeta\) to
\(K_\Sigma\) cancels the sum of the labels of the primes of \(U\).

\begin{lemma}\label{lem:cubic-deletion}
Let \(p\in U\), put \(U'=U\setminus\{p\}\), and let \(e_p\) be
the inertia generator of \(K_U/K_{U'}\). Then
\begin{equation}\label{eq:cubic-deletion}
 A_U(\zeta+e_p)+A_U(\zeta)
 =A_{U'}\bigl(\zeta|_{K_{U'}}+\operatorname{Frob}_p\bigr).
\end{equation}
The right side keeps the total label at \(\Sigma\) and the local
unit twist at every remaining split prime unchanged.
\end{lemma}
\begin{proof}
Use the trace operators of \(f_U=\theta B_U\). Choose an inert
element \(g\) with label \(\zeta\), a local inertia element
\(\tau_p\) mapping to \(e_p\), and a compatible Frobenius lift
\(\varphi_p\). Every eigensystem in the Hecke span of \(f_U\)
has conductor exponent one at \(p\) and ramified quadratic
nebentypus there. It is a principal series with one unramified
character: a Steinberg constituent would have unramified central
character, and a supercuspidal constituent would have conductor
exponent at least two. The \(U_p\)-eigenvalue is the Frobenius
eigenvalue on the unramified line. Every such eigensystem is
spanned by oldforms \(B_df_j(z)=f_j(dz)\) with \(p\nmid d\),
since both the level and newform conductor have \(p\)-exponent
one. The operator \(U_p\) commutes with these \(B_d\), hence acts
by that scalar on the entire eigensystem. It preserves the Hecke
span and its saturated lattice \(L_U\), since it preserves integral
\(q\)-expansions. Lemma~\ref{lem:inertia} can therefore be reduced
in \(L_U/\mathfrak mL_U\). Since the determinant reduces to one,
taking coefficient one gives
\[
 \bigl[(T(g\tau_p)+T(g))\overline f_U\bigr][1]
 =\bigl[(T(g\tau_p\varphi_p)+T(g\varphi_p))\overline f_U\bigr][p].
\]
These are coefficients in \(L_U/\mathfrak mL_U\); no action of
\(U_p\) on the cyclic submodule is needed. Chebotarev replaces the trace arguments by fresh
primes, avoiding all the coefficient indices.

Both arguments on the right are inert, since \(p\) is split in
\(F\). At a fresh inert prime \(q'\), the first shell is
\[
 B_U[pq']=\eta(p)\nu_p(U'q')\,2^{-\delta-|U'|}a(U'q'),\qquad
 \nu_p(D)=\frac{a_p-(D/p)-1}{2}.
\]
The two fresh primes have the same label in \(K_{U'}\), namely the
restriction of \(\zeta\) plus \(\operatorname{Frob}_p\). Thus their
normalized omitted-\(p\) coefficients have the same reduction, by
Proposition~\ref{prop:cubic-W0}. Their Legendre symbols at \(p\)
are opposite, giving respectively
\[
 \nu_p(D)=\frac{a_p-2}{2}\quad\text{and}\quad\frac{a_p}{2}.
\]
Their difference is one, and all phases reduce to one. This proves \eqref{eq:cubic-deletion}.

For \(q\in U'\), the coordinate of
\(\operatorname{Frob}_p\) at \(\sqrt q\) is
\((\frac qp)=(\frac{p^*}q)\), the contribution of \(p^*\) to
\(\kappa_q\). The right side of \eqref{eq:cubic-deletion}
therefore has the same total label on \(K_\Sigma\) and the same
bits \(\kappa_q\), \(q\in U'\), as the configurations on the left.
\end{proof}

\begin{lemma}\label{lem:cubic-least-divisibility}
Let \(m\) be the least split weight of an admissible set with
symbol one. If \(|U'|=m-1\), then \(B_{U'}/2\) is integral
at every coefficient.
\end{lemma}
\begin{proof}
Indices divisible by \(\Sigma\) have zero coefficient. Write every
other index as \(n=Irk^2\), with \(I\subseteq U'\), \(r\)
square-free and prime to \(U'\). If \(w(r)>0\),
Lemma~\ref{lem:cubic-B-integral} gives
\(\val(B_{U'}[n])\ge1\).
If \(w(r)=0\), the square-free source
\(D=(U'\setminus I)r\) has split weight \(m-1-|I|<m\).
Its coefficient is zero off the admissible support; on that support,
minimality and Remark~\ref{rem:integral-coefficient-valuations} give
\[
 \val(a(D))\ge\delta+m-1-|I|+1.
\]
The shell factors at \(I\) supply \(|I|\) powers of \(2\),
and all remaining factors, including those from \(k^2\), are
integral. After division by \(2^{\delta+m-1}\), the coefficient
is still divisible by \(2\).
\end{proof}

\begin{proposition}\label{prop:cubic-minimal-transfer}
Suppose \(A(\varnothing)=0\), and let \(m\) be the least split
weight of an admissible set with symbol one. Then \(m\ge1\).
At this weight the symbol depends on the split primes only through
the \(G_\Q\)-orbits of their oriented labels in
\(\Gal(F_\Sigma^{\rm el}/F)\), and on the inert primes only
through their total fixed-support label. A configuration with
symbol one consequently has disjoint fresh copies with its types
preserved, every finite choice of pair residues allowed by
Proposition~\ref{prop:cubic-matrix-realization}, and arbitrary
unramified twisting bits. Explicitly, if \(U_i\cup R_i\) are
admissible, \(U_i\) consists of \(m\) split primes, \(R_i\) is
inert, and the multisets of type orbits agree, then
\begin{equation}\label{eq:cubic-type-transfer}
 \{[\mathfrak h_p]:p\in U_1\}
 =\{[\mathfrak h_p]:p\in U_2\}
 \quad\Longrightarrow\quad
 A(U_1\cup R_1)=A(U_2\cup R_2).
\end{equation}
Here brackets denote \(G_\Q\)-orbits and the braces denote multisets.
\end{proposition}
\begin{proof}
The coefficient minimum supplies a symbol equal to one, and
Proposition~\ref{prop:cubic-W0} excludes split weight zero.
Every smaller-weight symbol vanishes. By
Lemma~\ref{lem:cubic-deletion}, changing an inert-product bit at
one of the \(m\) split primes therefore has zero derivative.
Only the restriction of that product to \(K_\Sigma\) matters.

\emph{Traces.} Fix \(|U'|=m-1\). By
Lemma~\ref{lem:cubic-least-divisibility},
\(\widetilde f_{U'}=\theta B_{U'}/2\) is integral. Take its
saturated integral Hecke lattice and reduced cyclic module.
Two split trace operators have product zero on this module.
For each coefficient index, represent their classes by distinct
fresh split primes avoiding that index. The shell estimate gains
two powers of two before the extra division, so every resulting
coefficient vanishes after reduction. Saturation proves the
operator identity on \(\overline{\widetilde f_{U'}}\); commutation extends
it to the cyclic module. The trace identity, applied to \(gu,u\),
now gives \(T(gu^2)=T(g)\) for \(g,u\in G_F\). Thus the split
trace factors through the elementary extension of \(F\) supported
on \(\Sigma\cup U'\).

\emph{Inertia.} Fix a fresh inert prime \(r\) whose support label
restores admissibility, and consider the class function
\[
 \Phi_r(g)=[T(g)T(\operatorname{Frob}_r)\overline{\widetilde f_{U'}}][1],
 \qquad g\in G_F.
\]
At a fresh split Frobenius it reads the weight-\(m\) symbol.
The inclusion \(K_\Sigma F\subset F_\Sigma^{\rm el}\) shows
that a prescribed oriented label fixes the rational support label,
so the same probe \(r\) works for every moving prime with that label.
At \(p\in U'\), apply the level-prime identity of
Lemma~\ref{lem:cubic-deletion} in the saturated lattice of
\(\widetilde f_{U'}\). Its \(U_p\)-stability follows from the
same oldform decomposition; the trace at the fresh prime \(r\)
commutes with \(U_p\). Taking coefficient \(1\), hence coefficient
\(p\) on the \(U_p\)-side, tests
\[
 \eta(p)\nu_p(D')\,2^{-\delta-|U'|}a(D'),
 \qquad D'=(U'\setminus\{p\})q'r,
\]
with \(q'\) split. The square-free source has split weight
\(m-1\), so minimality gives
\[
 \overline{2^{-\delta-(m-1)}a(D')}=0.
\]
Since \(\nu_p(D')\in\Z\), the displayed coefficient has zero
reduction. Hence \(\Phi_r\) is invariant under inertia at \(p\).
It is invariant under \(G_\Q\)-conjugation because \(T\) is a
class function; the same invariance holds at every place above
\(p\). Removing these inertia groups leaves
\(F_\Sigma^{\rm el}\).

\emph{Replacing primes.} Suppose two moving split primes have
matched oriented labels and avoid a fixed set \(U'\).
Choose one inert probe fresh for their union and \(U'\), with
the support label restoring admissibility. The first paragraph
replaces either configuration's original inert set by this probe.
The equality of \(\Phi_r\) on matched labels then replaces the
moving prime. Every intermediate configuration is admissible of
split weight \(m\).

For two realizations, choose a third split array disjoint from both,
with the same types. Compare the first array to the third one by
successive replacements, then the second to the third. At each step
the probe is fresh for the finite union and is common to the two
moving primes. This proves invariance of the symbol. Changing an orientation
conjugates its label, so matching the multisets of type orbits gives
\eqref{eq:cubic-type-transfer}.
Separately, Lemma~\ref{lem:ray-realization} realizes the matched
types and prescribed cross residues. The first paragraph permits
arbitrary unramified twisting bits at split weight \(m\).
\end{proof}

\subsection{Symmetry and the Selmer bound}

For types \(\mathfrak h,\mathfrak h'\), let
\(b(\mathfrak h,\mathfrak h')\) be the sum of the local cup
pairings of the corresponding classes at the support. At an active split prime, decompose each class into its Kummer
lift with the same inertia coordinate and an unramified difference.
The Kummer plane and the unramified plane are both isotropic, so
only the cross pairings remain. Global reciprocity gives
\begin{equation}\label{eq:cubic-reciprocity-defect}
 B_{qp}+B_{pq}^{t}=b(\mathfrak h_q,\mathfrak h_p),
\end{equation}
also on the diagonal: the same decomposition includes the inertia self-pairing in the isotropic Kummer term. The off-diagonal formula puts \(b\) in
\(J\F_2[\mathsf R]\); on equal types it is alternating. Choose matrices
\(L(\mathfrak h,\mathfrak h')\) in this same span with
\[
 L(\mathfrak h,\mathfrak h')+
 L(\mathfrak h',\mathfrak h)^t=b(\mathfrak h,\mathfrak h').
\]
For unequal types either direction can be chosen first. For equal
types use \(J\mathsf R+(J\mathsf R)^t=J\). Take \(L=0\) whenever one type is
primary, since its support pairings vanish. Then
\begin{equation}\label{eq:cubic-symmetric-matrix}
 M_{qp}=B_{qp}+L(\mathfrak h_q,\mathfrak h_p)
\end{equation}
is symmetric. Its off-diagonal blocks can be prescribed arbitrarily
in \(J\F_2[\mathsf R]\), subject to transposition in the reverse direction.
Fix a symmetric representative \(Q_{\mathfrak h}\) of each
type's self block modulo \(J\). Denote the primary type by \(0\)
and choose \(Q_0=0\). The corrections and representatives are
chosen once for this finite type set.

\begin{proposition}\label{prop:cubic-nullity}
For any configuration using a fixed finite set \(\mathfrak T\)
of types, and with admissible twisting parameter,
\[
 \dim_{\F_2}\Sel_2(E^h/\Q)\ge\dim\ker M-2|\mathfrak T|.
\]
If \(M\) has zero entries between prescribed blocks, its nullity
is the sum of their nullities. When the symbol is one the number of
singular blocks is consequently bounded independently of the array.
\end{proposition}
\begin{proof}
Intersect \(\ker M\) with the kernel of
\[
 \bigoplus_pV\longrightarrow\bigoplus_{\mathfrak h\in\mathfrak T}V,
 \qquad (x_p)\longmapsto
 \left(\sum_{p:\,\mathfrak h_p=\mathfrak h}x_p\right)_{\mathfrak h}.
\]
This has codimension at most \(2|\mathfrak T|\).
These type sums kill the \(L\)-correction in every row and kill the local
classes at \(\Sigma\), since those maps depend only on the type.
Thus \(\sum_pC_p(x_p)\) satisfies the active split conditions,
the support conditions, and the unramified conditions at all other
good primes. The local cohomology at each inert twisting prime is
zero. The sum is a Selmer class, and its independent inertia
coordinates make this map injective. This proves the inequality.
When the symbol is one, Lemma~\ref{lem:selmerbound} gives
\[
 \dim\ker M\le M_E+2|\mathfrak T|.
\]
For a block diagonal matrix every singular block contributes at
least one to its nullity; hence a value-one configuration has at
most \(M_E+2|\mathfrak T|\) singular blocks.
\end{proof}

\subsection{Three states and simultaneous calibration}

A site is one split prime. Its three states are absence, presence
with self block \(Q_{\mathfrak h}\), and presence with self block
\(Q_{\mathfrak h}+J\). For a finite array and a chosen state at
each site, choose an aggregate inert prime restoring the total
fixed-support label and the specified self blocks. These are
quadratic conditions, compatible with inertness because
\(F\cap K_U=\Q\). Proposition~\ref{prop:cubic-W0} makes the
coefficient-symbol value independent of that choice. The state with
every site absent has value \(A(\varnothing)\).

The three states are encoded in an auxiliary binary space by
\[
 \begin{aligned}
 \text{presence with }Q_{\mathfrak h}&\longleftrightarrow \mathbf e_0,\\
 \text{presence with }Q_{\mathfrak h}+J&\longleftrightarrow \mathbf e_1,\\
 \text{absence}&\longleftrightarrow \mathbf e_0+\mathbf e_1.
 \end{aligned}
\]
These vectors encode evaluations and are independent of \(V=E[2]\).
A \emph{block} is a set of sites with fixed pair matrices.
A \emph{state} on an array is a finitely supported assignment of
one of the three choices to every site; its support consists of the
present sites. Its evaluated matrix is the corresponding finite
principal matrix, with the chosen self blocks. A block state is
\emph{singular} if this matrix is singular. A \emph{collection of
singular states} is a state whose restriction to each block meeting
its support is singular.

\begin{lemma}\label{lem:cubic-three-state}
For a symmetric binary matrix, toggling the two off-diagonal entries
within one coordinate pair gives
\[
 \det M+\det(M+J_p)=\det M^{\widehat p},
\]
where \(M^{\widehat p}\) deletes that pair of coordinates. Moreover,
in a tensor product of two-dimensional binary spaces, the non-zero
pure tensors annihilated by a non-zero functional \(f\) span
\(\ker f\).
\end{lemma}
\begin{proof}
In the determinant expansion, a permutation and its inverse cancel
unless they agree. Only products of diagonal entries and squares
of disjoint off-diagonal entries remain. Changing the chosen entry
by one therefore changes the determinant by the complementary
principal determinant.

For the tensor assertion, suppose a functional \(g\) vanishes on
the pure zeros of \(f\). Slice in the last factor at its three
non-zero vectors, obtaining \(f_3=f_1+f_2\) and \(g_3=g_1+g_2\).
Induction gives \(g_i=\epsilon_i f_i\) when \(f_i\ne0\), and
\(g_i=0\) when \(f_i=0\). If all three slices are non-zero, the
first two are independent and the scalars agree. If one is zero,
the other two agree and again have one common scalar. Thus
\(g=0\) or \(g=f\), which proves the span assertion by duality.
\end{proof}

\begin{lemma}\label{lem:cubic-state-relation}
With the other states fixed, the three symbol values at a site sum
to zero. The same identity holds for the determinant of its
evaluated matrix, with empty determinant one.
\end{lemma}
\begin{proof}
Changing the aggregate-inert bit at \(p\) adds \(J\) to its
self block. Lemma~\ref{lem:cubic-deletion} identifies the sum of
these two symbol values with the deleted-site value. Its transferred
label keeps all other self blocks unchanged. The determinant identity is
Lemma~\ref{lem:cubic-three-state}.
\end{proof}

The three-state relations make both the symbol and
the determinant linear on the tensor product of these spaces.
For a finite set of sites in a block, let \(d\) be its determinant
functional and
\(v_\varnothing\) its all-absent state, so
\(d(v_\varnothing)=1\). If \(d(v)=1\), then
\(v+v_\varnothing\in\ker d\), and
Lemma~\ref{lem:cubic-three-state} gives
\[
 v=v_\varnothing+\sum_i v_i,\qquad d(v_i)=0,
\]
with singular pure states \(v_i\). This identity can be evaluated
by the symbol with all other blocks fixed.

The correction to each self block is linear in the family index.
One inert prime for each coordinate therefore realizes the prescribed
state values in a single quadratic twist family.

\begin{proposition}\label{prop:cubic-calibration}
Fix a finite array of split sites with prescribed pair blocks.
Give each site linear forms \(\alpha_p,\beta_p\) on \(X_b=\F_2^b\).
At \(x\), omit the site if \(\alpha_p(x)=0\), and otherwise give
it self block \(Q_{\mathfrak h_p}+\beta_p(x)J\). There is one
admissible quadratic twist family realizing all these evaluations. It
uses the split primes and at most \(b\) further inert primes,
all with linear ramification forms and all absent at zero.
\end{proposition}
\begin{proof}
After the split primes and generators have been chosen, the
self-block formula and the fixed symmetrizing correction give
\[
 M_{pp}=Q_{\mathfrak h_p}+(s_p+\kappa_p(h))J.
\]
The type fixes \(Q_{\mathfrak h_p}\); the observed self residues
fix \(s_p\in\F_2\). The inert twisting bit then selects the
remaining \(J\)-offset. Put
\(n_{pq}=[q^*/p]\in\F_2\), where brackets denote the additive
Legendre symbol. The required contribution from the inert primes
at \(p\) is the linear form
\begin{equation}\label{eq:cubic-calibration}
 r_p(x)=s_p\alpha_p(x)+\beta_p(x)
       +\sum_{q\ne p}n_{pq}\alpha_q(x).
\end{equation}
When \(p\) is active, this cancels its actual offset and the unit
twist from all the other split primes, leaving \(\beta_p(x)J\).
When \(p\) is absent its prescribed bit has no effect.

Let \(e_j\) be the basis vector dual to the coordinate \(x_j\).
For each \(j\), choose one fresh inert prime \(r_j\).
Prescribe its \(K_\Sigma\)-label to cancel
\(\sum_p\alpha_p(e_j)\lambda_\Sigma(p)\), and prescribe
\([r_j^*/p]=r_p(e_j)\) for every site. These specifications
give an element of the rational multiquadratic Galois group:
the fresh prime radical coordinates are independent of the
support radical coordinates. This multiquadratic field is disjoint
from \(F\), so Chebotarev permits the independent order-three
Frobenius condition. The primes can be chosen distinct and away
from the split array.

Activate \(r_j\) by \(x_j\). The total fixed-support label is
zero for every \(x\); signed prime discriminants then give
positivity and local squareness at \(\Sigma\). Formula
\eqref{eq:cubic-calibration} proves the self-block assertion.
At any \(x\), the active inert primes have exactly the aggregate
label defining the corresponding state value, so
Proposition~\ref{prop:cubic-W0} identifies the two symbols.
Primes with identically zero ramification form may be omitted. At zero
the inert product is empty, as required.
\end{proof}

\subsection{Selection on disjoint blocks}

Between distinct blocks we prescribe \(M_{qp}=0\), or equivalently
\(B_{qp}=L(\mathfrak h_q,\mathfrak h_p)\). Proposition~\ref{prop:cubic-matrix-realization}
permits these prescriptions. The value of a finitely supported state
is its coefficient symbol, with the aggregate inert label specified
above.

\begin{lemma}\label{lem:cubic-selector}
Let \(M_*\) be a singular block state. For \(b\ge2\) and
\(a\ne0\) in \(\F_2^b\), it can be extended by primary sites
with linear ramification forms so that its matrix is singular at
\(a\), invertible at every other non-zero assignment, and empty at
zero. At \(a\), precisely the original state remains.
\end{lemma}
\begin{proof}
By Lemma~\ref{lem:cubic-three-state}, the determinant as a function
of the within-site toggle bits is multilinear, with top coefficient
equal to the empty determinant, one. Thus there is a sum \(D\)
of within-site \(J\)-toggles for which \(M_*+D\) is invertible.

Put \(m=b-1\). Choose a linear form \(l\) with \(l(a)=1\),
and a basis \(y_1,\ldots,y_m\) of the forms vanishing at \(a\).
Activate every site of the original state by \(l\), retaining its
chosen self block. We construct a component that adds one prescribed
\(J\)-toggle whenever \(y\ne0\), and is absent when \(y=0\).

For this component, take primary sites \(u_i,v_i\), both activated
by \(y_i\), for \(1\le i\le m\). In the order
\(u_1,\ldots,u_m,v_1,\ldots,v_m\), prescribe their matrix to be
\[
 P\otimes J,\qquad
 P=\begin{pmatrix}0&T\\T^t&I_m\end{pmatrix},\qquad
 T_{ij}=\begin{cases}1,&i\ge j,\\0,&i<j.\end{cases}
\]
Couple the chosen old site to every \(u_i\) by \(J\), and to
every \(v_i\) by zero. The self blocks are zero at the \(u_i\)
and \(J\) at the \(v_i\); all these choices are permitted by
Proposition~\ref{prop:cubic-matrix-realization} and the simultaneous
calibration of Proposition~\ref{prop:cubic-calibration}.

At an assignment \(x\), put \(I=\{i:y_i(x)=1\}\). If \(I\ne\varnothing\), the induced
matrix \(T_I\), with indices in increasing order, is again lower
triangular with every entry on and below the diagonal equal to one.
Consequently \(P_I\) is invertible and
\[
 (P_I^{-1})_{uu}=T_I^{-t}T_I^{-1},\qquad
 T_I^{-1}\mathbf1_I=(1,0,\ldots,0)^t.
\]
Since \(J^2=I\), the Schur contribution to the old site is
\[
 \bigl(\mathbf1_I^t T_I^{-t}T_I^{-1}\mathbf1_I\bigr)J
 =J.
\]
If \(I=\varnothing\), the component is absent and contributes
zero. Use one disjoint component for each toggle occurring in
\(D\), with zero couplings between components.

At \(x=a\), all \(y_i\) vanish and \(l=1\), so precisely
the original state remains. At any other non-zero \(x\) with
\(l(x)=1\), one has \(y(x)\ne0\), and the Schur complement
on the old sites is \(M_*+D\); both it and the auxiliary matrix
are invertible. If \(l(x)=0\), only invertible auxiliary
components remain. At zero every site is absent.
The construction uses \(2(b-1)\) primary sites per toggle in \(D\).
\end{proof}

\begin{lemma}\label{lem:cubic-disjoint}
Fix \(b\ge2\) and put \(K=2^b-1\). There are \(K\) disjoint
collections \(\mathcal B_1,\dots,\mathcal B_K\), each of \(r\ge1\)
blocks, and a collection \(W_j\) of one singular state on each
block of \(\mathcal B_j\), such that
each \(W_j\) has value one and every collection of singular states
on more than \(r\) blocks of \(\bigcup_j\mathcal B_j\) has value
zero. Their available sites include the selectors of
Lemma~\ref{lem:cubic-selector} for every finitely supported singular
block state.
\end{lemma}
\begin{proof}
Choose once, independently of \(b\), the types of a least-weight
configuration with symbol one, together with the primary type, and put
\(B_E=M_E+2|\mathfrak T|\). Construct more than \(KB_E^2\) blocks,
each containing a fresh copy of this configuration. Let
\(\mathcal V_0\) be their finite site array, partitioned into these
blocks. Inductively form \(\mathcal V_{n+1}\) by adjoining, within
the same block, fresh primary sites for the selectors of every
singular block state supported in \(\mathcal V_n\), for every
non-zero target in \(\F_2^b\). Set all undesignated \(M\)-couplings
to zero, including those between blocks. At each stage there are
finitely many requests and residue conditions, so
Proposition~\ref{prop:cubic-matrix-realization} extends the actual
prime array. The countable union
\(\mathcal V_\infty=\bigcup_{n\ge0}\mathcal V_n\) is closed
under every finite selector request: each finitely supported state
already occurs at some stage. All maxima below are taken in this
fixed array. Every initial copy has value one by
\eqref{eq:cubic-type-transfer}; the final family uses only finitely
many primes of \(\mathcal V_\infty\).

For a nonempty set \(\mathcal S\) of these blocks, let
\(r_{\mathcal S}\) be the supremum of the numbers of blocks in a
collection consisting entirely of singular states, of value one,
supported on \(\mathcal S\).
Equation~\eqref{eq:cubic-selmer-bound} and
Proposition~\ref{prop:cubic-nullity} bound it by \(B_E\).
It is positive if a seed state is singular. Otherwise expand that
state as empty plus singular states, using
Lemmas~\ref{lem:cubic-state-relation} and~\ref{lem:cubic-three-state}.
The empty configuration has value zero, so some nonempty singular
state has value one. The possible sizes therefore form a nonempty
subset of \(\{1,\dots,B_E\}\), so their supremum is a maximum
attained by a finite configuration.

Starting with all blocks, choose a collection attaining
\(r_{\mathcal S}\), remove its blocks, and repeat \(KB_E\) times.
At most \(B_E\) blocks are removed at each step, so the remaining
set is always nonempty. The nonincreasing integers
\(r_{\mathcal S}\) lie between one and \(B_E\); some value
\(r\) occurs at least \(K\) times. The corresponding maximizing
collections are disjoint and lie in the first set with \(r_{\mathcal S}=r\).
They satisfy the assertion. Maximality ranges over all finitely
supported states in the sites already constructed, including every
later selector state.
\end{proof}

\begin{proof}[Proof of Theorem~\ref{thm:cubic-E}]
Number the non-zero indices as \(a_1,\dots,a_K\), and use
Lemma~\ref{lem:cubic-disjoint}. Address every block of \(W_j\)
to \(a_j\) by Lemma~\ref{lem:cubic-selector}, within its available
sites. At \(a_j\), precisely \(W_j\) remains in its designated
blocks and every other block is invertible. Expand those other
states as empty plus singular states. Choosing empty states leaves
\(W_j\), of value one; every other term has more than \(r\)
singular blocks and has value zero by the lemma. Thus every non-zero
index has symbol one, and every site is absent at zero.

Only finitely many sites have been used. Each has a linear
ramification form, and its present self block is fixed; its diagonal
choice is therefore a constant times that form.
Proposition~\ref{prop:cubic-calibration} realizes all these
evaluations in one admissible quadratic twist family. Restriction to a
line gives the assertion for \(b=1\).
\end{proof}

This proves Theorem~\ref{thm:families} in case~\CaseIII{}.
The deduction in the introduction completes the proof of
Theorem~\ref{thm:A}.

\appendix
\section{Further residual criteria}\label{sec:residual-appendix}

Proposition~\ref{prop:RF-general} and
Corollary~\ref{cor:RF-stability} give class-group criteria for
\RF{} beyond the abelian two-division case.\footnote{For the
relation between class-group and fine-Selmer \(\mu\)-invariants at
odd primes, see \cite{CS05}. The criterion here retains the
real-place condition at \(2\).} They extend the
uniform interpolation theorem to further irreducible residual
representations. In particular, Proposition~\ref{prop:RF-general}(2)
applies Theorem~\ref{thm:main} unconditionally to
\(y^2=x^3-x-1\), with \(S_3\) two-division field.
The family constructions of Sections~\ref{sec:rules}--\ref{sec:cyclic-cubic}
remain restricted to cases~\CaseI{}--\CaseIII{}.

\begin{proposition}\label{prop:RF-general}
Suppose \(E[2]\) is irreducible. For a non-zero \(P\in E[2]\), put
\(k=\Q(x(P))\) and \(L=k(i)\).
\begin{enumerate}
\item If the cyclotomic \(\Z_2\)-extension of \(L\) has class-group
\(\mu\)-invariant zero, then \(E\) satisfies \RF{}.
\item The same conclusion holds if \(L\) has odd class number and
a single prime above \(2\). This applies to
\(E:y^2=x^3-x-1\), whose two-division field has group \(S_3\).
\end{enumerate}
\end{proposition}
\begin{proof}
For (1), the permutation representation on the three non-zero points
of \(E[2]\) decomposes as
\[
 \F_2[G_\Q/G_k]\simeq\F_2\oplus E[2].
\]
Indeed, the augmentation splits because the degree is odd, and its
two-dimensional kernel is the natural representation. Put
\(k_n=k\Q_n\) and \(K_n=k_n(i)=k(\mu_{2^{n+2}})\).
Shapiro's lemma reduces the degree-two bound to constant coefficients
over \(k_n\). If \(k\) has complex places, the extra archimedean terms
in this comparison are concentrated in degree one and do not change
degree two.

Put \(S_{k_n}=\{v:v\text{ lies above }S\}\).
Lemma~\ref{lem:class-bound} gives
\(d_2\Cl^+(k_n)\le d_2\Cl(K_n)+1\).

The assumption \(\mu=0\) makes the inverse-limit class group of
\(K_n\) finitely generated over \(\Z_2\). After a finite initial
segment, a dyadic prime is totally ramified in every successive
layer. The class-group norm maps are then surjective: the upper
layer has trivial intersection with the lower layer's Hilbert class
field. Hence the \(2\)-ranks \(d_2\Cl(K_n)\) are bounded; see
\cite[Chapter~13]{Washington} for the Iwasawa module structure.

Lemma~\ref{lem:class-bound} now bounds
\(\dim H^2_+(k_n,\F_2)\), since the number of places above the
fixed support \(S\) is bounded. This proves (1).

For (2), put \(\widetilde k_n=k(\mu_{2^{n+2}})\). Suppose a layer has odd class number and a unique dyadic prime.
The next quadratic layer is unramified outside \(2\), and at
infinity because the fields are totally imaginary. If it were also
unramified at the dyadic prime, it would contradict odd class number.
Thus that prime ramifies and remains unique. The
ambiguous class number formula shows that the invariant class
group at the next step has odd order if the present class group
does. A non-zero finite \(2\)-group with an involution has a non-zero
fixed element, so induction proves odd class number at every step.
The inverse-limit \(2\)-class group is zero, and (1) applies.

For the stated example, the polynomial \(x^3-x-1\) has discriminant
\(-23\), so the cubic field \(k\) has one real place and one pair
of complex places, and its normal closure has group \(S_3\).
The Minkowski bound \(8\sqrt{23}/(9\pi)<2\) gives class number one.
The reduction \(x^3+x+1\) is irreducible over \(\F_2\), so \(2\)
is inert in \(k\); adjoining \(i\) ramifies its unique dyadic prime.
In \(k(i)/k\), exactly that finite prime and the real place ramify.
The ambiguous class number formula gives
\[
 \#\Cl(k(i))^{\Gal(k(i)/k)}
 =\frac{2}{[\mathcal O_k^\times:
                \mathcal O_k^\times\cap N_{k(i)/k}(k(i)^\times)]}.
\]
The denominator is at least two, since \(-1\) is not a norm at the
real place. The right side is a positive integer, hence is one.
Again the involution criterion shows that \(k(i)\) has odd class
number.

\end{proof}

\begin{corollary}\label{cor:RF-stability}
Suppose \(E[2]\) is irreducible, let \(k=\Q(x(P))\) for a non-zero
\(P\in E[2]\), and put \(L=k(i)\). Let \(L_n\) be the degree
\(2^n\) layer of the cyclotomic \(\Z_2\)-extension \(L_\infty/L\).
Choose \(n_0\) such that every prime ramified in this tower is
totally ramified in \(L_\infty/L_{n_0}\). If, for some \(n\ge n_0\),
\[
 d_2\Cl(L_n)=d_2\Cl(L_{n+1}),
\]
then \(E\) satisfies \RF{}.
\end{corollary}
\begin{proof}
Fukuda's rank-stability theorem \cite[Theorem~1(2)]{Fukuda94}
gives \(d_2\Cl(L_m)=d_2\Cl(L_n)\) for all \(m\ge n\).
Thus the class-group \(\mu\)-invariant of \(L_\infty/L\) is zero,
and Proposition~\ref{prop:RF-general}(1) applies.
\end{proof}

\section{The congruent-number curves}\label{sec:congruent-example}

We illustrate the families and the dimension bound for
\[
 E_d:\ y^2=x^3-d^2x,
\]
with \(d>0\) square-free. These curves have full rational
two-torsion, so case~\CaseI{} applies. Take \(\Sigma\) to be the
primes dividing \(2d\). A square-free parameter \(h\) is admissible
precisely when
\[
 h>0,\qquad (h,2d)=1,\qquad h\equiv1\pmod8,\qquad
 \Bigl(\frac h\ell\Bigr)=1\quad(\ell\mid d,\ \ell\text{ odd}).
\]
The twist is \(E_{dh}\), and every prime dividing \(h\) contributes
a Tamagawa factor \(4\). Thus \(c(h)=4^{\omega(h)}\), where
\(\omega(h)\) counts the prime divisors.

\subsection{The curve \texorpdfstring{\(E_1\)}{E1}}
Tunnell's form is
\[
 g_{\mathrm T}(q)=
 q\prod_{m\ge1}(1-q^{8m})(1-q^{16m})
 \sum_{r\in\Z}q^{2r^2}
 =\sum_{n\ge1}a_{\mathrm T}(n)q^n.
\]
In our period convention, for odd square-free \(n>0\),
\begin{equation}\label{eq:tunnell-normalization}
 \Omega^+_{E_n}=\frac1{\sqrt n}
 \int_1^\infty\frac{dx}{\sqrt{x^3-x}},\qquad
 \Lalg(E_n,1)=\frac{a_{\mathrm T}(n)^2}{4}
 \quad\text{\cite{Tunnell}}.
\end{equation}
For the admissible support \(n\equiv1\pmod8\), the form of
Theorem~\ref{thm:D} may therefore be chosen with
\(n_0=1\), \(\psi=1\), \(a(n)=a_{\mathrm T}(n)\) and \(c_g=2\).

Burungale--Flach prove the full BSD formula for CM elliptic curves
over \(\Q\) of analytic rank zero, including its \(2\)-part
\cite[Corollary~2]{BF24}. Since \(c_2(E_n)=2\) for
\(n\equiv1\pmod8\), and \(E_n(\Q)_{\mathrm{tors}}=E_n[2]\), their
formula and \eqref{eq:tunnell-normalization} give, when
\(a_{\mathrm T}(n)\ne0\),
\begin{equation}\label{eq:congruent-minimum}
 \val(a_{\mathrm T}(n))-\omega(n)
 =\tfrac12\val\bigl(\#\Sha(E_n/\Q)[2^\infty]\bigr)\ge0.
\end{equation}
As \(a_{\mathrm T}(1)=1\), the global minimum is \(\delta=0\).
The same calculation gives \(D(n)=-2\) for every admissible
\(n\) with non-zero central value.

A family on \(X_2\) is obtained from \(P=\{3,11,19\}\) and
\[
 \alpha_3=x_1,\qquad \alpha_{11}=x_2,\qquad
 \alpha_{19}=x_1+x_2.
\]
All three primes have label \((1,1)\) on \((-1,2)\), and the
three forms sum to zero. Every non-zero index activates two
primes congruent to \(3\) modulo \(8\); their signed discriminants
have positive product congruent to \(1\) modulo \(8\).
The coefficients and their normalization are
\[
\begin{array}{c@{\qquad}r@{\qquad}r@{\qquad}r@{\qquad}c}
 x&h_x&c(h_x)&a_{\mathrm T}(h_x)&
 \val(a_{\mathrm T}(h_x))-\omega(h_x)\\[2pt]\hline
 00&1&1&1&0\\
 10&57&16&4&0\\
 01&209&16&4&0\\
 11&33&16&4&0
\end{array}
\]
Consequently \(\Lalg(E_{h_x},1)/c(h_x)=1/4\) at all four
indices.\footnote{For odd \(n\), put
\(R_j(n)=\#\{(u,v,w)\in\Z^3:2u^2+v^2+jw^2=n\}\).
Then \(2a_{\mathrm T}(n)=2R_{32}(n)-R_8(n)\).
At \(n=1,57,209,33\), the pairs \((R_{32},R_8)\) are
\((2,2),(12,16),(44,80),(12,16)\), respectively.}

There is also a short Selmer calculation. For \(n=pq\), with
\(p,q\equiv3\pmod8\), work over \(\F_2\) and put
\[
 \epsilon=\Bigl[\Bigl(\frac qp\Bigr)=-1\Bigr],\qquad
 B=\begin{pmatrix}\epsilon&\epsilon\\1+\epsilon&1+\epsilon\end{pmatrix}.
\]
Monsky's matrix for the pure \(2\)-Selmer group is
\(\left(\begin{smallmatrix}B+I&I\\I&B\end{smallmatrix}\right)\)
\cite[Appendix]{HB94}. Since \(B^2=B\), its kernel is zero:
the first equation gives \(y=(B+I)x\), and the second then gives
\(x=0\). Thus \(\Sel_2(E_n/\Q)=E_n[2]\) and
\(\Sha(E_n/\Q)[2^\infty]=0\) at the three auxiliary members.

\subsection{A vanishing base value}
For \(E_{41}\), one has \(a_{\mathrm T}(41)=0\) and
\(\dim_{\F_2}\Sel_2(E_{41}/\Q)=4\).
Its \(2^\infty\)-Selmer group is infinite.\footnote{The point
\((1025/16,25215/64)\) on \(E_{41}\) is non-torsion by the
Lutz--Nagell theorem.}
For admissible \(h\) with \(a_{\mathrm T}(41h)\ne0\), the same
BSD calculation gives
\[
 \val(a_{\mathrm T}(41h))-\omega(h)
 =1+\tfrac12\val\bigl(\#\Sha(E_{41h}/\Q)[2^\infty]\bigr)\ge1.
\]
\Needspace{10\baselineskip}
This minimum is attained simultaneously in the family
\[
 \alpha_3=\alpha_{17}=x_1,\qquad
 \alpha_7=\alpha_{13}=x_2,\qquad
 \alpha_{43}=x_1+x_2:
\]
\[
\begin{array}{c@{\qquad}r@{\qquad}r@{\qquad}c}
 x&h_x&a_{\mathrm T}(41h_x)&
 \val(a_{\mathrm T}(41h_x))-\omega(h_x)\\[2pt]\hline
 10&2193&16&1\\
 01&3913&16&1\\
 11&4641&-32&1
\end{array}
\]
Each \(h_x\) is \(1\) modulo \(8\) and a quadratic residue
modulo \(41\). Monsky's matrices at \(41h_x\) have full rank,
so \(\Sel_2(E_{41h_x}/\Q)=E_{41h_x}[2]\).
Normalizing by \(a(h)=a_{\mathrm T}(41h)/16\), with
\(n_0=2193\), gives \(a(1)=0\) and \(\delta=-3\).

This family illustrates the construction at a vanishing base
value. Passing to the large families of Theorem~\ref{thm:E}
gives \(\mathscr L_{\mathcal F,0}(0)\ne0\), as in
Remark~\ref{rem:trivial-character}. Here the base has positive
Selmer corank and vanishing central \(L\)-value.

\subsection{The sufficient dimension}
Fix \(E_d\), its Beilinson--Kato normalization and the form of
Theorem~\ref{thm:D}. Put
\[
 H_d=2\delta-c_g+\beta_{E_d},\qquad
 b_d=2^{C_{E_d}+H_d}.
\]
Suppose \(s_2(E_d)=0\) and, for contradiction, \(L(E_d,1)=0\).
Monsky's parity theorem gives root number \(+1\).
Choose a seed with symbol one, of \(r\) primes. It is non-empty,
and the one-dimensional comparison gives \(C_{E_d}+H_d\ge0\).
For each non-zero form on \(\F_2^{b_d}\), take a copy of its
\(r\) labels. This gives \(r(2^{b_d}-1)\) terminals, with
\(r2^{b_d-1}\) active at each non-zero index and total label zero.
The construction of Section~\ref{sec:seeds} and the selectors of
Theorem~\ref{thm:address} give an admissible family attaining
\(\delta\) at every non-zero index. Equation~\eqref{eq:mechanism}
and Theorem~\ref{thm:B} then give \(L(E_d,1)\ne0\).
The seed size affects the number of primes, but not this dimension.

For \(E_1\), the fixed quantities are
\[
 a_{E_1}=2,\qquad \rho_{E_1}=2,\qquad \epsilon_\infty=1.
\]
Indeed \(E_1(\Q_2)[2^\infty]=E_1[2]\), since \(-1\) and \(2\)
are non-squares in \(\Q_2\). Hence
\[
 \beta_{E_1}=\alpha_{E_1}-\lambda_{E_1}+2,\qquad
 H_1=\alpha_{E_1}-\lambda_{E_1},\qquad
 b_1=2^{c_{\mathrm{int}}+\alpha_{E_1}-\lambda_{E_1}}.
\]
Here we take \(C_{E_1}=c_{\mathrm{int}}\), as in the proof of
Theorem~\ref{thm:main}. Only the lattice index and the local
reciprocity constants remain to be evaluated numerically.

\end{document}